# VEDIC MATHEMATICS -
# 'VEDIC' OR 'MATHEMATICS':
# A FUZZY & NEUTROSOPHIC ANALYSIS

**W. B. VASANTHA KANDASAMY**
**FLORENTIN SMARANDACHE**

**2006**

# VEDIC MATHEMATICS -
# 'VEDIC' OR 'MATHEMATICS':
# A FUZZY & NEUTROSOPHIC ANALYSIS

**W. B. VASANTHA KANDASAMY**
e-mail: **vasanthakandasamy@gmail.com**
web: **http://mat.iitm.ac.in/~wbv**
**www.vasantha.net**

**FLORENTIN SMARANDACHE**
e-mail: **smarand@unm.edu**

**2006**

# CONTENTS

# PREFACE

Religious extremism has been the root cause of most of the world problems since time immemorial. It has decided the fates of men and nations. In a vast nation like India, the imposition of religious dogma and discrimination upon the people has taken place after the upsurge of Hindu rightwing forces in the political arena. As a consequence of their political ascendancy in the northern states of India, they started to rewrite school textbooks in an extremely biased manner that was fundamentalist and revivalist. Not only did they meddle with subjects like history (which was their main area of operation), but they also imposed their religious agenda on the science subjects. There was a plan to introduce Vedic Astrology in the school syllabus across the nation, which was dropped after a major hue and cry from secular intellectuals.

This obsession with 'Vedic' results from the fundamentalist Hindu organizations need to claim their identity as Aryan (and hence of Caucasian origin) and hence superior to the rest of the native inhabitants of India. The 'Vedas' are considered 'divine' in origin and are assumed to be direct revelations from God. The whole corpus of Vedic literature is in Sanskrit. The Vedas are four in number: Rgveda, Saamaveda, Yajurveda and Atharvaveda. In traditional Hinduism, the Vedas as a body of knowledge were to be learnt only by the 'upper' caste Hindus and the 'lower castes' (Sudras) and so-called 'untouchables' (who were outside the Hindu social order) were forbidden from learning or even hearing to their recitation. For several centuries, the Vedas were not written down but passed from generation to generation through oral transmission. While religious significance is essential for maintaining Aryan supremacy and the caste system, the claims made about the Vedas were of the highest order of hyperbole. Murli Manohar Joshi, a senior Cabinet minister of the Bharatiya Janata Party (BJP) that ruled India from 1999-2004 went on to claim that a cure of the dreaded AIDS was available in the Vedas! In the continuing trend, last week a scientist has announced that NASA (of the USA) is using a Vedic formula to produce electricity. One such popular topic of Hindutva imposition was Vedic Mathematics. Much of the hype about this topic is based on one single book authored by the Sankaracharya (the highest Hindu pontiff) Jagadguru Swami Sri Bharati Krsna Tirthaji Maharaja titled Vedic Mathematics and published in the year 1965, and reprinted several times since the 1990s [51]. This book was used as the foundation and the subject was systematically introduced in schools across India. It was introduced in the official curriculum in the school syllabus in the states of Uttar Pradesh and Madhya Pradesh. Further, schools run by Hindutva sympathizers or trusts introduced it into their curriculum. In this juncture, the first author of this book started working on this topic five years back, and has since met over 1000 persons from various walks of life and collected their opinion on Vedic Mathematics. This book is the result of those interactions.

In this book the authors probe into Vedic Mathematics (a concept that gained renown in the period of the religious fanatic and revivalist Hindutva rule in India): and explore whether it is really 'Vedic' in origin or 'Mathematics' in content. The entire field of Vedic Mathematics is supposedly based on 16 one-to-three-word sutras (aphorisms) in Sanskrit, which they claim can solve all modern mathematical problems. However, a careful perusal of the General Editor's note in this book gives away the basic fact that the origin of these sutras are not 'Vedic' at all. The book's General Editor, V.S. Agrawala, (M.A., PhD. D.Litt.,) writes in page VI as follows:

> "It is the whole essence of his assessment of Vedic tradition that it is not to be approached from a factual standpoint but from the ideal standpoint viz., as the Vedas, as traditionally accepted in India as the repository of all knowledge, *should be* and not what they are in human possession. That approach entirely turns the table on all critics, for the authorship of Vedic mathematics need not be labouriously searched for in the texts as preserved from antiquity. […]

In the light of the above definition and approach must be understood the author's statement that the sixteen sutras on which the present volume is based from part of a Parisista of the Atharvaveda. We are aware that each Veda has its subsidiary apocryphal text some of which remain in manuscripts and others have been printed but that formulation has not closed. For example, some Parisista of the Atharvaveda were edited by G.M.Bolling and J. Von Negelein, Leipzig,1909-10. But this work of Sri Sankaracharyaji deserves to be regarded as a new Parisista by itself and it is not surprising that the Sutras mentioned herein do not appear in the hitherto known Parisistas.

A list of these main 16 Sutras and of their sub-sutras or corollaries is prefixed in the beginning of the text and the style of language also points to their discovery by Sri Swamiji himself. At any rate, it is needless to dwell longer on this point of origin since the vast merit of these rules should be a matter of discovery for each intelligent reader. Whatever is written here by the author stands on its own merits and is presented as such to the mathematical world. [emphasis supplied]"

The argument that Vedas means all knowledge and hence the fallacy of claiming even $20^{th}$ century inventions to belong to the Vedas clearly reveals that there is a hidden agenda in bestowing such an antiquity upon a subject of such a recent origin. There is an open admission that these sutras are the product of one man's imagination. Now it has become clear to us that the so-called Vedic Mathematics is not even Vedic in origin.

Next, we wanted to analyze the mathematical content and its ulterior motives using fuzzy analysis. We analyzed this problem using fuzzy models like Fuzzy Cognitive Maps (FCM), Fuzzy Relational Maps (FRM) and the newly constructed fuzzy dynamical system (and its Neutrosophic analogue) that can analyze multi-experts opinion at a time using a single model. The issue of Vedic Mathematics involves religious politics, caste supremacy, apart from elementary arithmetic—so we cannot use simple statistics for our analysis. Further any study, when scientifically carried out using fuzzy models has more value than a statistical approach to the same. We used linguistic questionnaires for our data collection; experts filled in these questionnaires. In many cases, we also recorded our interviews with the experts in case they did not possess the technical knowledge of working with our questionnaire. Apart from this, several group discussions and meetings with various groups of people were held to construct the fuzzy models used to analyze this problem.

This book has five chapters. In Chapter I, we give a brief description of the sixteen sutras invented by the Swamiji. Chapter II gives the text of select articles about Vedic Mathematics that appeared in the media. Chapter III recalls some basic notions of some Fuzzy and Neutrosophic models used in this book. This chapter also introduces a fuzzy model to study the problem when we have to handle the opinion of multi-experts. Chapter IV analyses the problem using these models. The final chapter gives the observations made from our study.

The authors thank everybody who gave their opinion about Vedic Mathematics. Without their cooperation, the book could not have materialized. We next thank Dr.K.Kandasamy for proof-reading the book. I thank Meena and Kama for the layout and formatting of this book. Our thanks are also due to Prof. Praveen Prakash, Prof. Subrahmaniyam, Prof. E. L. Piriyakumar, Mr. Gajendran, Mr. S. Karuppasamy, Mr. Jayabhaskaran, Mr. Senguttuvan, Mr. Tamilselvan, Mr. D. Maariappan, Mr. P. Ganesan, Mr. N. Rajkumar and Ms. Rosalyn for the help rendered in various ways that could convert this book into a solid reality. We also thank the students of All India Students Federation (AISF) and the Students Federation of India (SFI) for their help in my work.

The authors dedicate this book to the great philosopher and intellectual Rahul Sangridyayan who revealed and exposed to the world many of the truths about the Vedas.

We have given a long list of references to help the interested reader.

<div style="text-align: right;">W.B.VASANTHA KANDASAMY<br>FLORENTIN SMARANDACHE</div>

**Chapter One**

# INTRODUCTION TO VEDIC MATHEMATICS

In this chapter we just recall some notions given in the book on Vedic Mathematics written by Jagadguru Swami Sri Bharati Krsna Tirthaji Maharaja (Sankaracharya of Govardhana Matha, Puri, Orissa, India), General Editor, Dr. V.S. Agrawala. Before we proceed to discuss the Vedic Mathematics that he professed we give a brief sketch of his heritage [51].

He was born in March 1884 to highly learned and pious parents. His father Sri P Narasimha Shastri was in service as a Tahsildar at Tinnivelly (Madras Presidency) and later retired as a Deputy Collector. His uncle, Sri Chandrasekhar Shastri was the principal of the Maharajas College, Vizianagaram and his great grandfather was Justice C. Ranganath Shastri of the Madras High Court. Born Venkatraman he grew up to be a brilliant student and invariably won the first place in all the subjects in all classes throughout his educational career. During his school days, he was a student of National College Trichanapalli; Church Missionary Society College, Tinnivelli and Hindu College Tinnivelly in Tamil Nadu. He passed his matriculation examination from the Madras University in 1899 topping the list as usual. His extraordinary proficiency in Sanskrit earned him the title "Saraswati" from the Madras Sanskrit Association in July 1899. After winning the highest place in the B.A examination Sri Venkataraman appeared for the M.A. examination of the American College of Sciences, Rochester, New York from the Bombay center in 1903. His subject of examination was Sanskrit, Philosophy, English, Mathematics, History and Science. He had a superb retentive memory.

In 1911 he could not anymore resist his burning desire for spiritual knowledge, practice and attainment and therefore, tearing himself off suddenly from the work of teaching, he went back to Sri Satcidananda Sivabhinava Nrisimha Bharati Swami at Sringeri. He spent the next eight years in the profoundest study of the most advanced Vedanta Philosophy and practice of the Brahmasadhana.

After several years in 1921 he was installed on the pontifical throne of Sharada Peetha Sankaracharya and later in 1925 he became the pontifical head of Sri Govardhan Math Puri where he served the remainder of his life spreading the holy spiritual teachings of Sanatana Dharma.

In 1957, when he decided finally to undertake a tour of the USA he rewrote from his memory the present volume of Vedic Mathematics [51] giving an introductory account of the sixteen formulae reconstructed by him. This is the only work on mathematics that has been left behind by him.

Now we proceed on to give the 16 sutras (aphorisms or formulae) and their corollaries [51]. As claimed by the editor, the list of these main 16 sutras and of their sub-sutras or corollaries is prefixed in the beginning of the text and the style of language also points to their discovery by Sri Swamiji himself. This is an open acknowledgement that they are not from the Vedas. Further the editor feels that at any rate it is needless to dwell longer on this point of origin since the vast merit of these rules should be a matter of discovery for each intelligent reader.

Now having known that even the 16 sutras are the Jagadguru Sankaracharya's invention we mention the name of the sutras and the sub sutras or corollaries as given in the book [51] pp. XVII to XVIII.

## Sixteen Sutras and their corollaries

| Sl. No | Sutras | Sub sutras or Corollaries |
|---|---|---|
| 1. | Ekādhikena Pūrvena (also a corollary) | Ānurūpyena |
| 2. | Nikhilam Navataścaramam Daśatah | Śisyate Śesamjnah |
| 3. | Ūrdhva - tiryagbhyām | Ādyamādyenantyamantyena |
| 4. | Parāvartya Yojayet | Kevalaih Saptakam Gunỹat |
| 5. | Sūnyam Samyasamuccaye | Vestanam |
| 6. | (Ānurūpye) Śūnyamanyat | Yāvadūnam Tāvadūnam |
| 7. | Sankalana - vyavakalanābhyām | Yāvadūnam Tāvadūnīkrtya Vargañca Yojayet |
| 8. | Puranāpuranābhyām | Antyayordasake' pi |
| 9. | Calanā kalanābhyām | Antyayoreva |
| 10. | Yāvadūnam | Samuccayagunitah |
| 11. | Vyastisamastih | Lopanasthāpanabhyām |
| 12. | Śesānyankena Caramena | Vilokanam |
| 13. | Sopantyadvayamantyam | Gunitasamuccayah Samuccayagunitah |
| 14. | Ekanyūnena Pūrvena |  |
| 15. | Gunitasamuccayah |  |
| 16. | Gunakasamuccayah |  |

The editor further adds that the list of 16 slokas has been complied from stray references in the text. Now we give spectacular illustrations and a brief descriptions of the sutras.

### The First Sutra: Ekādhikena Pūrvena

The relevant Sutra reads *Ekādhikena Pūrvena* which rendered into English simply says "By one more than the previous one". Its application and "modus operandi" are as follows.

(1) The last digit of the denominator in this case being 1 and the previous one being 1 "one more than the previous one"

evidently means 2. Further the proposition 'by' (in the sutra) indicates that the arithmetical operation prescribed is either multiplication or division. We illustrate this example from pp. 1 to 3. [51]

Let us first deal with the case of a fraction say 1/19. 1/19 where denominator ends in 9.

By the Vedic one - line mental method.

**A. First method**

$$\frac{1}{19} = \frac{.0\ 5\ 2\ 6\ 3\ 1\ 5\ 7\ 8\ 9\ 4\ 7\ 3\ 6\ 8\ 4\ 2\ i}{1\quad 1\quad\ \ 1\ 1\ 1\ 1\quad 1\quad\ 1\ 1}$$

**B. Second Method**

$$\frac{1}{19} = \frac{.0\ 5\ 2\ 6\ 3\ 1\ 5\ 7\ 8\ /\ 9\ 4\ 7\ 3\ 6\ 8\ 4\ 2\ i}{1\ 1\quad\ 1\ 1\ 1\ 1\quad\ \ 1\ 1\ 1}$$

This is the whole working. And the modus operandi is explained below.

**A. First Method**

Modus operandi chart is as follows:

(i)   We put down 1 as the right-hand most digit                1
(ii)  We multiply that last digit 1 by 2 and put the 2
      down as the immediately preceding digit.              2 1
(iii) We multiply that 2 by 2 and put 4 down as the
      next previous digit.                                  4 2 1
(iv)  We multiply that 4 by 2 and put it down thus          8 4 2 1
(v)   We multiply that 8 by 2 and get 16 as the
      product. But this has two digits. We therefore
      put the product. But this has two digits we
      therefore put the 6 down immediately to the
      left of the 8 and keep the 1 on hand to be
      carried over to the left at the next step (as we

|                                                                                                                                                                                                                                                                                             |           |
| ------------------------------------------------------------------------------------------------------------------------------------------------------------------------------------------------------------------------------------------------------------------------------------------- | --------- |
| always do in all multiplication e.g. of 69 × 2 = 138 and so on).                                                                                                                                                                                                                            | 6 8 4 2 1 |
|                                                                                                                                                                                                                                                                                             | 1         |
| (vi) We now multiply 6 by 2 get 12 as product, add thereto the 1 (kept to be carried over from the right at the last step), get 13 as the consolidated product, put the 3 down and keep the 1 on hand for carrying over to the left at the next step.                                        | 3 6 8 4 2 1 |
|                                                                                                                                                                                                                                                                                             | 1 1       |
| (vii) We then multiply 3 by 2 add the one carried over from the right one, get 7 as the consolidated product. But as this is a single digit number with nothing to carry over to the left, we put it down as our next multiplicand.                                                          | 7 3 6 8 4 2 1 |
|                                                                                                                                                                                                                                                                                             | 1 1       |
| ((viii) and xviii) we follow this procedure continually until we reach the 18<sup>th</sup> digit counting leftwards from the right, when we find that the whole decimal has begun to repeat itself. We therefore put up the usual recurring marks (dots) on the first and the last digit of the answer (from betokening that the whole of it is a Recurring Decimal) and stop the multiplication there. | |

Our chart now reads as follows:

$\frac{1}{19}$ =   . 0 5 2 6 3 1 5 7 8 / 9 4 7 3 6 8 4 2 i̇ .
         1   1     1 1 1 1 /   1   1 1

**B. Second Method**

The second method is the method of division (instead of multiplication) by the self-same "*Ekādhikena Pūrvena*" namely 2. And as division is the exact opposite of multiplication it stands to reason that the operation of division should proceed not from right to left (as in the case of multiplication as expounded here in before) but in the exactly opposite direction; i.e. from left to right. And such is actually found to be the case. Its application and modus operandi are as follows:

(i) Dividing 1 (The first digit of the dividend) by 2, we see the quotient is zero and the remainder is 1. We therefore set 0 down as the first digit of the quotient and prefix the remainder 1 to that very digit of the quotient (as a sort of reverse-procedure to the carrying to the left process used in multiplication) and thus obtain 10 as our next dividend.

   0
   1

(ii) Dividing this 10 by 2, we get 5 as the second digit of the quotient, and as there is no remainder to be prefixed thereto we take up that digit 5 itself as our next dividend.

. 0 5
  1

(iii) So, the next quotient – digit is 2, and the remainder is 1. We therefore put 2 down as the third digit of the quotient and prefix the remainder 1 to that quotient digit 2 and thus have 12 as our next dividend.

. 0 5 2
  1  1

(iv) This gives us 6 as quotient digit and zero as remainder. So we set 6 down as the fourth digit of the quotient, and as there is no remainder to be prefixed thereto we take 6 itself as our next digit for division which gives the next quotient digit as 3.

. 0 5 2 6 3 1
  1 1   1

(v) That gives us 1 and 1 as quotient and remainder respectively. We therefore put 1 down as the 6$^{th}$ quotient digit prefix the 1 thereto and have 11 as our next dividend.

. 0 5 2 6 3 1 5
  1 1   1 1

(vi to xvii) Carrying this process of straight continuous division by 2 we get 2 as the 17<sup>th</sup> quotient digit and 0 as remainder.

(xviii) Dividing this 2 by 2 are get 1 as 18<sup>th</sup> quotient digit and 0 as remainder. But this is exactly what we began with. This means that the decimal begins to repeat itself from here. So we stop the mental division process and put down the usual recurring symbols (dots) on the 1<sup>st</sup> and 18<sup>th</sup> digit to show that the whole of it is a circulating decimal.

$$\left.\begin{array}{l} .0\ 5\ 2\ 6\ 3\ 1\ 5\ 7\ 8 \\ 1\quad 1\quad\quad 1\ 1\ 1\ 1 \\ 9\ 4\ 7\ 3\ 6\ 8\ 4\ 2\ i \\ \quad 1\quad 1\ 1 \end{array}\right\}$$

Now if we are interested to find 1/29 the student should note down that the last digit of the denominator is 9, but the penultimate one is 2 and one more than that means 3. Likewise for 1/49 the last digit of the denominator is 9 but penultimate is 4 and one more than that is 5 so for each number the observation must be memorized by the student and remembered.

The following are to be noted

1. Student should find out the procedure to be followed. The technique must be memorized. They feel it is difficult and cumbersome and wastes their time and repels them from mathematics.

2. "This problem can be solved by a calculator in a time less than a second. Who in this modernized world take so much strain to work and waste time over such simple calculation?" asked several of the students.

3. According to many students the long division method was itself more interesting.

**The Second Sutra: Nikhilam Navataścaramam Daśatah**

Now we proceed on to the next sutra "*Nikhilam sutra*" The sutra reads "*Nikhilam Navataścaramam Daśatah*", which literally translated means: all from 9 and the last from 10". We shall presently give the detailed explanation presently of the meaning and applications of this cryptical-sounding formula [51] and then give details about the three corollaries.

He has given a very simple multiplication.

Suppose we have to multiply 9 by 7.
1. We should take, as base for our calculations that power of 10 which is nearest to the numbers to be multiplied. In this case 10 itself is that power.  
   $$\begin{array}{r}(10)\\9-1\\\underline{7-3}\\6\,/\,3\end{array}$$
2. Put the numbers 9 and 7 above and below on the left hand side (as shown in the working alongside here on the right hand side margin);
3. Subtract each of them from the base (10) and write down the remainders (1 and 3) on the right hand side with a connecting minus sign (–) between them, to show that the numbers to be multiplied are both of them less than 10.
4. The product will have two parts, one on the left side and one on the right. A vertical dividing line may be drawn for the purpose of demarcation of the two parts.
5. Now, the left hand side digit can be arrived at in one of the 4 ways
   a) Subtract the base 10 from the sum of the given numbers (9 and 7 i.e. 16). And put (16 – 10) i.e. 6 as the left hand part of the answer $\qquad 9 + 7 - 10 = 6$
   or b) Subtract the sum of two deficiencies (1 + 3 = 4) from the base (10) you get the same answer (6) again $\qquad 10 - 1 - 3 = 6$
   or c) Cross subtract deficiency 3 on the second row from the original number 9 in the first row. And you find that you have got (9 – 3) i.e. 6 again $\qquad 9 - 3 = 6$
   or d) Cross subtract in the converse way (i.e. 1 from 7), and you get 6 again as the left hand side portion of the required answer $\qquad 7 - 1 = 6.$

*Note:* This availability of the same result in several easy ways is a very common feature of the Vedic system and is great advantage and help to the student as it enables him to test and verify the correctness of his answer step by step.

6. Now vertically multiply the two deficit figures (1 and 3). The product is 3. And this is the right hand side portion of the answer  (10) 9 – 1
7. Thus 9 × 7 = 63.  7 – 3
 ─────
 6 / 3

This method holds good in all cases and is therefore capable of infinite application. Now we proceed on to give the interpretation and working of the 'Nikhilam' sutra and its three corollaries.

**The First Corollary**

The first corollary naturally arising out of the Nikhilam Sutra reads in English "*whatever the extent of its deficiency lessen it still further to that very extent, and also set up the square of that deficiency*".

This evidently deals with the squaring of the numbers. A few elementary examples will suffice to make its meaning and application clear:
 Suppose one wants to square 9, the following are the successive stages in our mental working.

(i) We would take up the nearest power of 10, i.e. 10 itself as our base.
(ii) As 9 is 1 less than 10 we should decrease it still further by 1 and set 8 down as our left side portion of the answer  
 8/
(iii) And on the right hand we put down the square of that deficiency $1^2$  8/1.
(iv) Thus $9^2 = 81$   9 – 1
   9 – 1
   ─────
   8 / 1

Now we proceed on to give second corollary from (p.27, [51]).

**The Second Corollary**

The second corollary in applicable only to a special case under the first corollary i.e. the squaring of numbers ending in 5 and other cognate numbers. Its wording is exactly the same as that of the sutra which we used at the outset for the conversion of 'vulgar' fractions into their recurring decimal equivalents. The sutra now takes a totally different meaning and in fact relates to a wholly different setup and context.

Its literal meaning is the same as before (i.e. by one more than the previous one") but it now relates to the squaring of numbers ending in 5. For example we want to multiply 15. Here the last digit is 5 and the "previous" one is 1. So one more than that is 2.

Now sutra in this context tells us to multiply the previous digit by one more than itself i.e. by 2. So the left hand side digit is $1 \times 2$ and the right hand side is the vertical multiplication product i.e. 25 as usual.

$$\frac{1 \ / \ 5}{2 \ / \ 25}$$

Thus $15^2 = 1 \times 2 \ / \ 25 = 2 \ / \ 25$.

Now we proceed on to give the third corollary.

**The Third Corollary**

Then comes the third corollary to the Nikhilam sutra which relates to a very special type of multiplication and which is not frequently in requisition elsewhere but is often required in mathematical astronomy etc. It relates to and provides for multiplications where the multiplier digits consists entirely of nines.

The procedure applicable in this case is therefore evidently as follows:

i)  Divide the multiplicand off by a vertical line into a right hand portion consisting of as many digits as the multiplier;

and subtract from the multiplicand one more than the whole excess portion on the left. This gives us the left hand side portion of the product;

or take the Ekanyuna and subtract therefrom the previous i.e. the excess portion on the left; and

ii) Subtract the right hand side part of the multiplicand by the Nikhilam rule. This will give you the right hand side of the product.

The following example will make it clear:

$$
\begin{array}{r}
43 \times 9 \\
4 : 3 : \\
:-5 : 3 \\
\hline
3 : 8 : 7
\end{array}
$$

### The Third Sutra: Ūrdhva Tiryagbhyām

*Ūrdhva Tiryagbhyām* sutra which is the General Formula applicable to all cases of multiplication and will also be found very useful later on in the division of a large number by another large number.

The formula itself is very short and terse, consisting of only one compound word and means "vertically and cross-wise." The applications of this brief and terse sutra are manifold.

A simple example will suffice to clarify the modus operandi thereof. Suppose we have to multiply 12 by 13.

(i) We multiply the left hand most digit 1 of the multiplicand vertically by the left hand most digit 1 of the multiplier get their product 1 and set down as the left hand most part of the answer;

$$
\begin{array}{r}
12 \\
13 \\
\hline
1:3 + 2:6 = 156
\end{array}
$$

(ii) We then multiply 1 and 3 and 1 and 2 crosswise add the two get 5 as the sum and set it down as the middle part of the answer; and

(iii) We multiply 2 and 3 vertically get 6 as their product and put it down as the last the right hand most part of the answer. Thus $12 \times 13 = 156$.

## The Fourth Sutra: Parāvartya Yojayet

The term *Parāvartya Yojayet* which means "Transpose and Apply." Here he claims that the Vedic system gave a number is applications one of which is discussed here. The very acceptance of the existence of polynomials and the consequent remainder theorem during the Vedic times is a big question so we don't wish to give this application to those polynomials. However the four steps given by them in the polynomial division are given below: Divide $x^3 + 7x^2 + 6x + 5$ by $x - 2$.

i. $x^3$ divided by x gives us $x^2$ which is therefore the first term of the quotient
$$\frac{x^3 + 7x^2 + 6x + 5}{x - 2} \quad \therefore Q = x^2 + \ldots.$$

ii. $x^2 \times -2 = -2x^2$ but we have $7x^2$ in the divident. This means that we have to get $9x^2$ more. This must result from the multiplication of x by 9x. Hence the 2$^{nd}$ term of the divisor must be 9x
$$\frac{x^3 + 7x^2 + 6x + 5}{x - 2} \quad \therefore Q = x^2 + 9x + \ldots.$$

iii. As for the third term we already have $-2 \times 9x = -18x$. But we have 6x in the dividend. We must therefore get an additional 24x. Thus can only come in by the multiplication of x by 24. This is the third term of the quotient.
$\therefore \quad Q = x^2 + 9x + 24$

iv. Now the last term of the quotient multiplied by $-2$ gives us $-48$. But the absolute term in the dividend is 5. We have therefore to get an additional 53 from some where. But there is no further term left in the dividend. This means that the 53 will remain as the remainder $\therefore Q = x^2 + 9x + 24$ and $R = 53$.

This method for a general degree is not given. However this does not involve anything new. Further is it even possible that the concept of polynomials existed during the period of Vedas itself?

Now we give the 5$^{th}$ sutra.

## The Fifth Sutra: Sūnyam Samyasamuccaye

We begin this section with an exposition of several special types of equations which can be practically solved at sight with the aid of a beautiful special sutra which reads *Sūnyam Samyasamuccaye* and which in cryptic language which renders its applicable to a large number of different cases. It merely says "when the Samuccaya is the same that Samuccaya is zero i.e. it should be equated to zero."

Samuccaya is a technical term which has several meanings in different contexts which we shall explain one at a time.

Samuccaya firstly means a term which occurs as a common factor in all the terms concerned.

Samuccaya secondly means the product of independent terms.

Samuccaya thirdly means the sum of the denominators of two fractions having same numerical numerator.

Fourthly Samuccaya means combination or total.

Fifth meaning: With the same meaning i.e. total of the word (Samuccaya) there is a fifth kind of application possible with quadratic equations.

Sixth meaning – With the same sense (total of the word – Samuccaya) but in a different application it comes in handy to solve harder equations equated to zero.

Thus one has to imagine how the six shades of meanings have been perceived by the Jagadguru Sankaracharya that too from the Vedas when such types of equations had not even been invented in the world at that point of time. However the immediate application of the subsutra Vestnam is not given but extensions of this sutra are discussed.

So we next go to the sixth sutra given by His Holiness Sankaracharya.

## The Sixth Sutra: Ānurūpye Śūnyamanyat

As said by Dani [32] we see the $6^{th}$ sutra happens to be the subsutra of the first sutra. Its mention is made in {pp. 51, 74, 249 and 286 of [51]}. The two small subsutras (i) *Anurpyena* and (ii) *Adayamadyenantyamantyena* of the sutras 1 and 3 which mean "proportionately" and "the first by the first and the last by the last".

Here the later subsutra acquires a new and beautiful double application and significance. It works out as follows:

i. Split the middle coefficient into two such parts so that the ratio of the first coefficient to the first part is the same as the ratio of that second part to the last coefficient. Thus in the quadratic $2x^2 + 5x + 2$ the middle term 5 is split into two such parts 4 and 1 so that the ratio of the first coefficient to the first part of the middle coefficient i.e. 2 : 4 and the ratio of the second part to the last coefficient i.e. 1 : 2 are the same. Now this ratio i.e. $x + 2$ is one factor.

ii. And the second factor is obtained by dividing the first coefficient of the quadratic by the first coefficient of the factor already found and the last coefficient of the quadratic by the last coefficient of that factor. In other words the second binomial factor is obtained thus

$$\frac{2x^2}{x} + \frac{2}{2} = 2x + 1.$$

Thus $2x^2 + 5x + 2 = (x + 2)(2x + 1)$. This sutra has *Yavadunam Tavadunam* to be its subsutra which the book claims to have been used.

## The Seventh Sutra: Sankalana Vyavakalanābhyām

Sankalana Vyavakalan process and the Adyamadya rule together from the seventh sutra. The procedure adopted is one of alternate destruction of the highest and the lowest powers by a suitable multiplication of the coefficients and the addition or subtraction of the multiples.

A concrete example will elucidate the process.

Suppose we have to find the HCF (Highest Common factor) of $(x^2 + 7x + 6)$ and $x^2 - 5x - 6$.

$$x^2 + 7x + 6 = (x + 1)(x + 6) \text{ and}$$
$$x^2 - 5x - 6 = (x + 1)(x - 6)$$

$\therefore$ the HCF is $x + 1$

but where the sutra is deployed is not clear.

This has a subsutra *Yavadunam Tavadunikrtya*. However it is not mentioned in chapter 10 of Vedic Mathematics [51].

### The Eight Sutra: Puranāpuranābhyām

*Puranāpuranābhyām* means "by the completion or not completion" of the square or the cube or forth power etc. But when the very existence of polynomials, quadratic equations etc. was not defined it is a miracle the Jagadguru could contemplate of the completion of squares (quadratic) cubic and forth degree equation. This has a subsutra *Antyayor dasake'pi* use of which is not mentioned in that section.

### The Ninth Sutra: Calanā kalanābhyām

The term (*Calanā kalanābhyām*) means differential calculus according to Jagadguru Sankaracharya. It is mentioned in page 178 [51] that this topic will be dealt with later on. We have not dealt with it as differential calculus not pertaining to our analysis as it means only differential calculus and has no mathematical formula or sutra value.

### The Tenth Sutra: Yāvadūnam

*Yāvadūnam Sutra* (for cubing) is the tenth sutra. However no modus operandi for elementary squaring and cubing is given in this book [51]. It has a subsutra called *Samuccayagunitah*.

### The Eleventh Sutra: Vyastisamastih Sutra

Vyastisamastih sutra teaches one how to use the average or exact middle binomial for breaking the biquadratic down into a simple quadratic by the easy device of mutual cancellations of the odd powers. However the modus operandi is missing.

## The Twelfth Sutra: Śesānyankena Caramena

The sutra *Śesānyankena Caramena* means "The remainders by the last digit". For instance if one wants to find decimal value of 1/7. The remainders are 3, 2, 6, 4, 5 and 1. Multiplied by 7 these remainders give successively 21, 14, 42, 28, 35 and 7. Ignoring the left hand side digits we simply put down the last digit of each product and we get 1/7 = .14 28 57!

Now this 12$^{th}$ sutra has a subsutra *Vilokanam*. *Vilokanam* means "mere observation" He has given a few trivial examples for the same.

Next we proceed on to study the 13$^{th}$ sutra Sopantyadvayamantyam.

## The Thirteen Sutra: *Sopantyadvayamantyam*

The sutra *Sopantyadvayamantyam* means "the ultimate and twice the penultimate" which gives the answer immediately. No mention is made about the immediate subsutra.

The illustration given by them.

$$\frac{1}{(x+2)(x+3)} + \frac{1}{(x+2)(x+4)} = \frac{1}{(x+2)(x+5)} + \frac{1}{(x+3)(x+4)}.$$

Here according to this sutra L + 2P (the last + twice the penultimate)
= (x + 5) + 2 (x + 4) = 3x + 13 = 0
∴ x = $-4\frac{1}{3}$ .

The proof of this is as follows.

$$\frac{1}{(x+2)(x+3)} + \frac{1}{(x+2)(x+4)} = \frac{1}{(x+2)(x+5)} + \frac{1}{(x+3)(x+4)}$$

$$\therefore \frac{1}{(x+2)(x+3)} - \frac{1}{(x+2)(x+5)} = \frac{1}{(x+3)(x+4)} - \frac{1}{(x+2)(x+4)}$$

$$\therefore \frac{1}{(x+2)}\left[\frac{2}{(x+3)(x+5)}\right] = \frac{1}{(x+4)}\left[\frac{-1}{(x+2)(x+3)}\right]$$

Removing the factors (x + 2) and (x + 3);

$$\frac{2}{x+5} = \frac{-1}{x+4} \quad \text{i.e.} \quad \frac{2}{L} = \frac{-1}{P}$$
$$\therefore L + 2P = 0.$$

The General Algebraic Proof is as follows.
$$\frac{1}{AB} + \frac{1}{AC} = \frac{1}{AD} + \frac{1}{BC}$$
(where A, B, C and D are in A.P).

Let d be the common difference
$$\frac{1}{A(A+d)} + \frac{1}{A(A+2d)} = \frac{1}{A(A+3d)} + \frac{1}{(A+d)(A+2d)}$$
$$\therefore \frac{1}{A(A+d)} - \frac{1}{A(A+3d)} = \frac{1}{(A+d)(A+2d)} + \frac{1}{A(A+2d)}$$
$$\therefore \frac{1}{A}\left\{\frac{2d}{(A+d)(A+3d)}\right\} = \frac{1}{(A+2d)}\left\{\frac{-d}{A(A+d)}\right\}.$$

Canceling the factors $A(A+d)$ of the denominators and d of the numerators:
$$\therefore \frac{2}{A+3d} = \frac{-1}{A+2d} \quad \text{(p. 137)}$$
In other words $\frac{2}{L} = \frac{-1}{P}$
$$\therefore L + 2P = 0$$
It is a pity that all samples given by the book form a special pattern.

We now proceed on to present the 14$^{th}$ Sutra.

### The Fourteenth Sutra: Ekanyūnena Pūrvena

The *Ekanyūnena Pūrvena* Sutra sounds as if it were the converse of the Ekadhika Sutra. It actually relates and provides for multiplications where the multiplier the digits consists entirely of nines. The procedure applicable in this case is therefore evidently as follows.

For instance 43 × 9.

i. Divide the multiplicand off by a vertical line into a right hand portion consisting of as many digits as the multiplier; and subtract from the multiplicand one more than the whole excess portion on the left. This gives us the left hand side portion of the product or take the Ekanyuna and subtract it from the previous i.e. the excess portion on the left and
ii. Subtract the right hand side part of the multiplicand by the Nikhilam rule. This will give you the right hand side of the product

$$43 \times 9$$
$$4 : 3$$
$$:-5 : 3$$
$$3 : 8 : 7$$

This Ekanyuna Sutra can be utilized for the purpose of postulating mental one-line answers to the question.

We now go to the 15$^{th}$ Sutra.

### The Fifthteen Sutra: Gunitasamuccayah

*Gunitasamuccayah* rule i.e. the principle already explained with regard to the $S_c$ of the product being the same as the product of the $S_c$ of the factors.

Let us take a concrete example and see how this method (p. 81) [51] can be made use of. Suppose we have to factorize $x^3 + 6x^2 + 11x + 6$ and by some method, we know $(x + 1)$ to be a factor. We first use the corollary of the 3$^{rd}$ sutra viz. Adayamadyena formula and thus mechanically put down $x^2$ and 6 as the first and the last coefficients in the quotient; i.e. the product of the remaining two binomial factors. But we know already that the $S_c$ of the given expression is 24 and as the $S_c$ of $(x + 1) = 2$ we therefore know that the $S_c$ of the quotient must be 12. And as the first and the last digits thereof are already known to be 1 and 6, their total is 7. And therefore the middle term must be 12 – 7 = 5. So, the quotient $x^2 + 5x + 6$.

This is a very simple and easy but absolutely certain and effective process.

As per pp. XVII to XVIII [51] of the book there is no corollary to the 15$^{th}$ sutra i.e. to the sutra Gunitasamuccayah but in p. 82 [51] of the same book they have given under the title corollaries 8 methods of factorization which makes use of mainly the Adyamadyena sutra. The interested reader can refer pp. 82-85 of [51].

Now we proceed on to give the last sutra enlisted in page XVIII of the book [51].

### The Sixteen Sutra :Gunakasamuccayah.

"It means the product of the sum of the coefficients in the factors is equal to the sum of the coefficients in the product".

In symbols we may put this principle as follows:
 $S_c$ of the product = Product of the $S_c$ (in factors).
For example
$$(x + 7)(x + 9) = x^2 + 16x + 63$$
and we observe
$$(1 + 7)(1 + 9) = 1 + 16 + 63 = 80.$$

Similarly in the case of cubics, biquadratics etc. the same rule holds good.
For example
$$(x + 1)(x + 2)(x + 3) = x^3 + 6x^2 + 11x + 6$$
$$2 \times 3 \times 4 = 1 + 6 + 11 + 6$$
$$= 24.$$

Thus if and when some factors are known this rule helps us to fill in the gaps.

It will be found useful in the factorization of cubics, biquadratics and will also be discussed in some other such contexts later on.

In several places in the use of sutras the corollaries are subsutras are dealt separately. One such instance is the subsutra of the 11$^{th}$ sutra i.e., *Vyastisamastih* and its corollary viz. *Lapanasthapanabhyam* finds its mention in page 77 [51] which is cited verbatim here. The *Lapana Sthapana* subsutra however removes the whole difficulty and makes the factorization of a quadratic of this type as easy and simple as that of the ordinary quadratic already explained. The procedure is as follows: Suppose we have to factorise the following long quadratic.

$$2x^2 + 6y^2 + 6z^2 + 7xy + 11yz + 7zx$$

i. We first eliminate by putting z = 0 and retain only x and y and factorise the resulting ordinary quadratic in x and y with Adyam sutra which is only a corollary to the 3$^{rd}$ sutra viz. Urdhva tryyagbhyam.
ii. We then similarly eliminate y and retain only x and z and factorise the simple quadratic in x and z.
iii. With these two sets of factors before us we fill in the gaps caused by our own deliberate elimination of z and y respectively. And that gives us the real factors of the given long expression. The procedure is an argumentative one and is as follows:

If z = 0 then the given expression is $2x^2 + 7xy + 6y^2 = (x + 2y)(2x + 3y)$. Similarly if y = 0 then $2x^2 + 7xz + 3z^2 = (x + 3z)(2x + z)$.
Filling in the gaps which we ourselves have created by leaving out z and y, we get E = (x + 2y + 3z) (2x + 3y + z)

**Note:**

This Lopanasthapana method of alternate elimination and retention will be found highly useful later on in finding HCF, in solid geometry and in co-ordinate geometry of the straight line, the hyperbola, the conjugate hyperbola, the asymptotes etc.
In the current system of mathematics we have two methods which are used for finding the HCF of two or more given expressions.
The first is by means of factorization which is not always easy and the second is by a process of continuous division like the method used in the G.C.M chapter of arithmetic. The latter is a mechanical process and can therefore be applied in all cases. But it is rather too mechanical and consequently long and cumbrous.

The Vedic methods provides a third method which is applicable to all cases and is at the same time free from this disadvantage.

It is mainly an application of the subsutras or corollaries of the 11$^{th}$ sutra viz. *Vyastisamastih*, the corollary *Lapanasthapana* sutra the 7$^{th}$ sutra viz. *Sankalana Vyavakalanabhyam* process and the subsutra of the 3$^{rd}$ sutra viz. *Adyamādyenantyamantyena*.

The procedure adopted is one of alternate destruction of the highest and the lowest powers by a suitable multiplication of the coefficients and the addition or subtraction of the multiples.

A concrete example will elucidate the process.

Suppose we have to find the H.C.F of $x^2 + 7x + 6$ and $x^2 - 5x - 6$

i. $x^2 + 7x + 6 = (x + 1)(x + 6)$ and $x^2 - 5x - 6 = (x + 1)(x - 6)$. HCF is $(x + 1)$. This is the first method.
ii. The second method the GCM one is well-known and need not be put down here.
iii. The third process of 'Lopanasthapana' i.e. of the elimination and retention or alternate destruction of the highest and the lowest powers is explained below.

Let $E_1$ and $E_2$ be the two expressions. Then for destroying the highest power we should substract $E_2$ from $E_1$ and for destroying the lowest one we should add the two. The chart is as follows:

$$\left.\begin{array}{l} x^2 + 7x + 6 \\ x^2 - 5x - 6 \end{array}\right\} \text{subtraction} \qquad \left.\begin{array}{l} x^2 - 5x - 6 \\ x^2 + 7x + 6 \end{array}\right\} \text{addition}$$

$$\begin{array}{l} \overline{12x + 12} \\ 12)\,\underline{12x + 12} \\ \phantom{12)\,}x + 1 \end{array} \qquad \begin{array}{l} \overline{2x^2 + 2x} \\ 2x)\,\underline{2x^2 + 2x} \\ \phantom{2x)\,}x + 1 \end{array}$$

We then remove the common factor if any from each and we find $x + 1$ staring us in the face i.e. $x + 1$ is the HCF. Two things are to be noted importantly.

(1) We see that often the subsutras are not used under the main sutra for which it is the subsutra or the corollary. This is the main deviation from the usual mathematical principles of theorem (sutra) and corollaries (subsutra).

(2) It cannot be easily compromised that a single sutra (a Sanskrit word) can be mathematically interpreted in this manner even by a stalwart in Sanskrit except the Jagadguru Puri Sankaracharya.

We wind up the material from the book of Vedic Mathematics and proceed on to give the opinion/views of great personalities on Vedic Mathematics given by Jagadguru.
 Since the notion of integral and differential calculus was not in vogue in Vedic times, here we do not discuss about the authenticated inventor, further we have not given the adaptation of certain sutras in these fields. Further as most of the educated experts felt that since the Jagadguru had obtained his degree with mathematics as one of the subjects, most of the results given in book on Vedic Mathematics were manipulated by His Holiness.

**Chapter Two**

# ANALYSIS OF VEDIC MATHEMATICS BY MATHEMATICIANS AND OTHERS

In this chapter we give the verbatim opinion of mathematicians and experts about Vedic Mathematics in their articles, that have appeared in the print media. The article of Prof. S.G. Dani, School of Mathematics, Tata Institute of Fundamental Research happen to give a complete analysis of Vedic Mathematics.

We have given his second article verbatim because we do not want any bias or our opinion to play any role in our analysis [32].

However we do not promise to discuss all the articles. Only articles which show "How Vedic is Vedic Mathematics?" is given for the perusal of the reader. We thank them for their articles and quote them verbatim. The book on Vedic Mathematics by Jagadguru Sankaracharya of Puri has been translated into Tamil by Dr. V.S. Narasimhan, a Retired Professor of an arts college and C. Mailvanan, M.Sc Mathematics (Vidya Barathi state-level Vedic Mathematics expert) in two volumes. The first edition appeared in 1998 and the corrected second edition in 2003.

In Volume I of the Tamil book the introduction is as follows: "Why was the name Vedic Mathematics given? On the title "a trick in the name of Vedic Mathematics" though professors in mathematics praise the sutras, they argue that the title Vedic Mathematics is not well suited. According to them the sutras published by the Swamiji are not found anywhere in the Vedas. Further the branches of mathematics like algebra and calculus which he mentions, did not exist in the Vedic times. It may help school students but only in certain problems where shortcut methods can be used. The Exaggeration that, it can be used in all branches of mathematics cannot be accepted.

Because it gives answers very fast it can be called "speed maths". He has welcomed suggestions and opinions of one and all.

It has also become pertinent to mention here that Jagadguru Puri Sankaracharya for the first time visited the west in 1958. He had been to America at the invitation of the Self Realization Fellowship Los Angeles, to spread the message of Vedanta. The book Vedic Metaphysics is a compilation of some of his discourses delivered there. On 19 February 1958, he has given a talk and demonstration to a small group of student mathematicians at the California Institute of Technology, Pasadena, California.

This talk finds its place in chapter XII of the book Vedic Metaphysics pp. 156-196 [52] most of which has appeared later on, in his book on Vedic Mathematics [51]. However some experts were of the opinion, that if Swamiji would have remained as Swamiji 'or' as a 'mathematician' it would have been better. His intermingling and trying to look like both has only brought him less recognition in both Mathematics and on Vedanta. The views of Wing Commander Vishva Mohan Tiwari, under the titles conventional to unconventionally original speaks of Vedic Mathematics as follows:

"Vedic Mathematics mainly deals with various Vedic mathematical formulas and their applications of carrying out tedious and cumbersome arithmetical operations, and to a very large extent executing them mentally. He feels that in this field of mental arithmetical operations the works of the famous mathematicians Trachtenberg and Lester Meyers (High speed mathematics) are elementary compared to that of Jagadguruji … An attempt has been made in this note to explain the unconventional aspects of the methods. He then gives a very brief sketch of first four chapters of Vedic Mathematics".

This chapter has seven sections; Section one gives the verbatim analysis of Vedic Mathematics given by Prof. Dani in his article in Frontline [31].

A list of eminent signatories asking people to stop this fraud on our children is given verbatim in section two. Some views given about the book both in favour of and against is given in section three.

Section four gives the essay Vedas: Repositories of ancient lore. "A rational approach to study ancient literature" an article found in Current Science, volume 87, August 2004 is given in Section five. Section Six gives the "Shanghai Rankings and Indian Universities." The final section gives conclusion derived on Vedic Mathematics and calculation of Guru Tirthaji.

## 2.1 Views of Prof. S.G. Dani about Vedic Mathematics from Frontline

Views of Prof. S.G.Dani gave the authors a greater technical insight into Vedic Mathematics because he has written 2 articles in Frontline in 1993. He has analyzed the book extremely well and we deeply acknowledge the services of professor S.G.Dani to the educated community in general and school students in particular. This section contains the verbatim views of Prof. Dani that appeared in Frontline magazine. He has given a marvelous analysis of the book Vedic Mathematics and has daringly concluded.

"One would hardly have imagine that a book which is transparently not from any ancient source or of any great mathematical significance would one day be passed off as a storehouse of some ancient mathematical treasure. It is high time saner elements joined hands to educate people on the truth of this so-called Vedic Mathematics and prevent the use of public money and energy on its propagation, beyond the limited extent that may be deserved, lest the intellectual and educational life in the country should get vitiated further and result in wrong attitudes to both history and mathematics, especially in the coming generation."

**Myths and Reality: On 'Vedic Mathematics'**
S.G. Dani, School of Mathematics,
Tata Institute of Fundamental Research
*An updated version of the 2-part article in Frontline, 22 Oct. and 5 Nov. 1993*

We in India have good reasons to be proud of a rich heritage in science, philosophy and culture in general, coming to us down the ages. In mathematics, which is my own area of specialization, the ancient Indians not only took great strides long before the Greek advent, which is a standard reference point in the Western historical perspective, but also enriched it for a long period making in particular some very fundamental contributions such as the place-value system for writing numbers as we have today, introduction of zero and so on. Further, the sustained development of mathematics in India in the post-Greek period was indirectly instrumental in the revival in Europe after "its dark ages".

Notwithstanding the enviable background, lack of adequate attention to academic pursuits over a prolonged period, occasioned by several factors, together with about two centuries of Macaulayan educational system, has unfortunately resulted, on the one hand, in a lack of awareness of our historical role in actual terms and, on the other, an empty sense of pride which is more of an emotional reaction to the colonial domination rather than an intellectual challenge. Together they provide a convenient ground for extremist and misguided elements in society to "reconstruct history" from nonexistent or concocted source material to whip up popular euphoria.

That this anti-intellectual endeavour is counter-productive in the long run and, more important, harmful to our image as a mature society, is either not recognized or ignored in favour of short-term considerations. Along with the obvious need to accelerate the process of creating an awareness of our past achievements, on the strength of authentic information, a more urgent need has also arisen to confront and expose such baseless constructs before it is too late. This is not merely a question of setting the record straight. The motivated versions have a way of corrupting the intellectual processes in society and weakening their very foundations in the long run, which needs to be prevented at all costs. The so-called "Vedic Mathematics" is a case in point. A book by that name written by Jagadguru Swami Shri Bharati Krishna Tirthaji Maharaja (Tirthaji, 1965) is at the centre of this pursuit, which has now acquired wide following; Tirthaji was the Shankaracharya of Govardhan Math, Puri, from 1925 until he passed away in 1960. The book was published posthumously, but he had been carrying out a campaign on the theme for a long time, apparently for several decades, by means of lectures, blackboard demonstrations, classes and so on. It has been known from the beginning that there is no evidence of the contents of the book being of Vedic origin; the Foreword to the book by the General Editor, Dr. A.S.Agrawala, and an account of the genesis of the work written by Manjula Trivedi, a disciple of the swamiji, make this clear even before one gets to the text of the book. No one has come up with any positive evidence subsequently either.

There has, however, been a persistent propaganda that the material is from the Vedas. In the face of a false sense of national pride associated with it and the neglect, on the part of the knowledgeable, in countering the propaganda, even educated and well meaning people have tended to accept it uncritically. The vested interests have also involved politicians in the propaganda process to gain state support. Several leaders have lent support to the "Vedic Mathematics" over the years, evidently in the belief of its being from ancient scriptures. In the current environment, when a label as ancient seems to carry considerable premium irrespective of its authenticity or merit, the purveyors would have it going easy.

Large sums have been spent both by the Government and several private agencies to support this "Vedic Mathematics", while authentic Vedic studies continue to be neglected. People, especially children, are encouraged to learn and spread the contents of the book, largely on the baseless premise of their being from the Vedas. With missionary zeal several "devotees" of this cause have striven to take the "message" around the world; not surprisingly, they have even met with some success in the West, not unlike some of the gurus and yogis peddling their own versions of "Indian philosophy". Several people are also engaged in "research" in the new "Vedic Mathematics."

To top it all, when in the early nineties the Uttar Pradesh Government introduced "Vedic Mathematics" in school text books, the contents of the swamiji's book were treated as if they were genuinely from the Vedas; this also naturally seems to have led them to include a list of the swamiji's sutras on one of the opening pages (presumably for the students to learn them by heart and recite!) and to accord the swamiji a place of honour in the "brief history of Indian mathematics" described in the beginning of the textbook, together with a chart, which curiously has Srinivasa Ramanujan's as the only other name from the twentieth century!

For all their concern to inculcate a sense of national pride in children, those responsible for this have not cared for the simple fact that modern India has also produced several notable mathematicians and built a worthwhile edifice in mathematics (as also in many other areas). Harish Chandra's work is held in great esteem all over the world and several leading seats of learning of our times pride themselves in having members pursuing his ideas; (see, for instance, Langlands, 1993). Even among those based in India, several like Syamdas Mukhopadhyay, Ganesh Prasad, B.N.Prasad, K.Anand Rau, T.Vijayaraghavan, S.S.Pillai, S.Minakshisundaram, Hansraj Gupta, K.G.Ramanathan, B.S.Madhava Rao, V.V.Narlikar, P.L.Bhatnagar and so on and also many living Indian mathematicians have carved a niche for themselves on the international mathematical scene (see Narasimhan, 1991). Ignoring all this while introducing the swamiji's name in the "brief history" would inevitably create a warped perspective in children's minds, favouring gimmickry rather than professional work. What does the swamiji's "Vedic Mathematics" seek to do and what does it achieve? In his preface of the book, grandly titled" A Descriptive Prefatory Note on the astounding Wonders of Ancient Indian Vedic Mathematics," the swamiji tells us that he strove from his childhood to study the Vedas critically "to prove to ourselves (and to others) the correctness (or otherwise)"of the "derivational meaning" of "Veda" that the" Vedas should contain within themselves all the knowledge needed by the mankind relating not only to spiritual matters but also those usually described as purely 'secular', 'temporal' or 'worldly'; in other words, simply because of the meaning of the word 'Veda', everything that is worth knowing is expected to be contained in the vedas and the swamiji seeks to prove it to be the case!

It may be worthwhile to point out here that there would be room for starting such an enterprise with the word 'science'! He also describes how the "contemptuous or at best patronising" attitude of Orientalists, Indologists and so on strengthened his determination to unravel the too-long-hidden mysteries of philosophy and science contained in ancient India's Vedic lore, with the consequence that, "after eight years of concentrated contemplation in forest solitude, we were at long last able to recover the long lost keys which alone could unlock the portals thereof."

The mindset revealed in this can hardly be said to be suitable in scientific and objective inquiry or pursuit of knowledge, but perhaps one should not grudge it in someone from a totally different milieu, if the outcome is positive. One would have thought that with all the commitment and grit the author would have come up with at least a few new things which can be attributed to the Vedas, with solid evidence. This would have made a worthwhile contribution to our understanding of our heritage. Instead, all said and done there is only the author's certificate that "we were agreeably astonished and intensely gratified to find that exceedingly though mathematical problems can be easily and readily solved with the help of these ultra-easy Vedic sutras (or mathematical aphorisms) contained in the Parishishta (the appendix portion) of the Atharva Veda in a few simple steps and by methods which can be conscientiously described as mere 'mental arithmetic' "(paragraph 9 in the preface). That passing reference to the Atharva Veda is all that is ever said by way of source material for the contents. The sutras, incidentally, which appeared later scattered in the book, are short phrases of just about two to four words in Sanskrit, such as Ekadhikena Purvena or Anurupye Shunyam Anyat. (There are 16 of them and in addition there are 13 of what are called sub-sutras, similar in nature to the sutras).

The first key question, which would occur to anyone, is where are these sutras to be found in the Atharva Veda. One does not mean this as a rhetorical question. Considering that at the outset the author seemed set to send all doubting Thomases packing, the least one would expect is that he would point out where the sutras are, say in which part, stanza, page and so on, especially since it is not a small article that is being referred to. Not only has the author not cared to do so, but when Prof.K.S.Shukla, a renowned scholar of ancient Indian mathematics, met him in 1950, when the swamiji visited Lucknow to give a blackboard demonstration of his "Vedic Mathematics", and requested him to point out the sutras in question in the Parishishta of the Atharva Veda, of which he even carried a copy (the standard version edited by G.M.Bolling and J.Von Negelein), the swamiji is said to have told him that the 16 sutra demonstrated by him were not in those Parishishtas and that "they occurred in his own Parishishta and not any other" (Shukla, 1980, or Shukla, 1991). What justification the swamiji thought he had for introducing an appendix in the Atharva Veda, the contents of which are nevertheless to be viewed as from the Veda, is anybody's guess. In any case, even such a Parishishta, written by the swamiji, does not exist in the form of a Sanskrit text.

Let us suppose for a moment that the author indeed found the sutras in some manuscript of the Atharva Veda, which he came across. Would he not then have preserved the manuscript? Would he not have shown at least to some people where the sutras are in the manuscript? Would he not have revealed to some cherished students how to look for sutras with such profound mathematical implications as he attributes to the sutras in question, in that or other manuscripts that may be found? While there is a specific mention in the write-up of Manjula Trivedi, in the beginning of the book, about some 16volume manuscript written by the swamiji having been lost in 1956, there is no mention whatever (let alone any lamentation that would be due in such an event) either in her write-up nor in the swamiji's preface about any original manuscript having been lost. No one certainly has come forward with any information received from the swamiji with regard to the other questions above. It is to be noted that want of time could not be a factor in any of this, since the swamiji kindly informs us in the preface that "Ever since (i.e. since several decades ago), we have been carrying on an incessant and strenuous campaign for the India-wide diffusion of all this scientific knowledge".

The only natural explanation is that there was no such manuscript. It has in fact been mentioned by Agrawala in his general editor's foreword to the book, and also by Manjula Trivedi in the short account of the genesis of the work, included in the book together with a biographical sketch of the swamiji, that the sutras do not appear in hitherto known Parishishtas. The general editor also notes that the style of language of the sutras "point to their discovery by Shri Swamiji himself " (emphasis added); the language style being contemporary can be confirmed independently from other Sanskrit scholars as well. The question why then the contents should be considered 'Vedic' apparently did not bother the general editor, as he agreed with the author that "by definition" the Vedas should contain all knowledge (never mind whether found in the 20th century, or perhaps even later)! Manjula Trivedi, the disciple has of course no problem with the sutras not being found in the Vedas as she in fact says that they were actually reconstructed by her beloved "Gurudeva," on the basis of intuitive revelation from material scattered here and there in the Atharva Veda, after "assiduous research" and 'Tapas' for about eight years in the forests surrounding Shringeri." Isn't that adequate to consider them to be "Vedic"? Well, one can hardly argue with the devout! There is a little problem as to why the Gurudeva himself did not say so (that the sutras were reconstructed) rather than referring to them as sutras contained in the Parishishta of the Atharva Veda, but we will have to let it pass. Anyway the fact remains that she was aware that they could not actually be located in what we lesser mortals consider to be the Atharva Veda. The question of the source of the sutras is merely the first that would come to mind, and already on that there is such a muddle. Actually, even if the sutras were to be found, say in the Atharva Veda or some other ancient text, that still leaves open another fundamental question as to whether they mean or yield, in some cognisable way, what the author claims; in other words, we would still need to know whether such a source really contains the mathematics the swamiji deals with or merely the phrases, may be in some quite different context. It is interesting to consider the swamiji's sutras in this light. One of them, for instance, is Ekadhikena Purvena which literally just means "by one more than the previous one." In chapter I, the swamiji tells us that it is a sutra for finding the digits in the decimal expansion of numbers such as 1/19, and 1/29, where the denominator is a number with 9 in the unit's place; he goes on to give a page-long description of the procedure to be followed, whose only connection with the sutra is that it involves, in particular, repeatedly multiplying by one more than the previous one, namely 2, 3 and so on, respectively, the "previous one" being the number before the unit's place; the full procedure involves a lot more by way of arranging the digits which can in no way be read off from the phrase.

In Chapter II, we are told that the same sutra also means that to find the square of a number like 25 and 35, (with five in unit's place) multiply the number of tens by one more than itself and write 25 ahead of that; like 625, 1,225 and so on. The phrase Ekanyunena Purvena which means "by one less than the previous one" is however given to mean something which has neither to do with decimal expansions nor with squaring of numbers but concerns multiplying together two numbers, one of which has 9 in all places (like 99,999, so on.)!

Allowing oneself such unlimited freedom of interpretation, one can also interpret the same three-word phrase to mean also many other things not only in mathematics but also in many other subjects such as physics, chemistry, biology, economics, sociology and politics. Consider, for instance, the following "meaning": the family size may be allowed to grow, at most, by one more than the previous one. In this we have the family-planning message of the 1960s; the "previous one" being the couple, the prescription is that they should have no more than three children. Thus the lal trikon (red triangle) formula may be seen to be "from the Atharva Veda," thanks to the swamiji's novel technique (with just a bit of credit to yours faithfully). If you think the three children norm now outdated, there is no need to despair. One can get the two-children or even the one-child formula also from the same sutra; count only the man as the "previous one" (the woman is an outsider joining in marriage, isn't she) and in the growth of the family either count only the children or include also the wife, depending on what suits the desired formula!

Another sutra is Yavadunam, which means "as much less;" a lifetime may not suffice to write down all the things such a phrase could "mean," in the spirit as above. There is even a sub-sutra, Vilokanam (observation) and that is supposed to mean various mathematical steps involving observation! In the same vein one can actually suggest a single sutra adequate not only for all of mathematics but many many subjects: Chintanam (think)!

It may be argued that there are, after all, ciphers which convey more information than meets the eye. But the meaning in those cases is either arrived at from the knowledge of the deciphering code or deduced in one or other way using various kinds of contexual information. Neither applies in the present case. The sutras in the swamiji's book are in reality mere names for various steps to be followed in various contexts; the steps themselves had to be known independently. In other words, the mathematical step is not arrived at by understanding or interpreting what are given as sutras; rather, sutras somewhat suggestive of the meaning of the steps are attached to them like names. It is like associating the 'sutra' VIBGYOR to the sequence of colours in rainbow (which make up the white light). Usage of words in Sanskrit, a language which the popular mind unquestioningly associates with the distant past(!), lend the contents a bit of antique finish!

An analysis of the mathematical contents of Tirthaji's book also shows that they cannot be from the Vedas. Though unfortunately there is considerable ignorance about the subject, mathematics from the Vedas is far from being an unexplored area. Painstaking efforts have been made for well over a century to study the original ancient texts from the point of view of understanding the extent of mathematical knowledge in ancient times. For instance, from the study of Vedic Samhitas and Brahamanas it has been noted that they had the system of counting progressing in multiples of 10 as we have today and that they considered remarkably large numbers, even up to 14 digits, unlike other civilizations of those times. From the Vedanga period there is in fact available a significant body of mathematical literature in the form of Shulvasutras, from the period between 800 bc and 500 bc, or perhaps even earlier, some of which contain expositions of various mathematical principles involved in construction of sacrificial 'vedi's needed in performing' yajna's (see, for instance, Sen and Bag 1983).

Baudhyana Shulvasutra, the earliest of the extant Shulvasutras, already contains, for instance, what is currently known as Pythagoras' Theorem (Sen and Bag, 1983, page 78, 1.12). It is the earliest known explicit statement of the theorem in the general form (anywhere in the world) and precedes Pythagoras by at least a few hundred years. The texts also show a remarkable familiarity with many other facts from the so-called Euclidean Geometry and it is clear that considerable use was made of these, long before the Greeks formulated them. The work of George Thibaut in the last century and that of A.Burk around the turn of the century brought to the attention of the world the significance of the mathematics of the Shulvasutras. It has been followed up in this century by both foreign and Indian historians of mathematics. It is this kind of authentic work, and not some mumbo-jumbo that would highlight our rich heritage. I would strongly recommend to the reader to peruse the monograph, The Sulbasutras by S.N.Sen and A.K.Bag (Sen and Bag, 1983), containing the original sutras, their translation and a detailed commentary, which includes a survey of a number of earlier works on the subject. There are also several books on ancient Indian mathematics from the Vedic period.

The contents of the swamiji's book have practically nothing in common with what is known of the mathematics from the Vedic period or even with the subsequent rich tradition of mathematics in India until the advent of the modern era; incidentally, the descriptions of mathematical principles or procedures in ancient mathematical texts are quite explicit and not in terms of cryptic sutras. The very first chapter of the book (as also chapters XXVI to XXVIII) involves the notion of decimal fractions in an essential way. If the contents are to be Vedic, there would have had to be a good deal of familiarity with decimal fractions, even involving several digits, at that time. It turns out that while the Shulvasutras make extensive use of fractions in the usual form, nowhere is there any indication of fractions in decimal form. It is inconceivable that such an important notion would be left out, had it been known, from what are really like users manuals of those times, produced at different times over a prolonged period. Not only the Shulvasutras and the earlier Vedic works, but even the works of mathematicians such as Aryabhata, Brahmagupta and Bhaskara, are not found to contain any decimal fractions. Is it possible that none of them had access to some Vedic source that the swamiji could lay his hands on (and still not describe it specifically)? How far do we have to stretch our credulity?

The fact is that the use of decimal fractions started only in the 16th century, propagated to a large extent by Francois Viete; the use of the decimal point (separating the integer and the fractional parts) itself, as a notation for the decimal representation, began only towards the end of the century and acquired popularity in the 17th century following their use in John Napier's logarithm tables (see, for instance, Boyer, 1968, page 334).

Similarly, in chapter XXII the swamiji claims to give "sutras relevant to successive differentiation, covering the theorems of Leibnitz, Maclaurin, Taylor, etc. and a lot of other material which is yet to be studied and decided on by the great mathematicians of the present-day Western world;" it should perhaps be mentioned before we proceed that the chapter does not really deal with anything of the sort that would even remotely justify such a grandiloquent announcement, but rather deals with differentiation as an operation on polynomials, which is a very special case reducing it all to elementary algebra devoid of the very soul of calculus, as taught even at the college level.

Given the context, we shall leave Leibnitz and company alone, but consider the notions of derivative and successive differentiation. Did the notions exist in the Vedic times? While certain elements preliminary to calculus have been found in the works of Bhaskara II from the 12th century and later Indian mathematicians in the pre-calculus era in international mathematics, such crystallised notions as the derivative or the integral were not known. Though a case may be made that the developments here would have led to the discovery of calculus in India, no historians of Indian mathematics would dream of proposing that they actually had such a notion as the derivative, let alone successive differentiation; the question here is not about performing the operation on polynomials, but of the concept. A similar comment applies with regard to integration, in chapter XXIV. It should also be borne in mind that if calculus were to be known in India in the early times, it would have been acquired by foreigners as well, long before it actually came to be discovered, as there was enough interaction between India and the outside world.

If this is not enough, in Chapter XXXIX we learn that analytic conics has an "important and predominating place for itself in the Vedic system of mathematics," and in Chapter XL we find a whole list of subjects such as dynamics, statics, hydrostatics, pneumatics and applied mathematics listed alongside such elementary things as subtractions, ratios, proportions and such money matters as interest and annuities (!), discounts (!) to which we are assured, without going into details, that the Vedic sutras can be applied. Need we comment any further on this? The remaining chapters are mostly elementary in content, on account of which one does not see such marked incongruities in their respect. It has, however, been pointed out by Shukla that many of the topics considered in the book are alien to the pursuits of ancient Indian mathematicians, not only form the Vedic period but until much later (Shukla, 1989 or Shukla, 1991). These include many such topics as factorisation of algebraic expressions, HCF (highest common factor) of algebraic expressions and various types of simultaneous equations. The contents of the book are akin to much later mathematics, mostly of the kind that appeared in school books of our times or those of the swamiji's youth, and it is unthinkable, in the absence of any pressing evidence, that they go back to the Vedic lore. The book really consists of a compilation of tricks in elementary arithmetic and algebra, to be applied in computations with numbers and polynomials. By a "trick" I do not mean a sleight of hand or something like that; in a general sense a trick is a method or procedure which involves observing and exploring some special features of a situation, which generally tend to be overlooked; for example, the trick described for finding the square of numbers like 15 and 25 with 5 in the unit's place makes crucial use of the fact of 5 being half of 10, the latter being the base in which the numbers are written. Some of the tricks given in the book are quite interesting and admittedly yield quicker solutions than by standard methods (though the comparison made in the book are facetious and misleading). They are of the kind that an intelligent hobbyist experimenting with numbers might be expected to come up with. The tricks are, however, based on well-understood mathematical principles and there is no mystery about them.

Of course to produce such a body of tricks, even using the well-known is still a non-trivial task and there is a serious question of how this came to be accomplished. It is sometimes suggested that Tirthaji himself might have invented the tricks. The fact that he had a M.A.degree in mathematics is notable in this context. It is also possible that he might have learnt some of the tricks from some elders during an early period in his life and developed on them during those "eight years of concentrated contemplation in forest solitude:" this would mean that they do involve a certain element of tradition, though not to the absurd extent that is claimed. These can, however, be viewed only as possibilities and it would not be easy to settle these details. But it is quite clear that the choice is only between alternatives involving only the recent times.

It may be recalled here that there have also been other instances of exposition and propagation of such faster methods of computation applicable in various special situations (without claims of their coming from ancient sources). Trachtenberg's Speed System (see Arther and McShane, 1965) and Lester Meyers' book, High-Speed Mathematics (Meyers, 1947) are some well-known examples of this. Trachtenberg had even set up an Institute in Germany to provide training in high-speed mathematics. While the swamiji's methods are independent of these, for the most part they are similar in spirit.

One may wonder why such methods are not commonly adopted for practical purposes. One main point is that they turn out to be quicker only for certain special classes of examples. For a general example the amount of effort involved (for instance, the count of the individual operations needed to be performed with digits, in arriving at the final answer) is about the same as required by the standard methods; in the swamiji's book, this is often concealed by not writing some of the steps involved, viewing it as "mental arithmetic." Using such methods of fast arithmetic involves the ability or practice to recognize various patterns which would simplify the calculations. Without that, one would actually spend more time, in first trying to recognize patterns and then working by rote anyway, since in most cases it is not easy to find useful patterns.

People who in the course of their work have to do computations as they arise, rather than choose the figures suitably as in the demonstrations, would hardly find it convenient to carry them out by employing umpteen different ways depending on the particular case, as the methods of fast arithmetic involve. It is more convenient to follow the standard method, in which one has only to follow a set procedure to find the answer, even though in some cases this might take more time. Besides, equipment such as calculators and computers have made it unnecessary to tax one's mind with arithmetical computations. Incidentally, the suggestion that this "Vedic Mathematics" of the Shankaracharya could lead to improvement in computers is totally fallacious, since the underlying mathematical principles involved in it were by no means unfamiliar in professional circles.

One of the factors causing people not to pay due attention to the obvious questions about "Vedic Mathematics" seems to be that they are overwhelmed by a sense of wonderment by the tricks. The swamiji tells us in the preface how "the educationists, the cream of the English educated section of the people including highest officials (e.g. the high court judges, the ministers etc.) and the general public as such were all highly impressed; nay thrilled, wonder-struck and flabbergasted!" at his demonstrations of the "Vedic Mathematics." Sometimes one comes across reports about similar thrilling demonstrations by some of the present-day expositors of the subject. Though inevitably they have to be taken with a pinch of salt, I do not entirely doubt the truth of such reports. Since most people have had a difficult time with their arithmetic at school and even those who might have been fairly good would have lost touch, the very fact of someone doing some computations rather fast can make an impressive sight. This effect may be enhanced with well-chosen examples, where some quicker methods are applicable.

Even in the case of general examples where the method employed is not really more efficient than the standard one, the computations might appear to be fast, since the demonstrator would have a lot more practice than the people in the audience. An objective assessment of the methods from the point of view of overall use can only be made by comparing how many individual calculations are involved in working out various general examples, on an average, and in this respect the methods of fast arithmetic do not show any marked advantage which would offset the inconvenience indicated earlier. In any case, it would be irrational to let the element of surprise interfere in judging the issue of origin of "Vedic Mathematics" or create a dreamy and false picture of its providing solutions to all kinds of problems.

It should also be borne in mind that the book really deals only with some middle and high school level mathematics; this is true despite what appear to be chapters dealing with some notions in calculus and coordinate geometry and the mention of a few, little more advanced topics, in the book. The swamiji's claim that "there is no part of mathematics, pure or applied, which is beyond their jurisdiction" is ludicrous. Mathematics actually means a lot more than arithmetic of numbers and algebra of polynomials; in fact multiplying big numbers together, which a lot of people take for mathematics, is hardly something a mathematician of today needs to engage himself in. The mathematics of today concerns a great variety of objects beyond the high school level, involving various kinds of abstract objects generalising numbers, shapes, geometries, measures and so on and several combinations of such structures, various kinds of operations, often involving infinitely many en-tities; this is not the case only about the frontiers of mathematics but a whole lot of it, including many topics applied in physics, engineering, medicine, finance and various other subjects.

Despite all its pretentious verbiage page after page, the swamiji's book offers nothing worthwhile in advanced mathematics whether concretely or by way of insight. Modern mathematics with its multitude of disciplines (group theory, topology, algebraic geometry, harmonic analysis, ergodic theory, combinatorial mathematics-to name just a few) would be a long way from the level of the swamiji's book. There are occasionally reports of some "researchers" applying the swamiji's "Vedic Mathematics" to advanced problems such as Kepler's problem, but such work involves nothing more than tinkering superficially with the topic, in the manner of the swamiji's treatment of calculus, and offers nothing of interest to professionals in the area.

Even at the school level "Vedic Mathematics" deals only with a small part and, more importantly, there too it concerns itself with only one particular aspect, that of faster computation. One of the main aims of mathematics education even at the elementary level consists of developing familiarity with a variety of concepts and their significance. Not only does the approach of "Vedic Mathematics" not contribute anything towards this crucial objective, but in fact might work to its detriment, because of the undue emphasis laid on faster computation. The swamiji's assertion "8 months (or 12 months) at an average rate of 2 or 3 hours per day should suffice for completing the whole course of mathematical studies on these Vedic lines instead of 15 or 20 years required according to the existing systems of the Indian and also foreign universities," is patently absurd and hopefully nobody takes it seriously, even among the activists in the area. It would work as a cruel joke if some people choose to make such a substitution in respect of their children.

It is often claimed that "Vedic Mathematics" is well-appreciated in other countries, and even taught in some schools in UK etc.. In the normal course one would not have the means to examine such claims, especially since few details are generally supplied while making the claims. Thanks to certain special circumstances I came to know a few things about the St. James Independent School, London which I had seen quoted in this context. The School is run by the 'School of Economic Science' which is, according to a letter to me from Mr. James Glover, the Head of Mathematics at the School, "engaged in the practical study of Advaita philosophy". The people who run it have had substantial involvement with religious groups in India over a long period. Thus in essence their adopting "Vedic Mathematics" is much like a school in India run by a religious group adopting it; that school being in London is beside the point. (It may be noted here that while privately run schools in India have limited freedom in choosing their curricula, it is not the case in England). It would be interesting to look into the background and motivation of other institutions about which similar claims are made. At any rate, adoption by institutions abroad is another propaganda feature, like being from ancient source, and should not sway us.

It is not the contention here that the contents of the book are not of any value. Indeed, some of the observations could be used in teaching in schools. They are entertaining and could to some extent enable children to enjoy mathematics. It would, however, be more appropriate to use them as aids in teaching the related concepts, rather than like a series of tricks of magic. Ultimately, it is the understanding that is more important than the transient excitement, By and large, however, such pedagogical application has limited scope and needs to be made with adequate caution, without being carried away by motivated propaganda.

It is shocking to see the extent to which vested interests and persons driven by guided notions are able to exploit the urge for cultural self-assertion felt by the Indian psyche. One would hardly have imagined that a book which is transparently not from any ancient source or of any great mathematical significance would one day be passed off as a storehouse of some ancient mathematical treasure. It is high time saner elements joined hands to educate people on the truth of this so-called Vedic Mathematics and prevent the use of public money and energy on its propagation, beyond the limited extent that may be deserved, lest the intellectual and educational life in the country should get vitiated further and result in wrong attitudes to both history and mathematics, especially in the coming generation.

## 2.2 Neither Vedic Nor Mathematics

We, the undersigned, are deeply concerned by the continuing attempts to thrust the so-called `Vedic Mathematics' on the school curriculum by the NCERT (National Council of Educational Research and Training).

As has been pointed out earlier on several occasions, the so-called 'Vedic Mathematics' is neither 'Vedic' nor can it be dignified by the name of mathematics. 'Vedic Mathematics', as is well-known, originated with a book of the same name by a former Sankaracharya of Puri (the late Jagadguru Swami Shri Bharati Krishna Tirthaji Maharaj) published posthumously in 1965. The book assembled a set of tricks in elementary arithmetic and algebra to be applied in performing computations with numbers and polynomials. As is pointed out even in the foreword to the book by the General Editor, Dr. A.S. Agarwala, the aphorisms in Sanskrit to be found in the book have nothing to do with the Vedas. Nor are these aphorisms to be found in the genuine Vedic literature.

The term "Vedic Mathematics" is therefore entirely misleading and factually incorrect. Further, it is clear from the notation used in the arithmetical tricks in the book that the methods used in this text have nothing to do with the arithmetical techniques of antiquity. Many of the Sanskrit aphorisms in the book are totally cryptic (ancient Indian mathematical writing was anything but cryptic) and often so generalize to be devoid of any specific mathematical meaning. There are several authoritative texts on the mathematics of Vedic times that could be used in part to teach an authoritative and correct account of ancient Indian mathematics but this book clearly cannot be used for any such purpose. The teaching of mathematics involves both the teaching of the basic concepts of the subject as well as methods of mathematical computation. The so-called "Vedic Mathematics" is entirely inadequate to this task considering that it is largely made up of tricks to do some elementary arithmetic computations. Many of these can be far more easily performed on a simple computer or even an advanced calculator.

The book "Vedic Mathematics" essentially deals with arithmetic of the middle and high-school level. Its claims that "there is no part of mathematics, pure or applied, which is beyond their jurisdiction" is simply ridiculous. In an era when the content of mathematics teaching has to be carefully designed to keep pace with the general explosion of knowledge and the needs of other modern professions that use mathematical techniques, the imposition of "Vedic Mathematics" will be nothing short of calamitous.

India today has active and excellent schools of research and teaching in mathematics that are at the forefront of modern research in their discipline with some of them recognised as being among the best in the world in their fields of research. It is noteworthy that they have cherished the legacy of distinguished Indian mathematicians like Srinivasa Ramanujam, V. K. Patodi, S. Minakshisundaram, Harish Chandra, K. G. Ramanathan, Hansraj Gupta, Syamdas Mukhopadhyay, Ganesh Prasad, and many others including several living Indian mathematicians. But not one of these schools has lent an iota of legitimacy to 'Vedic Mathematics'. Nowhere in the world does any school system teach "Vedic Mathematics" or any form of ancient mathematics for that matter as an adjunct to modern mathematical teaching. The bulk of such teaching belongs properly to the teaching of history and in particular the teaching of the history of the sciences.

We consider the imposition of 'Vedic Mathematics' by a Government agency, as the perpetration of a fraud on our children, condemning particularly those dependent on public education to a sub-standard mathematical education. Even if we assumed that those who sought to impose 'Vedic Mathematics' did so in good faith, it would have been appropriate that the NCERT seek the assistance of renowned Indian mathematicians to evaluate so-called "Vedic Mathematics" before making it part of the National Curricular framework for School Education. Appallingly they have not done so. In this context we demand that the NCERT submit the proposal for the introduction of 'Vedic Mathematics' in the school curriculum to recognized bodies of mathematical experts in India, in particular the National Board of Higher Mathematics (under the Dept. of Atomic Energy), and the Mathematics sections of the Indian Academy of Sciences and the Indian National Science Academy, for a thorough and critical examination. In the meanwhile no attempt should be made to thrust the subject into the school curriculum either through the centrally administered school system or by trying to impose it on the school systems of various States.

We are concerned that the essential thrust behind the campaign to introduce the so-called 'Vedic Mathematics' has more to do with promoting a particular brand of religious majoritarianism and associated obscurantist ideas rather than any serious and meaningful development of mathematics teaching in India. We note that similar concerns have been expressed about other aspects too of the National Curricular Framework for School Education. We re-iterate our firm conviction that all teaching and pedagogy, not just the teaching of mathematics, must be founded on rational, scientific and secular principles.

[Many eminent scholars, researchers from renowned Indian foreign universities have signed this. See the end of section for a detailed list.]

We now give the article "Stop this Fraud on our Children!" from Peoples Democracy.

Over a hundred leading scientists, academicians, teachers and educationists, in a statement have protested against the attempts by the Vajpayee government to introduce Vedic Mathematics and Vedic Astrology courses in the education system. They have in one voice demanded "Stop this Fraud on our Children!"

The scientists and mathematicians are deeply concerned that the essential thrust behind the campaign to introduce the so-called 'Vedic Mathematics' in the school curriculum by the NCERT, and 'Vedic Astrology' at the university level by the University Grants Commission, has more to do with promoting a particular brand of religious majoritarianism and associated obscurantist ideas than with any serious development of mathematical or scientific teaching in India. In rejecting these attempts, they re-iterate their firm conviction that all teaching and pedagogy must be founded on rational, scientific and secular principles.

Pointing out that the so-called *"Vedic Mathematics"* is neither vedic nor mathematics, they say that the imposition of '*Vedic maths'* will condemn particularly those dependent on public education to a sub-standard mathematical education and will be calamitous for them.

"The teaching of mathematics involves both imparting the basic concepts of the subject as well as methods of mathematical computations. The so-called 'Vedic maths' is entirely inadequate to this task since it is largely made up of tricks to do some elementary arithmetic computations. Its value is at best recreational and its pedagogical use limited", the statement noted. The signatories demanded that the NCERT submit the proposal for the introduction of 'Vedic maths' in the school curriculum for a thorough and critical examination to any of the recognised bodies of mathematical experts in India.

Similarly, they assert that while many people may believe in astrology, this is in the realm of belief and is best left as part of personal faith. Acts of faith cannot be confused with the study and practice of science in the public sphere.

Signatories to the statement include award -winning scientists, Fellows of the Indian National Science Academy, the Indian Academy of Sciences, Senior Professors and eminent mathematicians. Prominent among the over 100 scientists who have signed the statement are:

1. Yashpal (Professor, Eminent Space Scientist, Former Chairman, UGC),
2. J.V.Narlikar (Director, Inter University Centre for Astronomy and Astrophysics, Pune)
3. M.S.Raghunathan (Professor of Eminence, School of Maths, TIFR and Chairman National Board for Higher Maths).
4. S G Dani, (Senior Professor, School of Mathematics, TIFR)
5. R Parthasarathy (Senior Professor, School of Mathematics, TIFR),
6. Alladi Sitaram (Professor, Indian Statistical Institute (ISI), Bangalore),
7. Vishwambar Pati (Professor, Indian Statistical Institute , Bangalore),
8. Kapil Paranjape (Professor, Institute of Mathematical Sciences (IMSc), Chennai),
9. S Balachandra Rao, (Principal and Professor of Maths, National College, Bangalore)
10. A P Balachandran, (Professor, Dept. of Physics, Syracuse University USA),
11. Indranil Biswas (Professor, School of Maths, TIFR)
12. C Musili (Professor, Dept. of Maths and Statistics, Univ. of Hyderabad),
13. V.S.Borkar (Prof., School of Tech. and Computer Sci., TIFR)

14. Madhav Deshpande (Prof. of Sanskrit and Linguistics, Dept. of Asian Languages and Culture, Univ. of Michigan, USA),
15. N. D. Haridass (Senior Professor, Institute of Mathematical Science, Chennai),
16. V.S. Sunder (Professor, Institute of Mathematical Sciences, Chennai),
17. Nitin Nitsure (Professor, School of Maths, TIFR),
18. T Jayaraman (Professor, Institute of Mathematical Sciences, Chennai),
19. Vikram Mehta (Professor, School of Maths, TIFR),
20. R. Parimala (Senior Professor, School of Maths, TIFR),
21. Rajat Tandon (Professor and Head, Dept. of Maths and Statistics, Univ. of Hyderabad),
22. Jayashree Ramdas (Senior Reseacrh Scientist, Homi Bhabha Centre for Science Education, TIFR) ,
23. Ramakrishna Ramaswamy (Professor, School of Physical Sciences, JNU), D P Sengupta (Retd. Prof. IISc., Bangalore),
24. V Vasanthi Devi (Former VC, Manonmaniam Sundaranar Univ. Tirunelveli),
25. J K Verma (Professor, Dept. of Maths, IIT Bombay),
26. Bhanu Pratap Das (Professor, Indian Institute of Astrophysics, Bangalore)
27. Pravin Fatnani (Head, Accelerator Controls Centre, Centre for Advanced Technology, Indore),
28. S.L. Yadava (Professor, TIFR Centre, IISc, Bangalore) ,
29. Kumaresan, S (Professor, Dept. of Mathematics, Univ. of Mumbai),
30. Rahul Roy (Professor, ISI ,Delhi)

and others….

## 2.3 Views about the Book in Favour and Against

The view of his Disciple Manjula Trivedi, Honorary General Secretary, Sri Vishwa Punarnirmana Sangha, Nagpur written on 16[th] March 1965 and published in a reprint and revised edition of the book on Vedic Mathematics reads as follows.

"I now proceed to give a short account of the genesis of the work published here. Revered Guruji used to say that he had reconstructed the sixteen mathematical formulae (given in this text) from the Atharveda after assiduous research and 'Tapas' for about eight years in the forests surrounding Sringeri. *Obviously these formulae are not to be found in the present recensions of Atharvaveda; they were actually reconstructed, on the basis of intuitive revelation, from materials scattered here and there in the Atharvaveda.* Revered Gurudeva used to say that he had written sixteen volumes on these sutras one for each sutra and that the manuscripts of the said volumes were deposited at the house of one of his disciples. Unfortunately the said manuscripts were lost irretrievably from the place of their deposit and this colossal loss was finally confirmed in 1956.

Revered Gurudeva was not much perturbed over this irretrievable loss and used to say that everything was there in his memory and that he would rewrite the 16 volumes!

In 1957, when he had decided finally to undertake a tour of the USA he rewrote from memory the present volume giving an introductory account of the sixteen formulae reconstructed by him …. The present volume is the only work on mathematics that has been left over by Revered Guruji.

The typescript of the present volume was left over by Revered Gurudeva in USA in 1958 for publication. He had been given to understand that he would have to go to the USA for correction of proofs and personal supervision of printing. But his health deteriorated after his return to India and finally the typescript was brought back from the USA after his attainment of Mahasamadhi in 1960."

A brief sketch from the Statesman, India dated 10[th] Jan 1956 read as follows. "Sri Shankaracharya denies any spiritual or miraculous powers giving the credit for his revolutionary knowledge to anonymous ancients, who in 16 sutras and 120 words laid down simple formulae for all the world's mathematical problems […]. I could read a short descriptive note he had prepared on, "The Astounding Wonders of Ancient Indian Vedic Mathematics". His Holiness, it appears, had spent years in contemplation, and while going through the Vedas had suddenly happened upon the key to what many historians, devotees and translators had dismissed as meaningless jargon. There, contained in certain Sutras, were the processes of mathematics, psychology, ethics and metaphysics.

"During the reign of King Kamsa" read a sutra, "rebellions, arson, famines and insanitary conditions prevailed". Decoded this little piece of libelous history gave decimal answer to the fraction 1/17, sixteen processes of simple mathematics reduced to one.

The discovery of one key led to another, and His Holiness found himself turning more and more to the astounding knowledge contained in words whose real meaning had been lost to humanity for generations. This loss is obviously one of the greatest mankind has suffered and I suspect, resulted from the secret being entrusted to people like myself, to whom a square root is one of life's perpetual mysteries. Had it survived, every – educated 'soul' would be a mathematical 'wizard' and maths 'masters' would "starve". For my note reads "Little children merely look at the sums written on the blackboard and immediately shout out the answers they have … [Pages 353-355 Vedic Mathematics]

We now briefly quote the views of S.C. Sharma, Ex Head of the Department of Mathematics, NCERT given in Mathematics Today, September 1986.

"The epoch-making and monumental work on Vedic Mathematics unfolds a new method of approach. It relates to the truth of numbers and magnitudes equally applicable to all sciences and arts.

The book brings to light how great and true knowledge is born of intuition, quite different from modern western method. The ancient Indian method and its secret techniques are examined and shown to be capable of solving various problems of mathematics. The universe we live in has a basic mathematical structure obeying the rules of mathematical measures and relations. All the subjects in mathematics – Multiplication, Division, Factorization Equations of calculus Analytical Conics etc. are dealt with in forty chapters vividly working out all problems, in the easiest ever method discovered so far. The volume more a magic is the result of institutional visualization of fundamental mathematical truths born after eight years of highly concentrated endeavor of Jagadguru Sri Bharati Krishna Tirtha.

Throughout this book efforts have been made to solve the problems in a short time and in short space also …, one can see that the formulae given by the author from Vedas are very interesting and encourage a young mind for learning mathematics as it will not be a bugbear to him".

This writing finds its place in the back cover of the book of Vedic Mathematics of Jagadguru. Now we give the views of Bibek Debroy, "The fundamentals of Vedic Mathematics" pp. 126-127 of Vedic Mathematics in Tamil volume II).

"Though Vedic Mathematics evokes Hindutva connotations, the fact is, it is a system of simple arithmetic, which can be used for intricate calculations.

The resurgence of interest in Vedic Mathematics came about as a result of Jagadguru Swami Sri Bharati Krishna Tirthaji Maharaj publishing a book on the subject in 1965. Then recently the erstwhile Bharatiya Janata Party governments in Uttar Pradesh, Madhya Pradesh and Himachal Pradesh introduced Vedic Mathematics into the school syllabus, but this move was perceived as an attempt to impose Hindutva, because Vedic philosophy was being projected as the repository of all human wisdom. The subsequent hue and cry over the teaching of Vedic Mathematics is mainly because it has come to be identified with, fundamentalism and obscurantism, both considered poles opposite of science. The critics argue that belief in Vedic Mathematics automatically necessitates belief in Hindu renaissance. But Tirtha is not without his critics, even apart from those who consider Vedic maths is "unscientific".

## 2.4 Vedas: Repositories of Ancient Indian Lore

Extent texts of the Vedas do not contain mathematical formulae but they have been found in later associated works. Jagadguru the author of Vedic Mathematics says he has discovered 16 mathematical formulae, …

A standard criticism is that the Vedic Mathematics text is limited to middle and high school formulations and the emphasis is on a series of problem solving tricks. The critics also point out that the Atharva Veda appendix containing

Tirtha's 16 mathematical formulae, is not to be found in any part of the existing texts. A third criticism is the most pertinent. The book is badly written. (p.127, Vedic Mathematics 2) [85]. We shall now quote the preface given by His Excellency Dr. L.M.Singhvi, High Commissioner for India in UK, given in pp. V to VI Reprint Vedic Mathematics 2005, Book 2, [51].

Vedic Mathematics for schools is an exceptional book. It is not only a sophisticated pedagogic tool but also an introduction to an ancient civilization. It takes us back to many millennia of India's mathematical heritage…

The real contribution of this book, "Vedic Mathematics for schools, is to demonstrate that Vedic Mathematics belongs not only to an hoary antiquity but is any day as modern as the day after tomorrow. What distinguishes it particularly is that it has been fashioned by British teachers for use at St.James independent schools in London and other British schools and that it takes its inspiration from the pioneering work of the late Sankaracharya of Puri…

Vedic Mathematics was traditionally taught through aphorisms or Sutras. A sutra is a thread of knowledge, a theorem, a ground norm, a repository of proof. It is formulated as a proposition to encapsulate a rule or a principle. Both Vedic Mathematics and Sanskrit grammar built on the foundations of rigorous logic and on a deep understanding of how the human mind works. The methodology of Vedic Mathematics and of Sanskrit grammar help to hone the human intellect and to guide and groom the human mind into modes of logical reasoning."

## 2.5 A Rational Approach to Study Ancient Literature

Excerpted from Current Science Vol. 87, No. 4, 25 Aug. 2004.
It was interesting to read about Hertzstark's hand-held mechanical calculator, which converted subtraction into addition. But I would like to comment on the 'Vedic Mathematics' referred to in the note. Bharati Krishna Tirtha is a good mathematician, but the term 'Vedic Mathematics' coined by him is misleading, because his mathematics has nothing to do with the *Vedas*. It is his 20th century invention, which should be called 'rapid mathematics' or 'Shighra Ganita'. He has disguised his intention of giving it an aura of discovering ancient knowledge with the following admission in the foreword of his book, which few people take the trouble to read. He says there that he saw (thought of) of his Sutras just like the Vedic Rishis saw (thought of) the Richas. That is why he has called his method 'Vedic Mathematics'. This has made it attractive to the ignorant and not-so ignorant public. I hope scientists will take note of this fact. Vedic astrology is another term, which fascinates people and captures their imagination about its ancient origin. Actually, there is no mention of horoscope and planetary influence in Vedic literature. It only talks of Tithis and Nakshatras as astronomical entities useful for devising a calendar controlled by a series of sacrifices. Astrology of planets originated in Babylon, where astronomers made regular observations of planets, but could not understand their complicated motions. Astrology spread from there to Greece and Europe in the west and to India in the east. There is nothing Vedic about it. It appears that some Indian intellectuals would use the word Vedic as a brand name to sell their ideas to the public. It is imperative that scientists should study ancient literature from a rational point of view, consistent with the then contemporary knowledge."

## 2.6 Shanghai Rankings and Indian Universities

This article is from Current Science Vol. 87, No. 4, 25 August 2004 [7].
"The editorial "The Shanghai Ranking" is a shocking revelation about the fate of higher education and a slide down of scientific research in India. None of the reputed '5 star' Indian universities qualifies to find a slot among the top 500 at the global level. IISc Bangalore and IITs at Delhi and Kharagpur provide some redeeming feature and put India on the score board with a rank between 250 and 500. Some of the interesting features of the Shanghai rankings are noteworthy: (i) Among the top 99 in the world, we have universities from USA (58), Europe (29), Canada (4), Japan (5), Australia (2) and Israel (1). (ii) On the Asia-Pacific list of top 90, we have maximum number of universities from Japan (35), followed by China (18) including Taiwan (5) and Hongkong (5), Australia (13), South Korea (8), Israel (6), India (3), New Zealand (3), Singapore (2) and Turkey (2). (iii) Indian universities lag behind even small Asian countries, viz. South Korea, Israel, Taiwan and Hongkong, in ranking. I agree with the remark, 'Sadly, the real universities in India are limping, with the faculty disinterested in research outnumbering those with an academic bent of mind'. The malaise is deep rooted and needs a complete overhaul of the Indian education system."

## 2.7 Conclusions derived on Vedic Mathematics and the Calculations of Guru Tirthaji - Secrets of Ancient Maths

This article was translated and revised by its author Jan Hogendijk from his original version published in Dutch in the *Nieuwe Wiskrant* vol. 23 no.3 (March 2004), pp. 49–52.

"The "Vedic" methods of mental calculations in the decimal system are all based on the book Vedic Mathematics by Jagadguru (world guru) Swami (monk) Sri (reverend) Bharati Krsna Tirthaji Maharaja, which appeared in 1965 and which has been reprinted many times [51].

The book contains sixteen brief sutras that can be used for mental calculations in the decimal place-value system. An example is the sutra Ekadhikena Purvena, meaning: by one more than the previous one. The Guru explains that this sutra can for example be used in the mental computation of the period of a recurring decimal fraction such as 1/19 = 0.052631578947368421. as follows:

The word "Vedic" in the title of the book suggests that these calculations are authentic Vedic Mathematics. The question now arises how the Vedic mathematicians were able to write the recurrent decimal fraction of 1/19, while decimal fractions were unknown in India before the seventeenth century. We will first investigate the origin of the sixteen sutras. We cite the Guru himself [51]:

"And the contemptuous or, at best, patronizing attitude adopted by some so-called orientalists, indologists, antiquarians, research-scholars etc. who condemned, or light heartedly, nay irresponsibly, frivolously and flippantly dismissed, several abstruse-looking and recondite parts of the Vedas as 'sheer nonsense' or as 'infant-humanity's prattle,' and so on … further confirmed and strengthened our resolute determination to unravel the too-long hidden mysteries of philosophy and science contained in ancient India's Vedic lore, with the consequence that, after eight years of concentrated contemplation in forest-solitude, we were at long last able to recover the long lost keys which alone could unlock the portals thereof.

"And we were agreeably astonished and intensely gratified to find that exceedingly tough mathematical problems (which the mathematically most advanced present day Western scientific world had spent huge lots of time, energy and money on and which even now it solves with the utmost difficulty and after vast labour involving large numbers of difficult, tedious and cumbersome 'steps' of working) can be easily and readily solved with the help of these ultra-easy Vedic Sutras (or mathematical aphorisms) contained in the Parisısta (the Appendix-portion) of the Atharvaveda in a few simple steps and by methods which can be conscientiously described as mere 'mental arithmetic.' "

Concerning the applicability of the sixteen sutras to all mathematics, we can consult the Foreword to Vedic Mathematics written by Swami Pratyagatmananda Saraswati. This Swami states that one of the sixteen sutras reads Calanakalana, which can be translated as Becoming. The Guru himself translates the sutra in question as "differential calculus"[4, p. 186]. Using this "translation" the sutra indeed promises applicability to a large area in mathematics; but the sutra is of no help in differentiating or integrating a given function such as $f(x) = 1/\sin x$.

Sceptics have tried to locate the sutras in the extant Parisista's (appendices) of the Atharva-Veda, one of the four Vedas. However, the sutras have never been found in authentic texts of the Vedic period. It turns out that the Guru had "seen" the sutras by himself, just as the authentic Vedas were, according to tradition, "seen" by the great Rishi's or seers of ancient India. The Guru told his devotees that he had "re-constructed" his sixteen sutras from the Atharva-Veda in the eight years in which he lived in the forest and spent his time on contemplation and ascetic practices. The book Vedic Mathematics is introduced by a General Editor's Note [51], in which the following is stated about the sixteen sutras: "[the] style of language also points to their discovery by Sri Swamiji (the Guru)himself."

Now we know enough about the authentic Katapayadi system to identify the origin of the Guru's verse about π / 10. Here is the verse: (it should be noted that the abbreviation r represents a vowel in Sanskrit):

gopi bhagya madhuvrata
srngiso dadhi sandhiga
Khala jivita Khatava
Gala hala rasandhara.

According to the guru, decoding the verse produces the following number:

31415 92653 58979 32384 62643 38327 92

In this number we recognize the first 31 decimals of π (the 32th decimal of π is 5). In the authentic Katapayadi system, the decimals are encoded in reverse order. So according to the authentic system, the verse is decoded as

29723 83346 26483 23979 85356 29514 13

We conclude that the verse is not medieval, and certainly not Vedic. In all likelihood, the guru is the author of the verse.

There is nothing intrinsically wrong with easy methods of mental calculations and mnemonic verses for π. However, it was a miscalculation on the part of the Guru to present his work as ancient Vedic lore. Many experts in India know that the relations between the Guru's methods and the Vedas are faked. In 1991 the supposed "Vedic" methods of mental calculation were introduced in schools in some cities, perhaps in the context of the political program of saffronisation, which emphasizes Hindu religious elements in society (named after the saffron garments of Hindu Swamis). After many protests, the "Vedic" methods were omitted from the programs, only to be reintroduced a few years later. In 2001, a group of intellectuals in India published a statement against the introduction of the Guru's "Vedic" mathematics in primary schools in India.

Of course, there are plenty of real highlights in the ancient and medieval mathematical tradition of India. Examples are the real Vedic sutras that we have quoted in the beginning of this paper; the decimal place-value system for integers; the concept of sine; the cyclic method for finding integer solutions x, y of the "equation of Pell" in the form $px^2 + 1 = y^2$ (for pa given integer); approximation methods for the sine and arctangents equivalent to modern Taylor series expansions; and so on. Compared to these genuine contributions, the Guru's mental calculation are of very little interest. In the same way, the Indian philosophical tradition has a very high intrinsic value, which does not need to be "proved" by the so-called applications invented by Guru Tirthaji.

**Chapter Three**

# INTRODUCTION TO BASIC CONCEPTS AND A NEW FUZZY MODEL

In this chapter we briefly the recall the mathematical models used in the chapter IV for analysis of, "Is Vedic Mathematics – vedas or mathematics?"; so as to make the book a self contained one. Also in this chapter we have introduced two new models called as new fuzzy dynamical system and new neutrosophic dynamical model to analyze the problem. This chapter has six sections. Section One just recalls the working of the Fuzzy Cognitive Maps (FCMs) model. Definition and illustration of the Fuzzy Relational Maps (FRMs) model is carried out in section two. Section three introduces the new fuzzy dynamical system. In section 4 we just recall the definition of Neutrosophic Cognitive Maps (NCMs), Neutrosophic Relational Maps (NRMs) are given in section 5 (for more about these notions please refer [143]). The final section for the first time introduces the new neutrosophic dynamical model, which can at a time analyze multi experts (n experts, n any positive integer) opinion using a single fuzzy neutrosophic matrix.

## 3.1 Introduction to FCM and the Working of this Model

In this section we recall the notion of Fuzzy Cognitive Maps (FCMs), which was introduced by Bart Kosko [68] in the year 1986. We also give several of its interrelated definitions. FCMs have a major role to play mainly when the data concerned is an unsupervised one. Further this method is most simple and an effective one as it can analyse the data by directed graphs and connection matrices.

**DEFINITION 3.1.1:** *An FCM is a directed graph with concepts like policies, events etc. as nodes and causalities as edges. It represents causal relationship between concepts.*

*Example 3.1.1*: In Tamil Nadu (a southern state in India) in the last decade several new engineering colleges have been approved and started. The resultant increase in the production of engineering graduates in these years is disproportionate with the need of engineering graduates. This has resulted in thousands of unemployed and underemployed graduate engineers. Using an expert's opinion we study the effect of such unemployed people on the society. An expert spells out the five major concepts relating to the unemployed graduated engineers as

$E_1$ – Frustration
$E_2$ – Unemployment
$E_3$ – Increase of educated criminals
$E_4$ – Under employment
$E_5$ – Taking up drugs etc.

The directed graph where $E_1, \ldots, E_5$ are taken as the nodes and causalities as edges as given by an expert is given in the following Figure 3.1.1:

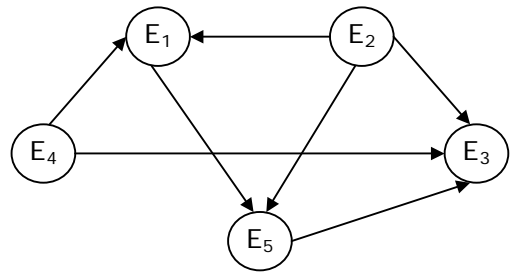

**FIGURE: 3.1.1**

According to this expert, increase in unemployment increases frustration. Increase in unemployment, increases the educated criminals. Frustration increases the graduates to take up to evils like drugs etc. Unemployment also leads to the increase in number of persons who take up to drugs, drinks etc. to forget their worries and unoccupied time. Under-employment forces them to do criminal acts like theft (leading to murder) for want of more money and so on. Thus one cannot actually get data for this but can use the expert's opinion for this unsupervised data to obtain some idea about the real plight of the situation. This is just an illustration to show how FCM is described by a directed graph.

{If increase (or decrease) in one concept leads to increase (or decrease) in another, then we give the value 1. If there exists no relation between two concepts the value 0 is given. If increase (or decrease) in one concept decreases (or increases) another, then we give the value –1. Thus FCMs are described in this way.}

**DEFINITION 3.1.2:** *When the nodes of the FCM are fuzzy sets then they are called as fuzzy nodes.*

**DEFINITION 3.1.3:** *FCMs with edge weights or causalities from the set {–1, 0, 1} are called simple FCMs.*

**DEFINITION 3.1.4:** *Consider the nodes / concepts $C_1, ..., C_n$ of the FCM. Suppose the directed graph is drawn using edge weight $e_{ij} \in \{0, 1, -1\}$. The matrix E be defined by $E = (e_{ij})$ where $e_{ij}$ is the weight of the directed edge $C_i C_j$. E is called the adjacency matrix of the FCM, also known as the connection matrix of the FCM.*

It is important to note that all matrices associated with an FCM are always square matrices with diagonal entries as zero.

**DEFINITION 3.1.5:** *Let $C_1, C_2, ..., C_n$ be the nodes of an FCM. $A = (a_1, a_2, ..., a_n)$ where $a_i \in \{0, 1\}$. A is called the instantaneous state vector and it denotes the on-off position of the node at an instant.*

$$a_i = 0 \text{ if } a_i \text{ is off and}$$
$$a_i = 1 \text{ if } a_i \text{ is on for } i = 1, 2, ..., n.$$

**DEFINITION 3.1.6:** *Let $C_1, C_2, ..., C_n$ be the nodes of an FCM. Let $\overrightarrow{C_1C_2}, \overrightarrow{C_2C_3}, \overrightarrow{C_3C_4}, ..., \overrightarrow{C_iC_j}$ be the edges of the FCM ($i \neq j$). Then the edges form a directed cycle. An FCM is said to be cyclic if it possesses a directed cycle. An FCM is said to be acyclic if it does not possess any directed cycle.*

**DEFINITION 3.1.7:** *An FCM with cycles is said to have a feedback.*

**DEFINITION 3.1.8:** *When there is a feedback in an FCM, i.e., when the causal relations flow through a cycle in a revolutionary way, the FCM is called a dynamical system.*

**DEFINITION 3.1.9:** *Let $\overrightarrow{C_1C_2}, \overrightarrow{C_2C_3}, ..., \overrightarrow{C_{n-1}C_n}$ be a cycle. When $C_i$ is switched on and if the causality flows through the edges of a cycle and if it again causes $C_i$, we say that the dynamical system goes round and round. This is true for any node $C_i$, for $i = 1, 2, ..., n$. The equilibrium state for this dynamical system is called the hidden pattern.*

**DEFINITION 3.1.10:** *If the equilibrium state of a dynamical system is a unique state vector, then it is called a fixed point.*

***Example 3.1.2:*** Consider a FCM with $C_1, C_2, ..., C_n$ as nodes. For example let us start the dynamical system by switching on $C_1$. Let us assume that the FCM settles down with $C_1$ and $C_n$ on i.e. the state vector remains as $(1, 0, 0, ..., 0, 1)$ this state vector $(1, 0, 0, ..., 0, 1)$ is called the fixed point.

**DEFINITION 3.1.11:** *If the FCM settles down with a state vector repeating in the form $A_1 \rightarrow A_2 \rightarrow ... \rightarrow A_i \rightarrow A_1$ then this equilibrium is called a limit cycle.*

Methods of finding the hidden pattern are discussed in the following.

**DEFINITION 3.1.12:** *Finite number of FCMs can be combined together to produce the joint effect of all the FCMs. Let $E_1$, $E_2$, ... , $E_p$ be the adjacency matrices of the FCMs with nodes $C_1$, $C_2$, ..., $C_n$ then the combined FCM is got by adding all the adjacency matrices $E_1$, $E_2$, ..., $E_p$.*

*We denote the combined FCM adjacency matrix by $E = E_1 + E_2 + ... + E_p$.*

**NOTATION:** Suppose $A = (a_1, ..., a_n)$ is a vector which is passed into a dynamical system E. Then $AE = (a'_1, ..., a'_n)$ after thresholding and updating the vector suppose we get $(b_1, ..., b_n)$ we denote that by

$$(a'_1, a'_2, ..., a'_n) \to (b_1, b_2, ..., b_n).$$

Thus the symbol '$\to$' means the resultant vector has been thresholded and updated.

FCMs have several advantages as well as some disadvantages. The main advantage of this method is; it is simple. It functions on expert's opinion. When the data happens to be an unsupervised one the FCM comes handy. This is the only known fuzzy technique that gives the hidden pattern of the situation. As we have a very well known theory, which states that the strength of the data depends on, the number of experts' opinion we can use combined FCMs with several experts' opinions.

At the same time the disadvantage of the combined FCM is when the weightages are 1 and –1 for the same $C_i$ $C_j$, we have the sum adding to zero thus at all times the connection matrices $E_1, ..., E_k$ may not be conformable for addition.

Combined conflicting opinions tend to cancel out and assisted by the strong law of large numbers, a consensus emerges as the sample opinion approximates the underlying population opinion. This problem will be easily overcome if the FCM entries are only 0 and 1.

We have just briefly recalled the definitions. For more about FCMs please refer Kosko [68]. Fuzzy Cognitive Maps (FCMs) are more applicable when the data in the first place is an unsupervised one. The FCMs work on the opinion of experts. FCMs model the world as a collection of classes and causal relations between classes. FCMs are fuzzy signed directed graphs with feedback. The directed edge $e_{ij}$ from causal concept $C_i$ to concept $C_j$ measures how much $C_i$ causes $C_j$. The time varying concept function $C_i(t)$ measures the non negative occurrence of some fuzzy event, perhaps the strength of a political sentiment, historical trend or military objective.

FCMs are used to model several types of problems varying from gastric-appetite behavior, popular political developments etc. FCMs are also used to model in robotics like plant control.

The edges $e_{ij}$ take values in the fuzzy causal interval $[-1, 1]$. $e_{ij} = 0$ indicates no causality, $e_{ij} > 0$ indicates causal increase $C_j$ increases as $C_i$ increases (or $C_j$ decreases as $C_i$ decreases). $e_{ij} < 0$ indicates causal decrease or negative causality. $C_j$ decreases as $C_i$ increases (and or $C_j$ increases as $C_i$ decreases). Simple FCMs have edge values in $\{-1, 0, 1\}$. Then if causality occurs, it occurs to a maximal positive or negative degree. Simple FCMs provide a quick first approximation to an expert stand or printed causal knowledge.

***Example 3.1.3:*** We illustrate this by the following, which gives a simple FCM of a Socio-economic model. A Socio-economic model is constructed with Population, Crime, Economic condition, Poverty and Unemployment as nodes or concept. Here the simple trivalent directed graph is given by the following Figure 3.1.2, which is the experts opinion.

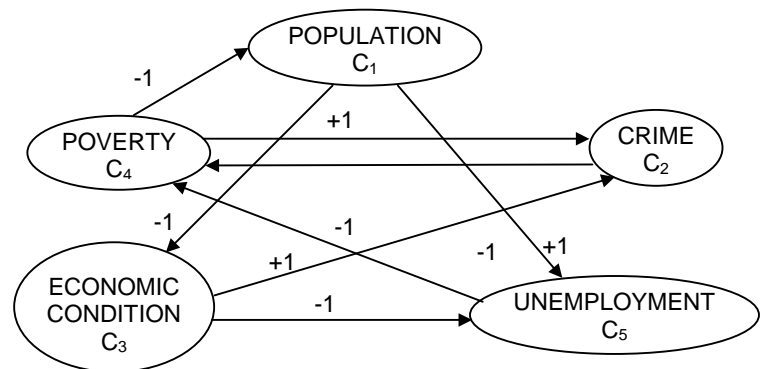

**FIGURE: 3.1.2**

Causal feedback loops abound in FCMs in thick tangles. Feedback precludes the graph-search techniques used in artificial-intelligence expert systems.

FCMs feedback allows experts to freely draw causal pictures of their problems and allows causal adaptation laws, infer causal links from simple data. FCM feedback forces us to abandon graph search, forward and especially backward chaining. Instead we view the FCM as a dynamical system and take its equilibrium behavior as a forward-evolved inference. Synchronous FCMs behave as Temporal Associative Memories (TAM). We can always, in case of a model, add two or more FCMs to produce a new FCM. The strong law of large numbers ensures in some sense that knowledge reliability increases with expert sample size.

We reason with FCMs. We pass state vectors C repeatedly through the FCM connection matrix E, thresholding or non-linearly transforming the result after each pass. Independent of the FCMs size, it quickly settles down to a temporal associative memory limit cycle or fixed point which is the hidden pattern of the system for that state vector C. The limit cycle or fixed-point inference summarizes the joint effects of all the interacting fuzzy knowledge.

Consider the 5 × 5 causal connection matrix E that represents the socio economic model using FCM given in figure in Figure 3.1.2.

$$E = \begin{bmatrix} 0 & 0 & -1 & 0 & 1 \\ 0 & 0 & 0 & -1 & 0 \\ 0 & -1 & 0 & 0 & -1 \\ -1 & 1 & 0 & 0 & 0 \\ 0 & 0 & 0 & 1 & 0 \end{bmatrix}$$

Concept nodes can represent processes, events, values or policies. Consider the first node $C_1 = 1$. We hold or clamp $C_1$ on the temporal associative memories recall process. Threshold signal functions synchronously update each concept after each pass, through the connection matrix E. We start with the concept population alone in the ON state, i.e., $C_1 = (1\ 0\ 0\ 0\ 0)$. The arrow indicates the threshold operation,

$$
\begin{aligned}
C_1 E &= (0\ 0\ -1\ 0\ 1) \rightarrow (1\ 0\ 0\ 0\ 1) \\
&= C_2 \\
C_2 E &= (0\ 0\ -1\ 1\ 1) \rightarrow (1\ 0\ 0\ 1\ 1) \\
&= C_3 \\
C_3 E &= (-1\ 1\ -1\ 1\ 1) \rightarrow (1\ 1\ 0\ 1\ 1) \\
&= C_4 \\
C_4 E &= (-1\ 1\ -1\ 0\ 1) \rightarrow (1\ 1\ 0\ 0\ 1) \\
&= C_5 \\
C_5 E &= (0\ 0\ -1\ 0\ 1) \rightarrow (1\ 0\ 0\ 0\ 1) \\
&= C_6 = C_2.
\end{aligned}
$$

So the increase in population results in the unemployment problem, which is a limit cycle. For more about FCM refer Kosko [67] and for more about these types of socio economic models refer [124, 132-3].

## 3.2 Definition and Illustration of Fuzzy Relational Maps (FRMS)

In this section, we introduce the notion of Fuzzy Relational Maps (FRMs); they are constructed analogous to FCMs described and discussed in the earlier sections. In FCMs we promote the correlations between causal associations among concurrently active units. But in FRMs we divide the very causal associations into two disjoint units, for example, the relation between a teacher and a student or relation between an employee and an employer or a relation between doctor and patient and so on. Thus for us to define a FRM we need a domain space and a range space which are disjoint in the sense of concepts. We further assume no intermediate relation exists within the domain elements or node and the range spaces elements. The number of elements in the range space need not in general be equal to the number of elements in the domain space.

Thus throughout this section we assume the elements of the domain space are taken from the real vector space of dimension n and that of the range space are real vectors from the vector space of dimension m (m in general need not be equal to n). We denote by R the set of nodes $R_1,\ldots, R_m$ of the range space, where R = {$(x_1,\ldots, x_m)$ | $x_j$ = 0 or 1 } for j = 1, 2,…, m. If $x_i$ = 1 it means that the node $R_i$ is in the ON state and if $x_i$ = 0 it means that the node $R_i$ is in the OFF state. Similarly D denotes the nodes $D_1, D_2,\ldots, D_n$ of the domain space where D = {$(x_1,\ldots, x_n)$ | $x_j$ = 0 or 1} for i = 1, 2,…, n. If $x_i$ = 1 it means that the node $D_i$ is in the ON state and if $x_i$ = 0 it means that the node $D_i$ is in the OFF state.

Now we proceed on to define a FRM.

**DEFINITION 3.2.1:** *A FRM is a directed graph or a map from D to R with concepts like policies or events etc, as nodes and causalities as edges. It represents causal relations between spaces D and R.*

*Let $D_i$ and $R_j$ denote that the two nodes of an FRM. The directed edge from $D_i$ to $R_j$ denotes the causality of $D_i$ on $R_j$ called relations. Every edge in the FRM is weighted with a number in the set {0, ±1}. Let $e_{ij}$ be the weight of the edge $D_iR_j$, $e_{ij} \in \{0, ±1\}$. The weight of the edge $D_i R_j$ is positive if increase in $D_i$ implies increase in $R_j$ or decrease in $D_i$ implies decrease in $R_j$, i.e., causality of $D_i$ on $R_j$ is 1. If $e_{ij} = 0$, then $D_i$ does not have any effect on $R_j$. We do not discuss the cases when increase in $D_i$ implies decrease in $R_j$ or decrease in $D_i$ implies increase in $R_j$.*

**DEFINITION 3.2.2:** *When the nodes of the FRM are fuzzy sets then they are called fuzzy nodes. FRMs with edge weights {0, ±1} are called simple FRMs.*

**DEFINITION 3.2.3:** *Let $D_1, \ldots, D_n$ be the nodes of the domain space D of an FRM and $R_1, \ldots, R_m$ be the nodes of the range space R of an FRM. Let the matrix E be defined as E = $(e_{ij})$ where $e_{ij}$ is the weight of the directed edge $D_iR_j$ (or $R_jD_i$), E is called the relational matrix of the FRM.*

*Note*: It is pertinent to mention here that unlike the FCMs the FRMs can be a rectangular matrix with rows corresponding to the domain space and columns corresponding to the range space. This is one of the marked difference between FRMs and FCMs.

**DEFINITION 3.2.4:** *Let $D_1, \ldots, D_n$ and $R_1, \ldots, R_m$ denote the nodes of the FRM. Let $A = (a_1, \ldots, a_n)$, $a_i \in \{0, \pm 1\}$. A is called the instantaneous state vector of the domain space and it denotes the on-off position of the nodes at any instant. Similarly let $B = (b_1, \ldots, b_m)$, $b_i \in \{0, \pm 1\}$. B is called instantaneous state vector of the range space and it denotes the on-off position of the nodes at any instant; $a_i = 0$ if $a_i$ is off and $a_i = 1$ if $a_i$ is on for $i = 1, 2, \ldots, n$. Similarly, $b_i = 0$ if $b_i$ is off and $b_i = 1$ if $b_i$ is on, for $i = 1, 2, \ldots, m$.*

**DEFINITION 3.2.5:** *Let $D_1, \ldots, D_n$ and $R_1, \ldots, R_m$ be the nodes of an FRM. Let $D_i R_j$ (or $R_j D_i$) be the edges of an FRM, $j = 1, 2, \ldots, m$ and $i = 1, 2, \ldots, n$. Let the edges form a directed cycle. An FRM is said to be a cycle if it posses a directed cycle. An FRM is said to be acyclic if it does not posses any directed cycle.*

**DEFINITION 3.2.6:** *An FRM with cycles is said to be an FRM with feedback.*

**DEFINITION 3.2.7:** *When there is a feedback in the FRM, i.e. when the causal relations flow through a cycle in a revolutionary manner, the FRM is called a dynamical system.*

**DEFINITION 3.2.8:** *Let $D_i R_j$ (or $R_j D_i$), $1 \leq j \leq m$, $1 \leq i \leq n$. When $R_i$ (or $D_j$) is switched on and if causality flows through edges of the cycle and if it again causes $R_i$ (or $D_j$), we say that the dynamical system goes round and round. This is true for any node $R_j$ (or $D_i$) for $1 \leq i \leq n$, (or $1 \leq j \leq m$). The equilibrium state of this dynamical system is called the hidden pattern.*

**DEFINITION 3.2.9:** *If the equilibrium state of a dynamical system is a unique state vector, then it is called a fixed point.*

*Consider an FRM with $R_1, R_2, ..., R_m$ and $D_1, D_2, ..., D_n$ as nodes. For example, let us start the dynamical system by switching on $R_1$ (or $D_1$). Let us assume that the FRM settles down with $R_1$ and $R_m$ (or $D_1$ and $D_n$) on, i.e. the state vector remains as (1, 0, ..., 0, 1) in R) or (1, 0, 0, ..., 0, 1) in D), This state vector is called the fixed point.*

**DEFINITION 3.2.10:** *If the FRM settles down with a state vector repeating in the form*
$A_1 \to A_2 \to A_3 \to ... \to A_i \to A_1$ *(or* $B_1 \to B_2 \to ... \to B_i \to B_1$*)*
*then this equilibrium is called a limit cycle.*

Here we give the methods of determining the hidden pattern.

Let $R_1, R_2, ..., R_m$ and $D_1, D_2, ..., D_n$ be the nodes of a FRM with feedback. Let E be the relational matrix. Let us find a hidden pattern when $D_1$ is switched on i.e. when an input is given as vector $A_1 = (1, 0, ..., 0)$ in $D_1$, the data should pass through the relational matrix E. This is done by multiplying $A_1$ with the relational matrix E. Let $A_1 E = (r_1, r_2, ..., r_m)$, after thresholding and updating the resultant vector we get $A_1 E \in R$. Now let $B = A_1 E$, we pass on B into $E^T$ and obtain $BE^T$. We update and threshold the vector $BE^T$ so that $BE^T \in D$. This procedure is repeated till we get a limit cycle or a fixed point.

**DEFINITION 3.2.11:** *Finite number of FRMs can be combined together to produce the joint effect of all the FRMs. Let $E_1, ..., E_p$ be the relational matrices of the FRMs with nodes $R_1, R_2, ..., R_m$ and $D_1, D_2, ..., D_n$, then the combined FRM is represented by the relational matrix $E = E_1 + ... + E_p$.*

Now we give a simple illustration of a FRM, for more about FRMs please refer [136-7, 143].

*Example 3.2.1:* Let us consider the relationship between the teacher and the student. Suppose we take the domain space as the concepts belonging to the teacher say $D_1, ..., D_5$ and the range space denote the concepts belonging to the student say $R_1, R_2$ and $R_3$.

We describe the nodes $D_1,\ldots, D_5$ and $R_1$, $R_2$ and $R_3$ in the following:

Nodes of the Domain Space

$D_1$ – Teaching is good
$D_2$ – Teaching is poor
$D_3$ – Teaching is mediocre
$D_4$ – Teacher is kind
$D_5$ – Teacher is harsh [or rude].

(We can have more concepts like teacher is non-reactive, unconcerned etc.)

Nodes of the Range Space

$R_1$ – Good Student
$R_2$ – Bad Student
$R_3$ – Average Student.

The relational directed graph of the teacher-student model is given in Figure 3.2.1.

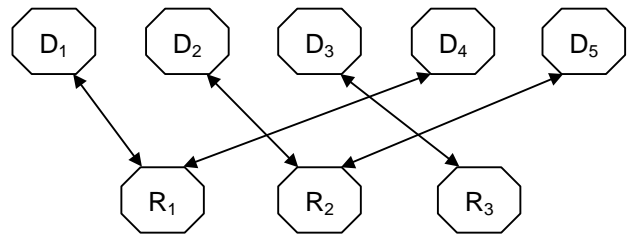

**FIGURE: 3.2.1**

The relational matrix E got from the above map is

$$E = \begin{bmatrix} 1 & 0 & 0 \\ 0 & 1 & 0 \\ 0 & 0 & 1 \\ 1 & 0 & 0 \\ 0 & 1 & 0 \end{bmatrix}$$

If A = (1 0 0 0 0) is passed on in the relational matrix E, the instantaneous vector, AE = (1 0 0) implies that the student is a good student . Now let AE = B, $BE^T$ = (1 0 0 1 0) which implies that the teaching is good and he / she is a kind teacher. Let $BE^T$ = $A_1$, $A_1E$ = (2 0 0) after thresholding we get $A_1E$ = (1 0 0) which implies that the student is good, so on and so forth.

### 3.3 Definition of the New Fuzzy Dynamical System

This new system is constructed when we have at hand the opinion of several experts. It functions more like an FRM but in the operations max min principle is used. We just describe how we construct it. We have n experts who give their opinion about the problem using p nodes along the column and m nodes along the rows. Now we define the new fuzzy system M = $(a_{ij})$ to be a m × p matrix with $(a_{ij}) \in [0, 1]$; $1 \leq i \leq m$ and $1 \leq j \leq p$, giving equal importance to the views of the n experts.

The only assumption is that all the n experts choose to work with the same p sets of nodes/ concepts along the columns and m sets of nodes/concepts along the rows. Suppose $P_1, \ldots, P_p$ denotes the nodes related with the columns and $C_1, \ldots, C_m$ denotes the nodes of the rows. Then $a_{ij}$ denotes how much or to which degree $C_i$ influences $P_j$ which is given a membership degree in the interval [0, 1] i.e., $a_{ij} \in [0, 1]$; $1 \leq i \leq m$ and $1 \leq j \leq p$ by any $t^{th}$ expert.

Now $M_t = (a^t_{ij})$ is a fuzzy m × p matrix which is defined as the new fuzzy vector matrix. We take the views of all the n experts and if $M_1, \ldots, M_n$ denotes the n number of fuzzy m × p matrices where $M_t = (a_{ij}^t)$; $1 \leq t \leq n$.

Let

$$M = \frac{M_1 + \ldots + M_n}{n}$$

$$= \frac{\left(a_{ij}^1\right) + \left(a_{ij}^2\right) + \ldots + \left(a_{ij}^n\right)}{n}$$

= $(a_{ij})$; $1 \leq i \leq m$ and $1 \leq j \leq p$.
i.e.,

$$a_{11} = \frac{a_{11}^1 + a_{11}^2 + ... + a_{11}^n}{n}$$

$$a_{12} = \frac{a_{12}^1 + a_{12}^2 + ... + a_{12}^n}{n}$$

and so on. Thus

$$a_{1j} = \frac{a_{1j}^1 + a_{1j}^2 + ... + a_{1j}^n}{n}.$$

The matrix $M = (a_{ij})$ is defined as the new fuzzy dynamical model of the n experts or the dynamical model of the multi expert n system. For it can simultaneously work with n experts view. Clearly $a_{ij} \in [0, 1]$, so M is called as the new fuzzy dynamical model. The working will be given in chapter IV.

### 3.4 Neutrosophic Cognitive Maps with Examples

The notion of Fuzzy Cognitive Maps (FCMs) which are fuzzy signed directed graphs with feedback are discussed and described in section 1 of this chapter. The directed edge $e_{ij}$ from causal concept $C_i$ to concept $C_j$ measures how much $C_i$ causes $C_j$. The time varying concept function $C_i(t)$ measures the non negative occurrence of some fuzzy event, perhaps the strength of a political sentiment, historical trend or opinion about some topics like child labor or school dropouts etc. FCMs model the world as a collection of classes and causal relations between them.

The edge $e_{ij}$ takes values in the fuzzy causal interval $[-1, 1]$ ($e_{ij} = 0$, indicates no causality, $e_{ij} > 0$ indicates causal increase; that $C_j$ increases as $C_i$ increases or $C_j$ decreases as $C_i$ decreases, $e_{ij} < 0$ indicates causal decrease or negative causality; $C_j$ decreases as $C_i$ increases or $C_j$, increases as $C_i$ decreases. Simple FCMs have edge value in $\{-1, 0, 1\}$. Thus if causality occurs it occurs to maximal positive or negative degree.

It is important to note that $e_{ij}$ measures only absence or presence of influence of the node $C_i$ on $C_j$ but till now any researcher has not contemplated the indeterminacy of any relation between two nodes $C_i$ and $C_j$. When we deal with unsupervised data, there are situations when no relation can be determined between some two nodes. So in this section we try to introduce the indeterminacy in FCMs, and we choose to call this generalized structure as Neutrosophic Cognitive Maps (NCMs). In our view this will certainly give a more appropriate result and also caution the user about the risk of indeterminacy [143].

Now we proceed on to define the concepts about NCMs.

**DEFINITION 3.4.1:** *A Neutrosophic Cognitive Map (NCM) is a neutrosophic directed graph with concepts like policies, events etc. as nodes and causalities or indeterminates as edges. It represents the causal relationship between concepts.*

Let $C_1, C_2, \ldots, C_n$ denote n nodes, further we assume each node is a neutrosophic vector from neutrosophic vector space V. So a node $C_i$ will be represented by $(x_1, \ldots, x_n)$ where $x_k$'s are zero or one or I (I is the indeterminate introduced in [...]) and $x_k = 1$ means that the node $C_k$ is in the ON state and $x_k = 0$ means the node is in the OFF state and $x_k = I$ means the nodes state is an indeterminate at that time or in that situation.

Let $C_i$ and $C_j$ denote the two nodes of the NCM. The directed edge from $C_i$ to $C_j$ denotes the causality of $C_i$ on $C_j$ called connections. Every edge in the NCM is weighted with a number in the set $\{-1, 0, 1, I\}$. Let $e_{ij}$ be the weight of the directed edge $C_iC_j$, $e_{ij} \in \{-1, 0, 1, I\}$. $e_{ij} = 0$ if $C_i$ does not have any effect on $C_j$, $e_{ij} = 1$ if increase (or decrease) in $C_i$ causes increase (or decreases) in $C_j$, $e_{ij} = -1$ if increase (or decrease) in $C_i$ causes decrease (or increase) in $C_j$. $e_{ij} = I$ if the relation or effect of $C_i$ on $C_j$ is an indeterminate.

**DEFINITION 3.4.2:** *NCMs with edge weight from {-1, 0, 1, I} are called simple NCMs.*

**DEFINITION 3.4.3:** *Let $C_1$, $C_2$, ..., $C_n$ be nodes of a NCM. Let the neutrosophic matrix $N(E)$ be defined as $N(E) = (e_{ij})$ where $e_{ij}$ is the weight of the directed edge $C_i C_j$, where $e_{ij} \in \{0, 1, -1, I\}$. $N(E)$ is called the neutrosophic adjacency matrix of the NCM.*

**DEFINITION 3.4.4:** *Let $C_1$, $C_2$, ..., $C_n$ be the nodes of the NCM. Let $A = (a_1, a_2, ..., a_n)$ where $a_i \in \{0, 1, I\}$. A is called the instantaneous state neutrosophic vector and it denotes the on – off – indeterminate state/ position of the node at an instant*

- $a_i$ = *0 if $a_i$ is off (no effect)*
- $a_i$ = *1 if $a_i$ is on (has effect)*
- $a_i$ = *I if $a_i$ is indeterminate(effect cannot be determined)*

*for $i = 1, 2, ..., n$.*

**DEFINITION 3.4.5:** *Let $C_1$, $C_2$, ..., $C_n$ be the nodes of the FCM. Let $\overrightarrow{C_1C_2}$, $\overrightarrow{C_2C_3}$, $\overrightarrow{C_3C_4}$, ..., $\overrightarrow{C_iC_j}$ be the edges of the NCM. Then the edges form a directed cycle. An NCM is said to be cyclic if it possesses a directed cycle. An NCM is said to be acyclic if it does not possess any directed cycle.*

**DEFINITION 3.4.6:** *An NCM with cycles is said to have a feedback. When there is a feedback in the NCM i.e. when the causal relations flow through a cycle in a revolutionary manner the NCM is called a neutrosophic dynamical system.*

**DEFINITION 3.4.7:** *Let $\overrightarrow{C_1C_2}, \overrightarrow{C_2C_3}, \cdots, \overrightarrow{C_{n-1}C_n}$ be a cycle, when $C_i$ is switched on and if the causality flows through the edges of a cycle and if it again causes $C_i$, we say that the dynamical system goes round and round. This is true for any node $C_i$, for $i = 1, 2, ..., n$. The equilibrium state for this dynamical system is called the hidden pattern.*

**DEFINITION 3.4.8:** *If the equilibrium state of a dynamical system is a unique state vector, then it is called a fixed point. Consider the NCM with $C_1$, $C_2$, ..., $C_n$ as nodes. For example let us start the dynamical system by switching on $C_1$. Let us assume* *that the NCM settles down with $C_1$ and $C_n$ on, i.e. the state vector remain as (1, 0,…, 1), this neutrosophic state vector (1,0,…, 0, 1) is called the fixed point.*

**DEFINITION 3.4.9:** *If the NCM settles with a neutrosophic state vector repeating in the form*
$$A_1 \to A_2 \to ... \to A_i \to A_1,$$
*then this equilibrium is called a limit cycle of the NCM.*

The methods of determining the hidden pattern is described in the following:

Let $C_1, C_2, …, C_n$ be the nodes of an NCM, with feedback. Let E be the associated adjacency matrix. Let us find the hidden pattern when $C_1$ is switched on, when an input is given as the vector $A_1 = (1, 0, 0,…, 0)$, the data should pass through the neutrosophic matrix N(E), this is done by multiplying $A_1$ by the matrix N(E). Let $A_1N(E) = (a_1, a_2,…, a_n)$ with the threshold operation that is by replacing $a_i$ by 1 if $a_i > k$ and $a_i$ by 0 if $a_i < k$ (k – a suitable positive integer) and $a_i$ by I if $a_i$ is not a integer. We update the resulting concept, the concept $C_1$ is included in the updated vector by making the first coordinate as 1 in the resulting vector. Suppose $A_1N(E) \to A_2$ then consider $A_2N(E)$ and repeat the same procedure. This procedure is repeated till we get a limit cycle or a fixed point.

**DEFINITION 3.4.10:** *Finite number of NCMs can be combined together to produce the joint effect of all NCMs. If $N(E_1)$, $N(E_2),…, N(E_p)$ be the neutrosophic adjacency matrices of a NCM with nodes $C_1, C_2,…, C_n$ then the combined NCM is got by adding all the neutrosophic adjacency matrices $N(E_1),…, N(E_p)$. We denote the combined NCMs adjacency neutrosophic matrix by $N(E) = N(E_1) + N(E_2)+…+ N(E_p)$.*

**NOTATION:** Let $(a_1, a_2, … , a_n)$ and $(a'_1, a'_2, … , a'_n)$ be two neutrosophic vectors. We say $(a_1, a_2, … , a_n)$ is equivalent to $(a'_1, a'_2, … , a'_n)$ denoted by $((a_1, a_2, … , a_n) \sim (a'_1, a'_2, …, a'_n)$ if $(a'_1, a'_2, … , a'_n)$ is got after thresholding and updating the vector ($a_1, a_2, \ldots, a_n$) after passing through the neutrosophic adjacency matrix N(E).

The following are very important:

***Note 1:*** The nodes $C_1, C_2, \ldots, C_n$ are not indeterminate nodes because they indicate the concepts which are well known. But the edges connecting $C_i$ and $C_j$ may be indeterminate i.e. an expert may not be in a position to say that $C_i$ has some causality on $C_j$ either will he be in a position to state that $C_i$ has no relation with $C_j$ in such cases the relation between $C_i$ and $C_j$ which is indeterminate is denoted by I.

***Note 2:*** The nodes when sent will have only ones and zeros i.e. ON and OFF states, but after the state vector passes through the neutrosophic adjacency matrix the resultant vector will have entries from {0, 1, I} i.e. they can be neutrosophic vectors, i.e., it may happen the node under those circumstances may be an indeterminate.

The presence of I in any of the coordinate implies the expert cannot say the presence of that node i.e. ON state of it, after passing through N(E) nor can we say the absence of the node i.e. OFF state of it, the effect on the node after passing through the dynamical system is indeterminate so only it is represented by I. Thus only in case of NCMs we can say the effect of any node on other nodes can also be indeterminates. Such possibilities and analysis is totally absent in the case of FCMs.

***Note 3:*** In the neutrosophic matrix N(E), the presence of I in the $a_{ij}^{th}$ place shows, that the causality between the two nodes i.e. the effect of $C_i$ on $C_j$ is indeterminate. Such chances of being indeterminate is very possible in case of unsupervised data and that too in the study of FCMs which are derived from the directed graphs.

Thus only NCMs helps in such analysis.
Now we shall represent a few examples to show how in this set up NCMs is preferred to FCMs. At the outset before we proceed to give examples we make it clear that all unsupervised data need not have NCMs to be applied to it. Only data which have the relation between two nodes to be an indeterminate need to be modeled with NCMs if the data has no indeterminacy factor between any pair of nodes, one need not go for NCMs; FCMs will do the best job.

***Example 3.4.1:*** The child labor problem prevalent in India is modeled in this example using NCMs.

Let us consider the child labor problem with the following conceptual nodes;

$C_1$ - Child Labor
$C_2$ - Political Leaders
$C_3$ - Good Teachers
$C_4$ - Poverty
$C_5$ - Industrialists
$C_6$ - Public practicing/encouraging Child Labor
$C_7$ - Good Non-Governmental Organizations (NGOs).

$C_1$ - Child labor, it includes all types of labor of children below 14 years which include domestic workers, rag pickers, working in restaurants / hotels, bars etc. (It can be part time or fulltime).
$C_2$ - We include political leaders with the following motivation: Children are not vote banks, so political leaders are not directly concerned with child labor but they indirectly help in the flourishing of it as industrialists who utilize child laborers or cheap labor; are the decision makers for the winning or losing of the political leaders. Also industrialists financially control political interests. So we are forced to include political leaders as a node in this problem.
$C_3$ - Teachers are taken as a node because mainly school dropouts or children who have never attended the school are child laborers. So if the motivation by the teacher is very good, there would be less school dropouts and therefore there would be a decrease in child laborers.

$C_4$ - Poverty which is the most important reason for child labor.

$C_5$ - Industrialists – when we say industrialists we include one and all starting from a match factory or beedi factory, bars, hotels rice mill, garment industries etc.

$C_6$ - Public who promote child labor as domestic servants, sweepers etc.

$C_7$ - We qualify the NGOs as good for some NGOs may not take up the issue fearing the rich and the powerful. Here "good NGOs" means NGOs who try to stop or prevent child labor.

Now we give the directed graph as well as the neutrosophic graph of two experts in the following Figures 3.4.1 and 3.4.2:

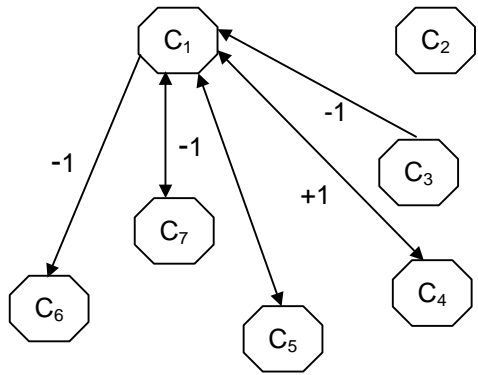

FIGURE: 3.4.1

Figure 3.4.1 gives the directed graph with $C_1$, $C_2$, …, $C_7$ as nodes and Figure 3.4.2 gives the neutrosophic directed graph with the same nodes.

The connection matrix E related to the directed neutrosophic graph given in Figure 3.4.1. which is the associated graph of NCM is given in the following:

$$E = \begin{bmatrix} 0 & 0 & 0 & 1 & 1 & 1 & -1 \\ 0 & 0 & 0 & 0 & 0 & 0 & 0 \\ -1 & 0 & 0 & 0 & 0 & 0 & 0 \\ 1 & 0 & 0 & 0 & 0 & 0 & 0 \\ 1 & 0 & 0 & 0 & 0 & 0 & 0 \\ 0 & 0 & 0 & 0 & 0 & 0 & 0 \\ -1 & 0 & 0 & 0 & 0 & 0 & 0 \end{bmatrix}.$$

According to this expert no connection however exists between political leaders and industrialists.

Now we reformulate a different format of the questionnaire where we permit the experts to give answers like the relation between certain nodes is indeterminable or not known. Now based on the expert's opinion also about the notion of indeterminacy we give the following neutrosophic directed graph of the expert:

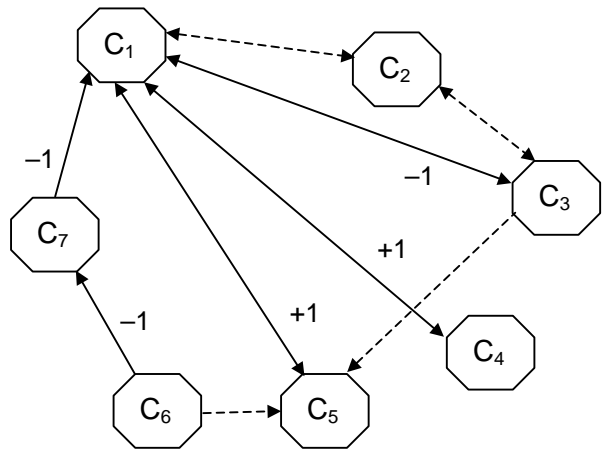

**FIGURE: 3.4.2**

The corresponding neutrosophic adjacency matrix N(E) related to the neutrosophic directed graph (Figure 3.4.2.) is given by the following:

$$N(E) = \begin{bmatrix} 0 & I & -1 & 1 & 1 & 0 & 0 \\ I & 0 & I & 0 & 0 & 0 & 0 \\ -1 & I & 0 & 0 & I & 0 & 0 \\ 1 & 0 & 0 & 0 & 0 & 0 & 0 \\ 1 & 0 & 0 & 0 & 0 & 0 & 0 \\ 0 & 0 & 0 & 0 & I & 0 & -1 \\ -1 & 0 & 0 & 0 & 0 & 0 & 0 \end{bmatrix}.$$

Suppose we take the state vector $A_1 = (1\ 0\ 0\ 0\ 0\ 0\ 0)$. We will see the effect of $A_1$ on E and on N(E).

$$\begin{aligned}
A_1 E &= (0\ 0\ 0\ 1\ 1\ 1\ -1) \\
&\to (1\ 0\ 0\ 1\ 1\ 1\ 0) \\
&= A_2. \\
A_2 E &= (2\ 0\ 0\ 1\ 1\ 1\ 0) \\
&\to (1\ 0\ 0\ 1\ 1\ 1\ 0) \\
&= A_3 = A_2.
\end{aligned}$$

Thus child labor flourishes with parents' poverty and industrialists' action. Public practicing child labor also flourish; but good NGOs are absent in such a scenario. The state vector gives the fixed point.

Now we find the effect of the state vector $A_1 = (1\ 0\ 0\ 0\ 0\ 0\ 0)$ on N(E).

$$\begin{aligned}
A_1 N(E) &= (0\ I\ -1\ 1\ 1\ 0\ 0) \\
&\to (1\ I\ 0\ 1\ 1\ 0\ 0) \\
&= A_2. \\
A_2 N(E) &= (I+2,\ I,\ -1+I,\ 1\ 1\ 0\ 0) \\
&\to (1\ I\ 0\ 1\ 1\ 0\ 0) \\
&= A_2.
\end{aligned}$$

Thus $A_2 = (1\ I\ 0\ 1\ 1\ 0\ 0)$, according to this expert the increase or the ON state of child labor certainly increases with the poverty of parents and other factors are indeterminate to him. This mainly gives the indeterminates relating to political leaders and teachers in the neutrosophic cognitive model and the parents poverty and industrialist activities become ON state.

However, the results by FCM give as if there is no effect by teachers and politicians for the increase in child labor. Actually the increase in school dropout increases the child labor hence certainly the role of teachers play a part. At least if it is termed as an indeterminate one would think or reflect about their (teachers) effect on child labor.

## 3.5 Description of Neutrosophic Relational Maps

Neutrosophic Cognitive Maps (NCMs) promote the causal relationships between concurrently active units or decides the absence of any relation between two units or the indeterminacy of any relation between any two units. But in Neutrosophic Relational Maps (NRMs) we divide the very causal nodes into two disjoint units. Thus for the modeling of a NRM we need a domain space and a range space which are disjoint in the sense of concepts. We further assume no intermediate relations exist within the domain and the range spaces. The number of elements or nodes in the range space need not be equal to the number of elements or nodes in the domain space.

Throughout this section we assume the elements of a domain space are taken from the neutrosophic vector space of dimension n and that of the range space are neutrosophic vector space of dimension m. (m in general need not be equal to n). We denote by R the set of nodes $R_1,\ldots, R_m$ of the range space, where $R = \{(x_1,\ldots, x_m) \mid x_j = 0$ or $1$ for $j = 1, 2, \ldots, m\}$.

If $x_i = 1$ it means that node $R_i$ is in the ON state and if $x_i = 0$ it means that the node $R_i$ is in the OFF state and if $x_i = I$ in the resultant vector it means the effect of the node $x_i$ is indeterminate or whether it will be OFF or ON cannot be predicted by the neutrosophic dynamical system.

It is very important to note that when we send the state vectors they are always taken as the real state vectors for we know the node or the concept is in the ON state or in the off state but when the state vector passes through the Neutrosophic dynamical system some other node may become indeterminate i.e. due to the presence of a node we may not be able to predict the presence or the absence of the other node i.e., it is indeterminate, denoted by the symbol I, thus the resultant vector can be a neutrosophic vector.

**DEFINITION 3.5.1:** *A Neutrosophic Relational Map (NRM) is a neutrosophic directed graph or a map from D to R with concepts like policies or events etc. as nodes and causalities as edges. (Here by causalities we mean or include the indeterminate causalities also). It represents Neutrosophic Relations and Causal Relations between spaces D and R.*

*Let $D_i$ and $R_j$ denote the nodes of an NRM. The directed edge from $D_i$ to $R_j$ denotes the causality of $D_i$ on $R_j$ called relations. Every edge in the NRM is weighted with a number in the set $\{0, +1, -1, I\}$. Let $e_{ij}$ be the weight of the edge $D_i R_j$, $e_{ij} \in \{0, 1, -1, I\}$. The weight of the edge $D_i R_j$ is positive if increase in $D_i$ implies increase in $R_j$ or decrease in $D_i$ implies decrease in $R_j$ i.e. causality of $D_i$ on $R_j$ is 1. If $e_{ij} = -1$ then increase (or decrease) in $D_i$ implies decrease (or increase) in $R_j$. If $e_{ij} = 0$ then $D_i$ does not have any effect on $R_j$. If $e_{ij} = I$ it implies we are not in a position to determine the effect of $D_i$ on $R_j$ i.e. the effect of $D_i$ on $R_j$ is an indeterminate so we denote it by I.*

**DEFINITION 3.5.2:** *When the nodes of the NRM take edge values from $\{0, 1, -1, I\}$ we say the NRMs are simple NRMs.*

**DEFINITION 3.5.3:** *Let $D_1, …, D_n$ be the nodes of the domain space D of an NRM and let $R_1, R_2,…, R_m$ be the nodes of the range space R of the same NRM. Let the matrix N(E) be defined as $N(E) = (e_{ij})$ where $e_{ij}$ is the weight of the directed edge $D_i R_j$ (or $R_j D_i$) and $e_{ij} \in \{0, 1, -1, I\}$. N(E) is called the Neutrosophic Relational Matrix of the NRM.*

The following remark is important and interesting to find its mention in this book [143].

**Remark**: Unlike NCMs, NRMs can also be rectangular matrices with rows corresponding to the domain space and columns corresponding to the range space. This is one of the marked difference between NRMs and NCMs. Further the number of entries for a particular model which can be treated as disjoint sets when dealt as a NRM has very much less entries than when the same model is treated as a NCM.

Thus in many cases when the unsupervised data under study or consideration can be spilt as disjoint sets of nodes or concepts; certainly NRMs are a better tool than the NCMs when time and money is a criteria.

**DEFINITION 3.5.4:** *Let $D_1, \ldots, D_n$ and $R_1, \ldots, R_m$ denote the nodes of a NRM. Let $A = (a_1, \ldots, a_n)$, $a_i \in \{0, 1, -I\}$ is called the Neutrosophic instantaneous state vector of the domain space and it denotes the on-off position or an indeterminate state of the nodes at any instant. Similarly let $B = (b_1, \ldots, b_n)$ $b_i \in \{0, 1, -I\}$, B is called instantaneous state vector of the range space and it denotes the on-off position or an indeterminate state of the nodes at any instant, $a_i = 0$ if $a_i$ is off and $a_i = 1$ if $a_i$ is on, $a_i = I$ if the state is an indeterminate one at that time for $i = 1, 2, \ldots, n$. Similarly, $b_i = 0$ if $b_i$ is off and $b_i = 1$ if $b_i$ is on, $b_i = I$ i.e., the state of $b_i$ is an indeterminate at that time for $i = 1, 2, \ldots, m$.*

**DEFINITION 3.5.5:** *Let $D_1, \ldots, D_n$ and $R_1, R_2, \ldots, R_m$ be the nodes of a NRM. Let $D_i R_j$ (or $R_j D_i$) be the edges of an NRM, $j = 1, 2, \ldots, m$ and $i = 1, 2, \ldots, n$. The edges form a directed cycle. An NRM is said to be a cycle if it possess a directed cycle. An NRM is said to be acyclic if it does not possess any directed cycle.*

**DEFINITION 3.5.6:** *A NRM with cycles is said to be a NRM with feedback.*

**DEFINITION 3.5.7:** *When there is a feedback in the NRM i.e. when the causal relations flow through a cycle in a revolutionary manner, the NRM is called a neutrosophic dynamical system.*

**DEFINITION 3.5.8:** *Let $D_i R_j$ (or $R_j D_i$), $1 \leq j \leq m$, $1 \leq i \leq n$, when $R_j$ (or $D_i$) is switched on and if causality flows through edges of a cycle and if it again causes $R_j$ (or $D_i$) we say that the neutrosophic dynamical system goes round and round. This is* *true for any node $R_j$ ( or $D_i$ ) for $1 \leq j \leq m$ (or $1 \leq i \leq n$). The equilibrium state of this neutrosophic dynamical system is called the Neutrosophic hidden pattern.*

**DEFINITION 3.5.9:** *If the equilibrium state of a neutrosophic dynamical system is a unique neutrosophic state vector, then it is called the fixed point. Consider an NRM with $R_1, R_2, ..., R_m$ and $D_1, D_2, ..., D_n$ as nodes. For example let us start the dynamical system by switching on $R_1$ (or $D_1$). Let us assume that the NRM settles down with $R_1$ and $R_m$ (or $D_1$ and $D_n$) on, or indeterminate on, i.e. the neutrosophic state vector remains as (1, 0, 0,..., 1) or (1, 0, 0,...I) (or (1, 0, 0,...1) or (1, 0, 0,...I) in D), this state vector is called the fixed point.*

**DEFINITION 3.5.10:** *If the NRM settles down with a state vector repeating in the form $A_1 \to A_2 \to A_3 \to ... \to A_i \to A_1$ (or $B_1 \to B_2 \to ... \to B_i \to B_1$) then this equilibrium is called a limit cycle.*

We describe the methods of determining the hidden pattern in a NRM.

Let $R_1, R_2,..., R_m$ and $D_1, D_2,..., D_n$ be the nodes of a NRM with feedback. Let $N(E)$ be the neutrosophic Relational Matrix. Let us find the hidden pattern when $D_1$ is switched on i.e. when an input is given as a vector; $A_1 = (1, 0, ..., 0)$ in D; the data should pass through the relational matrix $N(E)$. This is done by multiplying $A_1$ with the neutrosophic relational matrix $N(E)$. Let $A_1 N(E) = (r_1, r_2,..., r_m)$ after thresholding and updating the resultant vector we get $A_1 E \in R$, Now let $B = A_1 E$ we pass on B into the system $(N(E))^T$ and obtain $B(N(E))^T$. We update and threshold the vector $B(N(E))^T$ so that $B(N(E))^T \in D$.

This procedure is repeated till we get a limit cycle or a fixed point.

**DEFINITION 3.5.11:** *Finite number of NRMs can be combined together to produce the joint effect of all NRMs. Let $N(E_1)$, $N(E_2),..., N(E_r)$ be the neutrosophic relational matrices of the NRMs with nodes $R_1,..., R_m$ and $D_1,...,D_n$, then the combined*

*NRM* is represented by the neutrosophic relational matrix $N(E)$
$= N(E_1) + N(E_2) + \ldots + N(E_r)$.
Now we give a simple illustration of a NRM.

***Example 3.5.1:*** Now consider the example given in the section two of this chapter. We take $D_1, D_2, \ldots, D_5$ and the $R_1, R_2$ and $R_3$ as in Example 3.2.1:

$D_1$ – Teacher is good
$D_2$ – Teaching is poor
$D_3$ – Teaching is mediocre
$D_4$ – Teacher is kind
$D_5$ – Teacher is harsh (or Rude).

$D_1, \ldots, D_5$ are taken as the 5 nodes of the domain space, we consider the following 3 nodes to be the nodes of the range space.

$R_1$ – Good student
$R_2$ – Bad student
$R_3$ – Average student.

The Neutrosophic relational graph of the teacher student model is as follows:

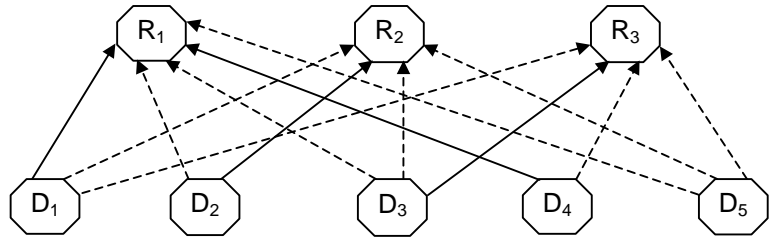

**FIGURE: 3.5.2**

$$N(E) = \begin{bmatrix} 1 & I & I \\ I & 1 & 0 \\ I & I & 1 \\ 1 & 0 & I \\ I & I & I \end{bmatrix}.$$

If $A_1 = (1, 0, 0, 0, 0)$ is taken as the instantaneous state vector and is passed on in the relational matrix $N(E)$, $A_1 N(E) = (1, I, I) = A_2$.
Now

$$A_2(N(E))^T = (1 + I, 1 + I, I\ 1 + I\ I)$$
$$\rightarrow (1\ 1\ I\ 1\ I)$$
$$= B_1$$
$$B_1 N(E) = (2 + I, I + 1, I)$$
$$\rightarrow (1\ I\ I)$$
$$= A_3$$

$$A_3 N(E) = (1 + I, I, I, 1 + I, I)$$
$$\rightarrow (1\ I\ I\ 1\ I)$$
$$= B_2 = B_1.$$
$$B_1 N(E) = (1\ I\ I).$$

Thus we see from the NRM given that if the teacher is good it implies it produces good students but nothing can be said about bad and average students. The bad and average students remain as indeterminates.

On the other hand in the domain space when the teacher is good the teaching quality of her remains indeterminate therefore both the nodes teaching is poor and teaching is mediocre remains as indeterminates but the node teacher is kind becomes in the ON state and the teacher is harsh is an indeterminate, (for harshness may be present depending on the circumstances).

### 3.6 Description of the new Fuzzy Neutrosophic model

In this section we for the first time introduce the new model which can evaluate the opinion of multiexperts say (n experts, n a positive integer) at a time (i.e., simultaneously). We call this the new fuzzy neutrosophic dynamical n expert system. This is constructed in the following way.

We assume I is the indeterminate and $I^2 = I$. We further define the fuzzy neutrosophic interval as $N_I = [0, 1] \cup [0, I]$ i.e.,

elements x of $N_I$ will be of the form $x = a + bI$ ($a, b \in [0, 1]$); x will be known as the fuzzy neutrosophic number.

A matrix $M = (a_{ij})$ where $a_{ij} \in N_I$ i.e., $a_{ij}$ are fuzzy neutrosophic numbers, will be called as the fuzzy neutrosophic matrix. We will be using only fuzzy neutrosophic matrix in the new fuzzy neutrosophic multiexpert system. Let us consider a problem P on which say some n experts give their views. In the first place the data related with the problem is an unsupervised one. Let the problem have m nodes taken as the rows and p nodes takes as the columns of the fuzzy neutrosophic matrix. Suppose we make the two assumptions;

1. All the n experts work only with same set of m nodes as rows and p nodes as columns.

2. All the experts have their membership function only from the fuzzy neutrosophic interval $N_I$. Let $M_t = (a_{ij}^t)$ be the fuzzy neutrosophic matrix given by the $t^{th}$ expert $t = 1, 2, \ldots, n$ i.e., $a_{ij}^t$ represent to which fuzzy neutrosophic degree the node $m_i$ is related with the node $p_j$ for $1 \leq i \leq m$ and $1 \leq j \leq p$. Thus $M_t = (a_{ij}^t)$ is the fuzzy neutrosophic matrix given by the $t^{th}$ expert. Let $M_1 = (a_{ij}^1)$, $M_2 = (a_{ij}^2)$, $\ldots$, $M_n = (a_{ij}^n)$ be the set of n fuzzy neutrosophic matrices given by the n experts. The new fuzzy neutrosophic multi n expert system $M = (a_{ij})$; $a_{ij} \in N_I$ is defined as follows:

Define
$$M = \frac{M_1 + M_2 + \ldots + M_n}{n}$$

$$= \frac{(a_{ij}^1) + (a_{ij}^2) + \ldots + (a_{ij}^n)}{n}$$

$= (a_{ij})$; $1 \leq i \leq m$ and $1 \leq j \leq p$.

i.e.,
$$a_{11} = \frac{a_{11}^1 + a_{11}^2 + \ldots + a_{11}^n}{n}$$

$$a_{12} = \frac{a_{12}^1 + a_{12}^2 + \ldots + a_{12}^n}{n}$$

and so on.

Now this system functions similar to the fuzzy dynamical system described in 3.3. of this book; the only difference is that their entries are from the fuzzy interval [0, 1] and in case of fuzzy neutrosophic dynamical system the entries are from $N_I$ = [0, 1] $\cup$ [0, I].

**Chapter Four**

# MATHEMATICAL ANALYSIS OF THE VIEWS ABOUT VEDIC MATHEMATICS USING FUZZY MODELS

In this chapter we use fuzzy and neutrosophic analysis to study the ulterior motives of imposing Vedic Mathematics in schools. The subsequent study led up to the question, "Is Vedic Mathematics, Vedic (derived from the Vedas) or Mathematics?" While trying to analyze about Vedic Mathematics from five different categories of people: students, teachers, parents, educationalists and public we got the clear picture that Vedic Mathematics does not contain any sound exposition to Vedas, nor is it mathematics. All these groups unanimously agreed upon the fact that the Vedic Mathematics book authored by the Swamiji contained only simple arithmetic of primary school standard. All the five categories of people could not comment on its Vedic content for it had no proper citation from the Vedas. And in some of the groups, people said that the book did not contain any Vedas of standard. Some people acknowledged that the content of Vedas itself was an indeterminate because in their opinion the Vedas itself was a trick to ruin the non-Brahmins and elevate the Brahmins. They pointed out that the Vedic Mathematics book also does it very cleverly. They said that the mathematical contents in Vedic Mathematics was zero and the Vedic contents was an indeterminate. This argument was substantiated because cunning and ulterior motives are richly present in the book where Kamsa is described and decried as a Sudra king of arrogance! This is an instance to show the 'charm' of Vedic Mathematics.

This chapter has five sections. In section one we give the views of the students about the use of Vedic Mathematics in their mathematics curriculum. In section two we analyze the feelings of teachers using FRM and NRM described in chapter three.

In section three we give the opinion of parents about Vedic Mathematics. The group (parents) was heterogeneous because some were educated, many were uneducated, some knew about Vedic Mathematics and some had no knowledge about it. So, we could not use any mathematical tool, and as in the case of students, the data collected from them could not be used for mathematical analysis because majority of them used a 'single term' in their questionnaire; hence any attempt at grading became impossible. The fourth section of our chapter uses the new fuzzy dynamic multi-expert model described in chapter 3, section 3 to analyze the opinion of the educated people about Vedic Mathematics. Also the fuzzy neutrosophic multi n-expert model described in section 3.6 is used to analyze the problem. The final section uses both FCMs and NCMs to study and analyze the public opinion on Vedic Mathematics.

In this chapter, the analysis of 'How 'Vedic' is Vedic Mathematics' was carried out using fuzzy and neutrosophic theory for the 5 peer groups. The first category is students who had undergone at least some classes in Vedic Mathematics. The second category consisted of teachers followed by the third group which constituted of parents.

The fourth group was made up of educationalists who were aware of Vedic Mathematics. The final group, that is, the public included politicians, heads of other religions, rationalists and so on. We have been first forced to use students as they are the first affected, followed by parents and teachers who are directly related with the students. One also needs the opinion of educationalists. Further, as this growth and imposition of Vedic Mathematics is strongly associated with a revivalist, political party we have included the views of both the public and the politicians.

## 4.1 Views of students about the use of Vedic Mathematics in their curriculum

We made a linguistic questionnaire for the students and asked them to fill and return it to us. Our only criteria was that these students must have attended Vedic Mathematics classes. We prepared 100 photocopies of the questionnaire. However, we could get back only 92 of the filled-in forms. The main questions listed in the questionnaire are given below; we have also given the gist of the answers provided by them.

1. What is the standard of the mathematics taught to you in Vedic Mathematics classes?

   The mathematics taught to us in Vedic Mathematics classes was very elementary (90 out of 92 responses). They did only simple arithmetical calculations, which we have done in our primary classes (16 of them said first standard mathematics). Two students said that it was the level of sixth standard.

2. Did you like the Vedic Mathematics classes?

   The typical answer of the students was: "Utterly boring! Just like UKG/ LKG students who repeat rhymes we were asked to say the sutras loudly everyday before the classes started, we could never get the meaning of the sutras!"

3. Did you attend Vedic Mathematics classes out of interest or out of compulsion?

   Everybody admitted that they studied it out of "compulsion"; they said, "if we don't attend the classes, our parents will be called and if we cut classes we have to pay a hefty fine and write the sentences like "I won't repeat this" or "I would not be absent for Vedic Mathematics classes" some 100 times and get this countersigned by our parents." They shared the opinion that nothing was 'interesting' about Vedic Mathematics classes and only simple tricks of elementary arithmetic was taught.

4. Did you pay any fees or was the Vedic Mathematics classes free?

   In a year they were asked to pay Rs.300/- (varies from school to school) for Vedic Mathematics classes. In some schools, the classes were for one month duration, in some schools 3 months duration. Only in a few schools was the subject taught throughout the year (weekly one class). The students further added, "We have to buy the Vedic Mathematics textbooks compulsorily. A salesman from the bookstore Motilal Banarsidass from Mylapore, Chennai sold these books."

5. Who took Vedic Mathematics classes?

   In some schools, the mathematics teachers took the classes. In some schools new teachers from some other schools or devotees from religious *mathas* took classes.

6. Did you find any difference between Vedic Mathematics classes and your other classes?

   At the start of the Vedic Mathematics class they were made to recite long Sanskrit slokas. They also had to end the class with recitation of Sanskrit slokas! A few students termed this a "Maha-bore".

7. Did the Vedic Mathematics teacher show any partiality or discrimination in the class?

   "Some teachers unnecessarily scolded some of our friends and punished them. They unduly scolded the Christian boys and non-Brahmin friends who had dark complexion. Discrimination was explicit." Some teachers had asked openly in the class, "How many of you have had the *upanayana* (sacred thread) ceremony?"

8. How useful is Vedic Mathematics in doing your usual mathematical courses?

   Absolutely no use (89 out of 92 mentioned so).

9. Does Vedic Mathematics help in the competitive exams?

   No connection or relevance (90 out 92).

10. Do you feel Vedic Mathematics can be included in the curriculum?

    It is already taught in primary classes under General Mathematics so there is no need to waste our time rereading it under the title of Vedic Mathematics was the answer from the majority of the students (89 out of 92).

11. Do you find any true relation between the sutra they recite and the problem solved under that sutra?

    No. No sutra looks like a formula or a theorem. So we don't see any mathematics or scientific term or formula in them.

12. Can Vedic Mathematics help you in any other subject?

    Never. Because it is very elementary and useless.

13. Is Vedic Mathematics high level (or advanced) mathematics?

    No it is only very simple arithmetic.

14. Were you taught anything like higher-level Vedic Mathematics ?

    No. Every batch was not taught any higher level Vedic Mathematics, only elementary calculations were taught to all of us. Only in the introductory classes they had given a

long lecture about how Vedic Mathematics is used in all 'fields' of mathematics but students were utterly disappointed to learn this simple arithmetic.

Here we wish to state that only after we promised to keep their identity anonymous, the students filled the questionnaire. Only 5 students out of the 92 respondents filled in their name and classes. They were probably afraid of their teachers and the school administration. Though they spoke several things orally (with a lot of enthusiasm) they did not wish to give in writing. The questionnaire had linguistic terms like: very useful, just useful, somewhat useful, cannot say, useless, absolutely useless and so on.

In majority of the cases they ticked useless or absolutely useless. Other comments were filled by phrases like 'Boring', 'Maha Bore', 'Killing our time', 'We are back in primary class' and so on. The composition of the students was heterogeneous: that is, it was drawn from both Brahmins and non-Brahmins. Some Christian students had remarked that it was only like Vedic Hindu classes and their parents had expressed objections to it.

The most important thing to be observed is that these classes were conducted unofficially by the schools run by Hindu trusts with BJP/RSS background. None of the schools run by the Government, or Christian or Muslim trusts ever conducted such classes.

***Remark:*** We supplied the students with a linguistic questionnaire with 57 questions, and students were asked to select a linguistic phrase as answer, or in some cases, express their opinion in short sentences. But to our disappointment they had ticked in the questionnaire choices like useless, absolutely useless, nothing, no use in a very careless way which only reflected their scant regard or interest in those classes. So using these response we found it impossible to apply any form of fuzzy tool to analyze the data mathematically, so we had no other option except to give their overall feelings in the last chapter on conclusions.
100

The final question "any other information or any other suggestion" elicited these responses:

They wanted this class to be converted into a computer class, a karate class or a class which prepared them for entrance exams, so that they could be benefited by it. What is the use when we have calculators for all calculations? Some said that their cell phone would serve the purpose of Vedic Mathematics. They feel that in times of modernity, these elementary arithmetic techniques are an utter waste. We have listed the observations not only from the contents of the filled-in questionnaire but also from our discussions. We have also included discussions with students who have not undergone Vedic Mathematics classes. The observations from them will also be given in the last chapter. The views of rural students who have not been taught Vedic Mathematics, but to whom we explained the techniques used are also given.

## 4.2 Teachers views on Vedic Mathematics and its overall influence on the Students Community

We held discussion with nearly 200 teachers from urban schools, rural schools and posh city schools. Also teachers from corporation schools and government schools were interviewed. We could not ask them to fill a questionnaire or ask them to give any write up. Some of them had not even seen the Vedic Mathematics book.

Only very few of them had seen it and some had taught it to students. So the crowd which we had to get views from was an heterogeneous one and they belonged to different types of schools some of which promoted Vedic Mathematics and some of which strongly opposed Vedic Mathematics. Thus we got their views through discussions and noted the vital points which will be used to draw conclusions about the course on Vedic Mathematics to the students.

The majority of them spoke about these 8 concepts in one way or other in their discussion.

    $D_1$ - The mathematical content of Vedic Mathematics.
    $D_2$ - Vedic value of Vedic Mathematics.
    $D_3$ - Religious values of Vedic Mathematics.
    $D_4$ - Use of Vedic Mathematics in higher learning.
    $D_5$ - Why is it called Vedic Mathematics?
    $D_6$ - Vedic Mathematics is a waste for school children.
    $D_7$ - Vedic Mathematics is used to globalize Hindutva.
    $D_8$ - Vedic Mathematics will induce caste and discrimination among children and teachers.

These eight attributes are given by majority of the teachers which is taken as the nodes or concepts related with the domain space.

The following were given by majority of the teachers about the standard and use of Vedic Mathematics.

    $R_1$ - Vedic Mathematics is very elementary
    $R_2$ - Vedic Mathematics is primary school level mathematics
    $R_3$ - Vedic Mathematics is secondary school level mathematics
    $R_4$ - Vedic Mathematics is high school level mathematics
    $R_5$ - Nil (No use in Vedic Mathematics education)
    $R_6$ - Hindutva imposition through Vedic Mathematics
    $R_7$ - Imposition of Brahminism and caste systems
    $R_8$ - Training young minds in religion without their knowledge
    $R_9$ - Has some Vedic value
    $R_{10}$ - Has no mathematical value
    $R_{11}$ - It has neither Vedic value nor Mathematical value
    $R_{12}$ - It has Hindutva / religious fundamentalist agenda
    $R_{13}$ - Absolutely no educational value only religious value

We make use of the FRM model to analyze this problem.

These 13 nodes / attributes are taken as the nodes of the range space. All these nodes in the domain and range space are self-explanatory so we have not described them. The following directed graph was given by the first expert.

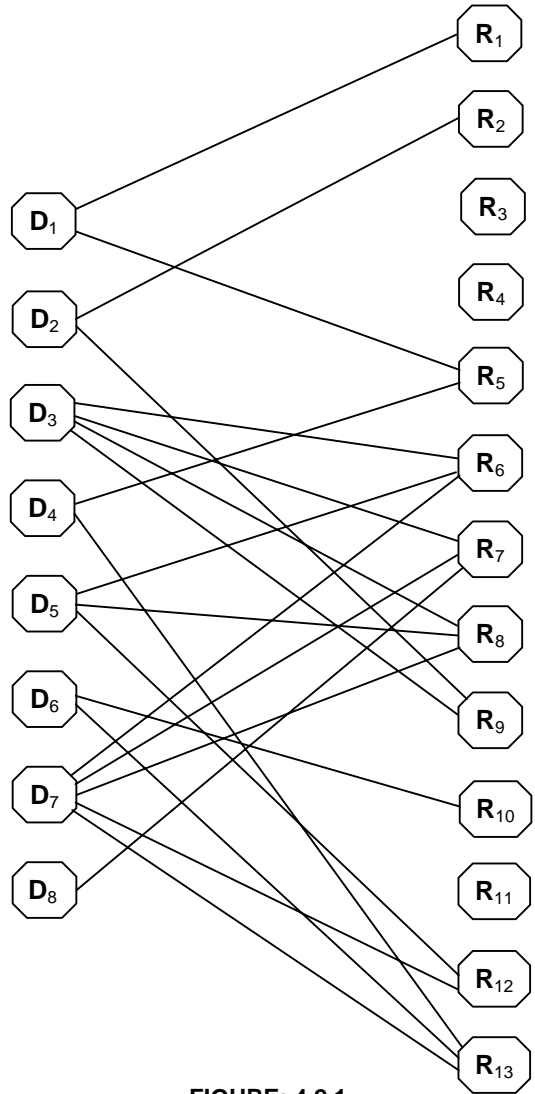

**FIGURE: 4.2.1**

The expert is a teacher working in a school run by pro-religious revivalist Hindutva trust. We use the directed graph of the Fuzzy Relational Maps and obtain the 8 × 13 connection matrix. The attributes related with the domain space are along the rows of the matrix and that of the range space attributes are taken along the column. Let us denote the 8 × 13 matrix by $M_1$.

$$M_1 = \begin{array}{c} \\ D_1 \\ D_2 \\ D_3 \\ D_4 \\ D_5 \\ D_6 \\ D_7 \\ D_8 \end{array} \begin{array}{c} R_1\ R_2\ R_3\ R_4\ R_5\ R_6\ R_7\ R_8\ R_9\ R_{10}\ R_{11}\ R_{12}\ R_{13} \\ \begin{bmatrix} 1 & 0 & 0 & 0 & 1 & 0 & 0 & 0 & 0 & 0 & 0 & 0 & 0 \\ 0 & 1 & 0 & 0 & 0 & 0 & 0 & 0 & 1 & 0 & 0 & 0 & 0 \\ 0 & 0 & 0 & 0 & 0 & 1 & 1 & 1 & 1 & 0 & 0 & 0 & 0 \\ 0 & 0 & 0 & 0 & 1 & 0 & 0 & 0 & 0 & 0 & 0 & 0 & 1 \\ 0 & 0 & 0 & 0 & 0 & 1 & 0 & 1 & 0 & 0 & 0 & 1 & 0 \\ 0 & 0 & 0 & 0 & 0 & 0 & 0 & 0 & 0 & 1 & 0 & 0 & 1 \\ 0 & 0 & 0 & 0 & 0 & 1 & 1 & 1 & 0 & 0 & 0 & 1 & 1 \\ 0 & 0 & 0 & 0 & 0 & 0 & 1 & 0 & 0 & 0 & 0 & 0 & 0 \end{bmatrix} \end{array}$$

Now we study the effect of the state vector X from the domain space in which, only the node $D_4$ alone is in the ON state and all other nodes are in the OFF state. Now we study the effect of X = (0 0 0 1 0 0 0 0) on the dynamical system $M_1$.

$$\begin{array}{lcll} XM_1 & = & (0\ 0\ 0\ 0\ 1\ 0\ 0\ 0\ 0\ 0\ 0\ 0\ 1) & = Y \\ YM_1^T & = & (1\ 0\ 0\ 1\ 0\ 1\ 1\ 0) & = X_1\ \text{(say)} \\ X_1M_1 & \rightarrow & (1\ 0\ 0\ 0\ 1\ 1\ 1\ 1\ 0\ 1\ 0\ 1\ 1) & = Y_1\ \text{(say)} \end{array}$$

'→' denotes after updating and thresholding the resultant vector got from $X_1M_1$. Now

$$\begin{array}{lcll} Y_1M_1^T & \rightarrow & (1\ 0\ 1\ 1\ 1\ 1\ 1\ 1) & = X_2\ \text{(say)} \\ X_2M_1 & \rightarrow & (1\ 0\ 0\ 0\ 1\ 1\ 1\ 1\ 1\ 1\ 0\ 1\ 1) & = Y_2\ \text{(say)} \\ Y_2M_1^T & = & (1\ 1\ 1\ 1\ 1\ 1\ 1\ 1) & = X_3\ \text{(say)} \\ X_3M_1 & = & (1\ 1\ 0\ 0\ 1\ 1\ 1\ 1\ 1\ 1\ 0\ 1\ 1) & = Y_3\ \text{(say)} \\ Y_3M_1^T & \rightarrow & (1\ 1\ 1\ 1\ 1\ 1\ 1\ 1) & = X_4\ (= X_3). \end{array}$$

Thus the hidden pattern of the dynamical system given by vector X = (0 0 0 1 0 0 0 0) is a binary pair which is a fixed binary pair of the dynamical system $M_1$.

When only the node ($D_4$) i.e. use of Vedic Mathematics in higher learning is on we see all the nodes in the domain space come to ON state.

In the range space all nodes except the nodes Vedic Mathematics is secondary school education level node $R_3$, Vedic Mathematics is high school level node $R_4$ and $R_{11}$, it has neither Vedic value nor mathematical value alone remain in the OFF state. The binary pair is given by {(1 1 1 1 1 1 1 1), (1 1 0 0 1 1 1 1 1 0 1 1)}.

Suppose we consider a state vector Y = (0 0 0 0 0 0 1 0 0 0 0 0) where only the node $R_7$ is in the ON state and all other nodes are in the OFF state; Y is taken from the range space. We study the effect of Y on the dynamical system $M_1$.

$$
\begin{aligned}
YM_1^T &= (0\ 0\ 1\ 0\ 0\ 0\ 1\ 1) &= X\ (\text{say}) \\
XM_1 &\rightarrow (0\ 0\ 0\ 0\ 0\ 1\ 1\ 1\ 1\ 0\ 0\ 1\ 1) &= Y_1\ (\text{say}) \\
Y_1 M_1^T &\rightarrow (0\ 1\ 1\ 1\ 1\ 1\ 1\ 1) &= X_1\ (\text{say}) \\
X_1 M_1 &\rightarrow (0\ 1\ 0\ 0\ 1\ 1\ 1\ 1\ 1\ 1\ 0\ 1\ 1) &= Y_2\ (\text{say}) \\
Y_2 M_1^T &= (1\ 1\ 1\ 1\ 1\ 1\ 1\ 1) &= X_2\ (\text{say}) \\
X_2 M_1 &\rightarrow (1\ 1\ 0\ 0\ 1\ 1\ 1\ 1\ 1\ 1\ 0\ 1\ 1) &= Y_3\ (\text{say}) \\
Y_3 M_1^T &= (1\ 1\ 1\ 1\ 1\ 1\ 1\ 1) &= X_3\ (= X_2).
\end{aligned}
$$

Thus resultant of the state vector Y = (0 0 0 0 0 0 1 0 0 0 0 0) is the binary pair which is a fixed point given by {(1 1 0 0 1 1 1 1 1 1 0 1 1), (1 1 1 1 1 1 1 1)} when only the node $R_7$ in the range space is in the ON state and all other nodes were in the OFF state.

Thus we can work with the ON state of any number of nodes from the range space or domain space and find the resultant binary pair and comment upon it (interpret the resultant vector).

Next we take the 2$^{nd}$ expert as a retired teacher who is even now active and busy by taking coaching classes for 10$^{th}$, 11$^{th}$ and 12$^{th}$ standard students. He says in his long span of teaching for over 5 decades he has used several arithmetical means which are shortcuts to multiplication, addition and division. He says that if he too had ventured he could have written a book like

Vedic Mathematics of course baring the sutras. We now give the directed graph given by him.

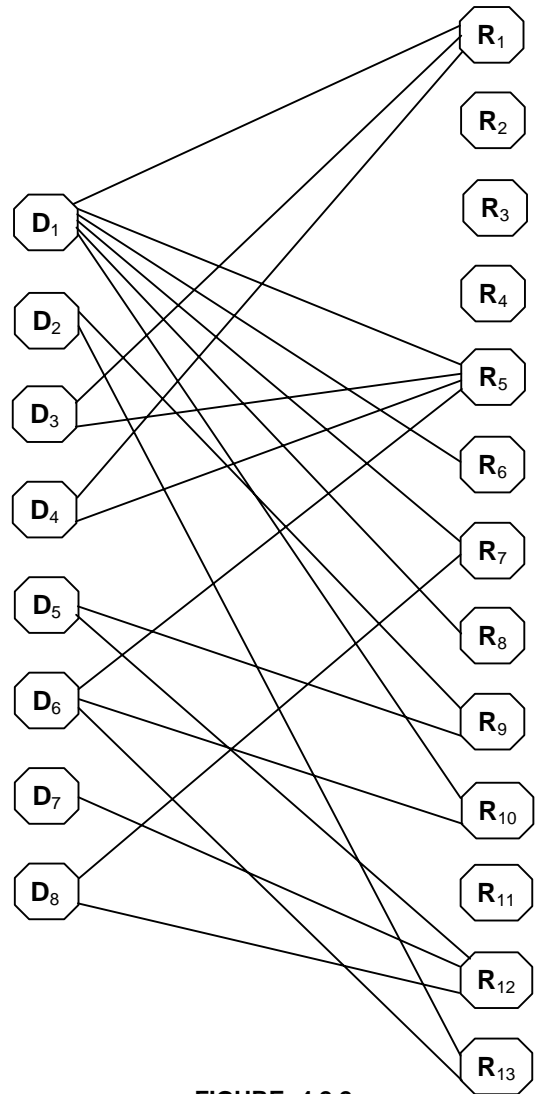

**FIGURE: 4.2.2**

Now using the directed graph given by him we have obtained the fuzzy relational matrix $M_2$.

$$M_2 = \begin{array}{c} \\ D_1 \\ D_2 \\ D_3 \\ D_4 \\ D_5 \\ D_6 \\ D_7 \\ D_8 \end{array} \begin{array}{c} R_1\ R_2\ R_3\ R_4\ R_5\ R_6\ R_7\ R_8\ R_9\ R_{10}\ R_{11}\ R_{12}\ R_{13} \\ \begin{bmatrix} 1 & 0 & 0 & 0 & 1 & 1 & 1 & 1 & 0 & 1 & 0 & 0 & 0 \\ 0 & 0 & 0 & 0 & 0 & 0 & 0 & 0 & 1 & 0 & 0 & 0 & 1 \\ 1 & 0 & 0 & 0 & 1 & 0 & 0 & 0 & 0 & 0 & 0 & 0 & 0 \\ 1 & 0 & 0 & 0 & 1 & 0 & 0 & 0 & 0 & 0 & 0 & 0 & 0 \\ 0 & 0 & 0 & 0 & 0 & 0 & 0 & 1 & 0 & 0 & 1 & 0 \\ 0 & 0 & 0 & 0 & 1 & 0 & 0 & 0 & 0 & 1 & 0 & 0 & 1 \\ 0 & 0 & 0 & 0 & 0 & 0 & 0 & 0 & 0 & 0 & 0 & 1 & 0 \\ 0 & 0 & 0 & 0 & 0 & 0 & 1 & 0 & 0 & 0 & 0 & 1 & 0 \end{bmatrix} \end{array}$$

Now this expert wants to study the effect of $X = (0\ 0\ 0\ 1\ 0\ 0\ 0\ 0)$ on $M_2$.

$$\begin{aligned} XM_2 &= (1\ 0\ 0\ 0\ 1\ 0\ 0\ 0\ 0\ 0\ 0\ 0\ 0) \\ &= Y_1\ (\text{say}) \\ Y_1 M_2^T &\to (1\ 0\ 1\ 1\ 0\ 1\ 0\ 0) \\ &= X_1\ (\text{say}) \\ X_1 M_2 &\to (1\ 0\ 0\ 0\ 1\ 1\ 1\ 1\ 0\ 1\ 0\ 0\ 1) \\ &= Y_2\ (\text{say}) \\ Y_2 M_2^T &= (1\ 1\ 1\ 1\ 0\ 1\ 0\ 1) \\ &= X_2\ (\text{say}) \\ X_2 M_2 &\to (1\ 0\ 0\ 0\ 1\ 1\ 1\ 1\ 1\ 0\ 1\ 1). \end{aligned}$$

Thus the resultant is a fixed point given by the binary pair $\{(1\ 1\ 1\ 1\ 0\ 1\ 0\ 1), (1\ 0\ 0\ 0\ 1\ 1\ 1\ 1\ 1\ 0\ 1\ 1)\}$.

Now we consider the same state vector of the range space given by the first expert.

Let $Y_1 = (0\ 0\ 0\ 0\ 0\ 0\ 1\ 0\ 0\ 0\ 0\ 0\ 0)$.

Now we study the effect of Y on the dynamical system $M_2$.

$$
\begin{array}{rcl}
Y_1 M_2^T &=& (1\ 0\ 0\ 0\ 0\ 0\ 0\ 1) \\
&=& X_1 \text{ (say)} \\
X_1 M_2 &\to& (1\ 0\ 0\ 0\ 1\ 1\ 1\ 1\ 0\ 1\ 0\ 1\ 0) \\
&=& Y_1 \text{ (say)} \\
Y_2 M_2^T &=& (1\ 0\ 1\ 1\ 1\ 1\ 1\ 1) \\
&=& X_2 \text{ (say)} \\
X_2 M_2 &\to& (1\ 0\ 0\ 0\ 1\ 1\ 1\ 1\ 1\ 1\ 0\ 1\ 1) \\
&=& Y_2 \text{ (say)} \\
Y_2 M_2^T &=& (1\ 1\ 1\ 1\ 1\ 1\ 1\ 1) \\
&=& X_3\ (=X_2) \\
X_3 M_2 &\to& (1\ 0\ 0\ 0\ 1\ 1\ 1\ 1\ 1\ 1\ 0\ 1\ 1) \\
&=& Y_3 \text{ (say)}.
\end{array}
$$

Thus resultant is a fixed binary pair given by {(1 0 0 0 1 1 1 1 1 1 0 1 1), (1 1 1 1 1 1 1 1)}. From the teacher's view-point we see that they are least bothered about the primary level or secondary level or high school level in Vedic Mathematics or whether it has a Vedic value or any mathematical value because what they are interested is whether Vedic Mathematics has no mathematical value or even any true Vedic value, that is why they remain zero at all stages. What is evident is that the introduction of Vedic Mathematics has ulterior motives and it only has a Hindutva background that is why in the dynamical system itself all these terms $R_2$, $R_3$, $R_4$ and $R_{11}$ are zero.

Now we have used several other experts to derive the conclusions using the C program given in [143].

The set of experts were given an option to work with NRM described in section 3.5 of this book. Most of them were reluctant to work with it. Only seven of them gave the NRM for the same sets of attributes. All the seven of them gave the relation between the node $D_2$ and $R_{11}$ as I. Some gave $D_2$ with $R_{10}$ as I and some other gave $D_2$ with $R_9$ as I. All these NRMs were constructed and using these NRM connection neutrosophic matrices hidden patterns of the suggested ON state of nodes as given by the experts were found and included in the chapter 5.

## 4.3 Views of Parents about Vedic Mathematics

In this section we give the views of parents. The parents from whom we could get the views happened to be a very heterogeneous crowd. Some educated parents had some notion about Vedic Mathematics, whereas some did not know about it, some were unconcerned and so on. Already in the earlier chapter, we have given important views about Vedic Mathematics that were obtained from parents. We met over 120 parents. Some had in fact met us for getting our views about their child attending the Vedic Mathematics classes and the uses of Vedic Mathematics in their child's curriculum.

The consolidated views from discussions find its place in the last chapter on observations. A few factors worth mentioning are :

1. Most of the non-Brahmin parents felt their child was ill-treated in Vedic Mathematics classes on the basis of caste. They were discriminated by the Vedic Mathematics teachers and were called as idiots, brainless, dull-head and so on.

2. A few parents said the pangs of caste discrimination had ruined their child psychologically due to the Vedic Mathematics classes. As a result, some parents had got special request from the educational officers to permit their child to remain absent for these classes.

3. Even most of the Brahmin parents felt that the Vedic Mathematics classes was only waste of time and that their children were forced to recite certain sutras which was meaningless both mathematically and scientifically. They felt that there was no visible improvement in their child's mathematical skill or knowledge.

4. Some parents were ignorant of what was happening in Vedic Mathematics classes.

5. All of them uniformly felt that these classes were an additional monetary commitment and of no use to their children.

6. Most parents feel that their duty is over once they pay the fees and give them the required money for travel and food so were unaware about what was taught in Vedic Mathematics classes.

7. Some parents felt that the school administration was perfect so they made the child attend the classes in spite of their child's displeasure and dislike in doing it; only in our discussion they found that they should have listened to their child and in fact some parents even said that this Vedic Mathematics classes have brought down the percentage of their marks in other subjects. They realized their ward had some mental conflicts due to different or discriminatory treatment in Vedic Mathematics classes. They openly repented that at the appropriate time they did not listen to their children.

8. Parents have been well informed by their children that, Vedic Mathematics classes were utter waste and the syllabus covered was very elementary. It is only the parents who failed to heed to the children because they were afraid to face any friction with the teachers or the authorities of the school. They were very apologetic towards their acts which they admitted during our open discussions (in several discussion the child was also present with the parents).

9. Some parents said that the Vedic Mathematics classes gave problems of primary school level and the recitation of sutras took their time and energy.

10. Some parents said "My child is a shy type. After coming to high school if they ask him to recite loudly some sutras which do not took like mathematical formula and that too not in English; it makes the teenagers feel bad. Some teachers punish them, that too like standing on the bench etc. Some teachers ask them to recite it individually; They feel so shy to pronounce meaningless Sanskrit words which is difficult to run smoothly through their mouth. For this act they become a laughing stock in the class and the Vedic Mathematics teachers take it as an insult and doubly punish the children."

11. Some of the uneducated and not-so-literate parents said, "after all my son is going to become a computer engineer, how is this sutra in Sanskrit going to help him?" The children say that the mathematical content is elementary arithmetic of primary level. One lady said, "they waste our money and our children's time by these Vedic Mathematics classes" though she has only studied up to $5^{th}$ standard. The questions she put to us about Vedic Mathematics was very pertinent. She laughed and said, "in temples they blabber something like this and get money, that too like beggars in a plate; now they have started to come to this school and get money in hundreds by saying some meaningless sutras." She further added that she was happy because her second son is studying in a Convent. She says in that school no such sloka-stuff is taught. Only after enquiring this, she put him in a different school. She says only Hindu schools teach Vedic Mathematics. Convents and Corporation schools or Government-run schools do not teach Vedic Mathematics. She says "I am uneducated. I want my children to get good education." She asked us, "Why is Vedic Mathematics having slokas? Are they training them as temple priests?" We have put this mainly to show how even uneducated parents take interest in their children's education!

12. Most of the parents said Vedic Mathematics teachers do not have tolerance or patience, they easily punish children for very simple things like laughing or not concentrating or attentive in the class by looking at the teachers. Only this atmosphere made the classes noisy, uncontrollable and unruly. The Vedic Mathematics teachers do not appear to be well-trained teachers. Some ask the students in Vedic

Mathematics classes whether they take bath daily and so on which is irrelevant, apart from being too personal.

13. Some parents said Vedic Mathematics teachers speak of epics and characters like Mahabaratha's 'Kamsa' and so on. They feel a mathematics class cannot have place for epics; why Krishna or Kamsa should come while teaching mathematics? One may adore Krishna, some other person may worship Kamsa it is after all individual freedom, choice and taste! No one should preach Hinduism in Mathematics class because there are Christian and Muslim boys who might feel offended! Also some teachers gave long lecture on Vedas and Vedic tradition, which they consider as the high heritage of Indians. Some parents said, "Are not Christians and Muslims living in India; Indians? Why did they become or converted to Christianity and Islam? They were humiliated and treated worse than animals by the Brahmins so to live and lead a life of self respect they sought Christianity or Islam." Some parents asked us, "if alone Christians and Muslims had not entered India; can ever a non-Brahmin dream of education?" They felt Vedic Mathematics was imposing brahminism i.e., casteism on children so they strongly objected to it. Some parents had already changed the school (and many had plans for changing their wards to a different school) because they felt it was unbearable to impose "Hindutva" in the name of Vedic Mathematics. (Several other charges were made which we have not given fully).

14. Vedic Mathematics classes had become the seed of discrimination on the basis of caste in schools! This was a view shared by non-Brahmin parents.

15. A tiny section of the educated parents said they have read the book on Vedic Mathematics and they had found it very elementary. Yet they felt that it was a powerful means of establishing the supremacy of the Aryans over the entire world. We wonder why they need mathematics to do this dirty trade?

16. Most of the parents felt it is fortunate that the Tamil Nadu state government has not made Vedic Mathematics as a part of the syllabi in schools because if this is imposed as in a few other Indian states, the school will be the breeding place of caste by birth, Aryan domination and so on.

17. Several parents said they wonder how these Brahmin use mathematics as a means to promote and spread "Hindutva" all over the world. One parent wondered why a Sankaracharya (Swamiji) of Puri mutt should be involved. Some people asked us, "Are they going to ultimately say that Vedic Mathematics is just like Vedas, so Sudras and Panchamas should not read mathematics?" But those who had read the book raised a point that the book has more ulterior motives than the elementary primary level mathematics displayed in it.

18. Uniformly, all parents appreciated the non-Hindutva schools that did not recommend Vedic Mathematics. They have fortunately not fallen a prey to this concept. However, they felt that because of extensive propaganda a few of the western schools have taken up Vedic Mathematics, but soon they too will realize the motivation behind the book. It is a mission to globalize 'Hindutva' and nothing more, they said.

19. This final point is not related with Vedic Mathematics but with the interrelation between parents and their children which is universally true. If this sort of relation continues in due course of time the bondage between parents with children would become very weak. The fault lies not with children but only with parents. We obtained this idea after our discussions with over 75 parents. Almost all the parents felt that their duty was over once they pay the fees to the children and provide them with all basic needs like transport, food and books. They fail to understand what the child needs is not all this, but above all these is their "time" that is they should make it a point to spend some time with their children finding their problems, progress and so on.

This should not be done sitting before a TV. or listening to news or music, this should be done whole-heartedly with no distractions. Most parents said or felt that their duty is over once they pay them fees and provide them their basic needs. This is not a recommendable attitude of the parents. So we requested parents to spend sometime individually on their children.

### 4.4 Views of Educationalists About Vedic Mathematics

We had discussed about Vedic Mathematics with judges, bank officers, vice-chancellors, directors, industrialists, engineers, doctors and others. We have categorized them under the broad title of 'educated elite'/ 'educationalists' because in the next section where the public have given their opinion many of them view it in the political angle, party angle and so on. Thus these educationalists have given their views on the social structure or changes that Vedic Mathematics could inculcate on the mindset of children (students), the psychological impact and so on.

They share the view that Vedic Mathematics may not only influence the students but to some extent may also strain the student-teacher relationship. Thus when we had to gather opinion it was more on why the Swamiji who said that it was just a simple arithmetic course to help students to do mathematical calculations mentally named it as Vedic Mathematics. Was the motivation behind it religious, casteist or both? Many questions were raised and several types of analysis were done. It was feared that such a topic would further kindle caste and discrimination at the very core, that too among students (who were just children.) If such discrimination is practiced, what progress will the nation make? Both caste superiority and caste discrimination are negative energies. Some of the respondents were worried about it and some felt that the way in which Vedic Mathematics was publicized was wrong.

Another perspective put forth by respondents was that Vedic Mathematics had become a moneymaking machine for some people. Thus many diverse opinions were received. As said by the first author's note it is pertinent that these people not only stayed at the putting questions about Vedic Mathematics but they also gave a lot of cooperation by sharing their thoughts and views. We can say with pride that over 90% of them had purchased the book on Vedic Mathematics and small group discussions were held with them only after they had thoroughly gone through the book. So this group seriously took up the topic of Vedic Mathematics and its ulterior motives in the context of society at large and the younger generation in particular.

They all uniformly feel that the 'Vedas' came into India only after the Aryans stepped into India. Their entry into India did more harm to the Indians (natives) than any good. A viewpoint shared by many members was that the Muslim conquest of India did not have such a bad impact because the Muslims treated the Indians as humans. But the Aryans followed their Vedas and treated the majority of the people as untouchables and un-seeables.

Secondly, the widely held opinion is that the British who ruled us were benevolent. The introduction of modern education system opened the doors to education for the so-called lower caste peoples, who were denied education according to the Vedas. By employing the native people as butlers, cooks, watchmen and helpers in their homes, they didn't practice discrimination. They dined equally with the Indians (natives). While the Aryans denied education and imposed curbs on the lower castes from becoming literate (lettered) the British helped the natives to become educated and self-sufficient. Thus within the span of a few generations, the indigenous people became more educated and more socially and economically powerful. The missionaries who came to India provided the people with good education. "The Aryans (Brahmins) who knew little English and little more educated than us tried to create a misunderstanding between us and the British" they said. They started to do this when they saw us getting education because they were not able to tolerate us getting educated and economically better. So they wanted the British to leave India. So they organized protests by falsely talking ill of the British on one side and on the other side, giving the feedback to the British that the lay people wanted freedom from them. Their double-stand ruined us because the politicians were power hungry and didn't bother about the well being of common man. Several people said that the Tamil Rationalist leader Periyar was very much against our independence because he rightly feared that we would be totally controlled and discriminated by the Brahmins. He declared 'independence day' to be a black day in the history of India. "Thus the Aryans crept in and the Vedas ruined us. We are now unaware of the real consequences that Vedic Mathematics has in store for us" they feared.

Now we have to be careful and above all rationalistic because it is not just mathematics but it is politically motivated and has several ulterior motives according to several of the respondents in this category. Thus we took their vital points about Vedic Mathematics as nodes / concepts.

$W_1$ - Vedic Mathematics: the ulterior motive is imposition of religion among the youth.
$W_2$ - Vedic Mathematics: ulterior motive is imposition of caste, based on birth (in Vedas) in the mindset of youth.
$W_3$ - Vedic Mathematics motivates the supremacy of Brahmins (Aryans) in the minds of the youth.
$W_4$ - Vedic Mathematics psychologically imposes Sanskrit as a better language in the minds of the youth.
$W_5$ - Vedic Mathematics tries to establish in the mindset of youth that all sciences and technologies are in Vedas!
$W_6$ - Vedic Mathematics develops complexes in young minds like caste difference and so on.
$W_7$ - Vedic Mathematics ruins the teacher-student relationship.
$W_8$ - Vedic Mathematics will develop the practice of caste differences (forms of untouchability) even among children.
$W_9$ - Vedic Mathematics has no real mathematical content.
$W_{10}$ - Vedic Mathematics has no real Vedic content.

$W_{11}$ - Vedic Mathematics is not an alternative for mathematics or arithmetic.
$W_{12}$ - Vedic Mathematics is a tool of the revivalist Hindutva.
$W_{13}$ - Vedic Mathematics is used to globalize Hindutva.
$W_{14}$ - Vedic Mathematics is an attempt to Brahminization of entire India.

We divided the educated respondents in this category into eight sub-categories. They are given below along with a brief description.

$E_1$ - People from the legal profession: includes judges, senior counsels, lawyers, professors who teach law and law college students.
$E_2$ - Educationalists: includes Vice chancellors, Directors, Principals, Headmasters and Headmistresses, non-mathematics teachers, professors in different fields, educational officers and inspectors of school etc.
$E_3$ - Technical Experts: this list includes engineers, technicians in different fields, all technically qualified persons like computer scientists, IT specialists, and teachers and researchers in those fields.
$E_4$ - Medical experts: Doctors, professors who teach in Medical colleges, Deans of Medical Colleges and researchers in medicine.
$E_5$ - Industrial experts: includes educated people who hold senior managerial positions in major industries.
$E_6$ - Government Staff: includes bank employees, government secretariat staff and clerical employees of government-run institutions.
$E_7$ - Businesspersons: includes people running private businesses like printing presses, magazines, export companies and so on.

E_8    -    Religious people: includes students of religion (theology) or philosophy who take up religious work, research scholars who study religion as their subject.

E_9    -    Social analysts: includes sociologists, social workers, teachers of social work, and others interested in studying social aspects and changes that influence the social setup.

Now the number of people in each group varied. The biggest group was educationalists numbering 41 and the least were the social scientists numbering only seven. Since all of them were educated, we placed before them the 14 conceptual nodes and asked them to give scores between 0 and 1. We took the groups and took their opinion on the 14 nodes. For the sake of uniformity if n people from a group gave the opinion we added the n terms against each node and divided it by n. This always gives a number between 0 and 1. Now taking along the rows the category people and along the columns the 14 concepts given by them on Vedic Mathematics we formed a $9 \times 14$ matrix which will be called as the New Fuzzy Dynamical System. Now using max-min operations we found the effect of any state vector on the dynamical system.

We had also explained to the groups about their values: when they give; zero, it suggests no influence, if they give positive small value say 0.01 it denotes a very small influence but something like 0.9 denotes a very large positive influence. We felt it difficult to educate all of them on the concept of negative, small negative and large negative values and so on. Therefore, we advised them to give values from 0 to 1.

Now we use all the experts opinion and have obtained the new fuzzy vector matrix M which we call as the New Fuzzy Dynamical System described in chapter 3 section 3.3. As most of the people gave the values only up to first decimal place we have worked with all the experts and have approximated the entries to first decimal place. Thus our dynamical system forms a fuzzy vector matrix with gradations. M is a $9 \times 14$ matrix with entries from the closed interval [0, 1]. Expert opinion will be given in the form of fit vectors that we have described in [68].

Using the experts opinion we find the resultant state vector, using the new dynamical system M.

$$\begin{bmatrix} 0.8 & 0.7 & 0.9 & 0.6 & 0 & 0.6 & 0.8 & 0.7 & 0.0 & 0 & 0 & 0.6 & 0.8 & 0.7 \\ 0.6 & 0.8 & 0.3 & 0.7 & 0.8 & 0.2 & 0.6 & 0 & 0.9 & 0 & 0.8 & 0.3 & 0.2 & 0.6 \\ 0.7 & 0.6 & 0.8 & 0 & 0.9 & 0 & 0 & 0.6 & 0.6 & 0 & 0.7 & 0.6 & 0.6 & 0.7 \\ 0.6 & 0.7 & 0.6 & 0.8 & 0.6 & 0.6 & 0.7 & 0 & 0 & 0 & 0.7 & 0.5 & 0.5 & 0.8 \\ 0.6 & 0.7 & 0.6 & 0.5 & 0.5 & 0.6 & 0.8 & 0.7 & 0 & 0 & 0 & 0.7 & 0.8 & 0.9 \\ 0.5 & 0.8 & 0.6 & 0.6 & 0.4 & 0.3 & 0.9 & 0.8 & 0 & 0 & 0 & 0.6 & 0.7 & 0.8 \\ 0.6 & 0.6 & 0.7 & 0.8 & 0 & 0.5 & 0.8 & 0.7 & 0 & 0 & 0 & 0.7 & 0.6 & 0.5 \\ 0.7 & 0.8 & 0.6 & 0.5 & 0.9 & 0.6 & 0.7 & 0.6 & 0 & 0 & 0 & 0.7 & 0.6 & 0.6 \\ 0.6 & 0.5 & 0.6 & 0.8 & 0.7 & 0.6 & 0.5 & 0.2 & 0 & 0 & 0 & 0.8 & 0.6 & 0.5 \end{bmatrix}$$

Suppose B = (1 0 0 0 0 1 0 0 0) is the state vector given by the expert. To find the effect of B on the new dynamical system M.

BM     =     $\max_{ij} \min (b_j, m_{ij})$
         =     (0.8, 0.8, 0.9, 0.6, 0.4, 0.6, 0.9, 0.8, 0, 0, 0, 0.6, 0.8, 0.8)
         =     A.

Now

M $A^T$     =     max min {$m_{ij}, a_i$}
         =     (0.9, 0.8, 0.8, 0.8, 0.8, 0.8, 0.8, 0.8, 0.6)
         =     $B_1$ (say).

B M     =     (0.8, 0.8, 0.9, 0.8, 0.8, 0.6, 0.8, 0.8, 0.8, 0, 0.8, 0.7, 0.8, 0.8)
         =     $A_1$ (say).

M $A^T_1$     =     (0.9, 0.8, 0.8, 0.8, 0.8, 0.8, 0.8, 0.8, 0.8,)
         =     $B_2$ (say).

$B_2$ M     =     (0.8, 0.8, 0.9, 0.8, 0.8, 0.6, 0.8, 0.8, 0.8, 0, 0.8, 0.8, 0.8, 0.8)
         =     $A_2$ (say).

$MA^T_2$ = $B_3$
     = (0.9, 0.8, 0.8, 0.8, 0.8, 0.8, 0.8, 0.8, 0.8)
     = $B_2$.

Thus we arrive at a fixed point. When the views of the educated from the legal side (1) together with the secretarial staff views (6) are given by the expert for analysis we see that they cannot comment about the Vedic content, so the node 10 is zero. However, to ones surprise they feel that Vedic Mathematics has no mathematical value because that node takes the maximum value 0.9. Further the study reveals that all others also feel the same, the nodes related to everyone is 0.8.

Now the expert wishes to work with the nodes 1, 3, 9 and 14 to be in the ON state. Let the fuzzy vector related with it be given by

A = (1 0 1 0 0 0 0 0 1 0 0 0 0 1).

The effect of A on the new dynamical system M is given by

$MA^T$ = (0.9, 0.9, 0.8, 0.8, 0.9, 0.8, 0.7, 0.7, 0.6)
     = B (say).

BM = (0.8, 0.8, 0.9, 0.8, 0.8, 0.6, 0.8, 0.8, 0.8, 0, 0.8, 0.7, 0.8, 0.8)
   = $A_1$ (say).

$MA^T_1$ = (0.9, 0.8, 0.8, 0.8, 0.8, 0.8, 0.8, 0.8, 0.8)
       = $B_1$ (say).

$B_1 M$ = (0.8, 0.8, 0.9, 0.8, 0.8, 0.6, 0.8, 0.8, 0.8, 0, 0.8, 0.8, 0.8, 0.8)
       = $A_2$ (say).

$MA^T_2$ = (0.9, 0.8, 0.8, 0.8, 0.8, 0.8, 0.8, 0.8, 0.8)
       = $B_2$ (say) = $B_1$.

Now $B_2 = B_1$. Thus we arrive at a fixed binary pair which says that when nodes 1, 3, 9 and 14 alone are in the ON state all nodes in B get the same value 0.8 except the node 1 which gets 0.9. There by showing that all educated groups feel and think alike about Vedic Mathematics. Further we see the views held as same as before with $10^{th}$ node, which comes as 0.

Now the expert wants to analyze only the views held by the educated religious people i.e. only the node 8 is in the ON state in the state vector B and all other nodes are in the off state, i.e. B = (0 0 0 0 0 0 0 1 0).

BM  = (0.7, 0.8, 0.6, 0.5, 0.9, 0.6, 0.7, 0.6, 0, 0, 0, 0.7, 0.6, 0.6)
    = A (say).

$MA^T$ = (0.7, 0.8, 0.9, 0.7, 0.7, 0.8, 0.7, 0.8, 0.7)
    = $B_1$ (say).

$B_1 M$ = (0.7, 0.8, 0.8, 0.7, 0.9, 0.6, 0.8, 0.8, 0.8, 0, 0.8, 0.7, 0.7, 0.8)
    = $A_1$ (say).

$M A^T_1$ = (0.8, 0.8, 0.9, 0.8, 0.8, 0.8, 0.8, 0.9, 0.7)
    = $B_2$ (say).

$B_2 M$ = (0.8, 0.8, 0.8, 0.8, 0.8, 0.6, 0.8, 0.8, 0.8, 0, 0.8, 0.7, 0.8, 0.8)
    = $A_2$ (say)

$M A^T_2$ = (0.8, 0.8, 0.9, 0.8, 0.8, 0.8, 0.8, 0.9, 0.8)
    = $B_3$ (say).

$B_3 M$ = (0.8, 0.8, 0.8, 0.8, 0.8, 0.6, 0.8, 0.8, 0.8, 0, 0.8, 0.8, 0.8, 0.8)
    = $A_3$ (say).

$M A^T_3$ = (0.8, 0.8, 0.9, 0.8, 0.8, 0.8, 0.8, 0.9, 0.8)
    = $B_4 = B_3$.

Thus we arrive at the fixed point. Everybody is of the same view as the religious people. One can derive at any state vector and draw conclusions. Further we see they do not in general differ in grades because they hold the even same degree of opinion about Vedic Mathematics. Thus we have given in the last chapter on observations about the results worked out using the new dynamical system. Further we cannot dispose of with the resultant vector for they hold a high degree viz. 0.8, in the interval [0, 1]. Also we see that the educated masses as a whole did not want to comment about the Vedic content in Vedic Mathematics.

Now we asked the experts if they thought there was any relation between concepts that cannot be given value from [0, 1] and remained as an indeterminate relationship. Some of them said yes and their opinion alone was taken and the following new fuzzy neutrosophic dynamical system $M_n$ was formed.

$$\begin{bmatrix} 0.8 & 0.7 & 0.9 & 0.6 & 0 & 0.6 & 0.8 & 0.7 & 0 & I & 0 & 0.6 & 0.8 & 0.7 \\ 0.6 & 0.8 & 0.3 & 0.7 & 0.8 & 0.2 & 0.6 & 0 & 0.9 & 0 & 0.8 & 0.3 & 0.2 & 0.6 \\ 0.7 & 0.6 & 0.8 & 0 & 0.9 & I & 0 & 0.6 & 0.6 & 0 & 0.7 & 0.6 & 0.6 & 0.7 \\ 0.6 & 0.7 & 0.6 & 0.8 & 0.6 & 0.6 & 0.7 & 0 & 0 & 0 & 0.7 & 0.5 & 0.5 & 0.8 \\ 0.6 & 0.7 & 0.6 & 0.5 & 0.5 & 0.6 & 0.8 & 0.7 & 0 & I & 0 & 0.7 & 0.8 & 0.9 \\ 0.5 & 0.8 & 0.6 & 0.6 & 0.4 & 0.3 & 0.9 & 0.8 & 0 & 0 & 0 & 0.6 & 0.7 & 0.8 \\ 0.6 & 0.6 & 0.7 & 0.8 & 0 & 0.5 & 0.8 & 0.7 & I & 0 & 0 & 0.7 & 0.6 & 0.5 \\ 0.7 & 0.8 & 0.6 & 0.5 & 0.9 & 0.6 & 0.7 & 0.6 & 0 & 0 & 0 & 0.7 & 0.6 & 0.6 \\ 0.6 & 0.5 & 0.6 & 0.8 & 0.7 & 0.6 & 0.5 & 0.2 & 0 & I & 0 & 0.8 & 0.6 & 0.5 \end{bmatrix}$$

As in case of the new dynamical system we worked with the state vectors given by the experts. They felt that because they were unaware of the Vedic language Sanskrit and the Vedas they restrained from commenting about it. Uniformly they shared the opinion that teaching such a subject may develop caste differences among children had a node value of 0.6 only.

### 4.5. Views of the Public about Vedic Mathematics

When we spoke about Vedic Mathematics to students, teachers, educated people and parents we also met several others who were spending their time for public cause, some were well educated, some had a school education and some had no formal education at all. Apart from this, there were many N.G.O volunteers and social workers and people devoted to some social cause. So, at first we could not accommodate them in any of the four groups. But they were in the largest number and showed more eagerness and enthusiasm than any other group to discuss about Vedic Mathematics and its ulterior motives. So, by the term 'public' we mean only this group which at large has only minimum or in some cases no overlap with the other four groups.

Here it has become pertinent to state that they viewed Vedic Mathematics entirely in a different angle: not as mathematics or as Vedas; but as a tool of the revivalist, Hindu-fundamentalist forces who wanted to impose Aryan supremacy. Somehow, majority of them showed only dislike and hatred towards Vedic Mathematics. The causes given by them will be enlisted and using experts' opinions, fuzzy mathematical analysis will be carried out and the observations would be given in the last chapter. Several of these people encouraged us to write this book.

The first edition of the book on Vedic Mathematics was published in 1965, five years after the death of its author, His Holiness Jagadguru Sankaracharya of Puri. The author says he had written sixteen volumes and his disciple lost them. So in this book he claims to have put the main gist of the 16 volumes.

The book remained in cold storage for nearly two decades. Slowly it gathered momentum. For instance, S.C.Sharma, Ex-Head of the Department of Mathematics, NCERT [National Council of Educational Research and Training—which formulates the syllabus for schools all over the nation] spoke about this book in Mathematics Today September 1986. Some of the excerpts from S.C.Sharma are, "This book brings to light how great and true knowledge is born of initiation, quite different from modern western methods. The ancient Indian method and its secret techniques are examined and shown to be capable of solving various problems of mathematics…"

The volume more a 'magic is the result of notational visualization of fundamental mathematical truths born after eight years of highly concentrated endeavour of Jagadguru Sri

Bharati…. The formulae given by the author from Vedas are very interesting and encourage a young mind for learning mathematics as it will not be a bugbear to him".

Part of this statement also appeared as a blurb on the back cover of Vedic Mathematics (Revised Ed. 1992) [51].

It is unfortunate that just like the 16 lost volumes of the author, the first edition [which they claim to have appeared in 1965] is not available. We get only the revised edition of 1992 and reprints have been made in the years 1994, 1995, 1997, 1998, 2000 and 2001. The people we interviewed in this category say that just like the Vedas, this book has also undergone voluminous changes in its mathematical contents. Several of the absurdities have been corrected. The questions and views put forth to us by the respondents are given verbatim. First, they say a responsible person like S.C. Sharma, who served, as Head of Department of Mathematics in the NCERT cannot use words like "magic" in the context of mathematics. Can mathematics be magic? It is the most real and accurate science right from the school level.

Secondly, they heavily criticized the fact that it took eight long years to publish such an elementary arithmetic mathematics book. Further they are not able to understand why S.C.Sharma uses the phrase "secret techniques" when westerners are so open about any discovery. If the discovery from Vedas had been worthwhile they would not keep it as a secret. The term "secret techniques" itself reveals the standard of the work.

One may even doubt whether these terms have any ulterior motives because the standard of Vedic Mathematics is itself just primary school level arithmetic. That is why, most people in this category held that only after the rightwing and revivalist Bharatiya Janata Party (BJP) picked up some political status in India, Vedic Mathematics became popular. It has achieved this status in one and a half decades. Because of their political power, they have gone to the extent of prescribing Vedic Mathematics in the syllabi of all schools in certain states ruled by BJP and this move is backed by the RSS (Rashtriya Swayamsevak Sangh) and VHP (Vishwa Hindu Parishad) (Hindu fanatic groups). They have their own vested interests for upholding and promoting Vedic Mathematics. The very act of waiting for the fanatic Hindutva Government to come to power and then forcing the book on innocent students shows that this Vedic Mathematics does not have any mathematical content or mathematical agenda but is the only evidence of ulterior motives of Hindutvaizing the nation.

It is a means to impose Brahmin supremacy on the non-Brahmins and nothing more. Further they added that 16 sutras said in Sanskrit are non-mathematical. One of the interviewed respondents remarked that it was a duty of the educated people to hold awareness meetings to let the masses know the ulterior motives of the Brahmins who had come to India as migrants through the Khyber Pass and now exploit the natives of the land. Discussion and debates over Vedic Mathematics will give us more information about the ulterior motives. It is apparently an effort to globalize Hindutva. All of them asked a very pertinent question: when the Vedas denied education to the non-Brahmins how can we learn Vedic Mathematics alone? They said one point of the agenda is that they have made lots of money by selling these books at very high prices. Moreover, people look at Vedic Mathematics as "magic" or "tricks" and so on.

They don't view Vedic Mathematics as mathematics, an organized or logical way of thinking. One respondent said, "They have done enough 'magic' and 'tricks' on us; that is why we are in this status. Why should a person with so high a profile use 'magic' to teach mathematics that too to very young children? These simple methods of calculations were taught in schools even before the advent of Vedic Mathematics. Each mathematics teacher had his own ingenious way of solving simple arithmetic problems. All the cunningness lies in the title itself: "Vedic Mathematics."

They said that when a person dies, a Brahmin carries out the death ceremony and rituals because he claims only he has the magical power to send the dead to heavens. So soon after the death he performs some rituals (collects money, rice and other things depending on the economic status of the dead). Not only this after 16 days he once again performs the ritual for the dead saying that only when he throws the rice and food in the sky it reaches them! Instead of stopping with this, he performs the same sort of ritual for the same dead person on every anniversary of the death.

Now in Vedas, it is said that after his death a man is reborn, he may be reborn as a bird or animal or human depending on the karma (deeds) of his past life. So according to this Brahmin theory, the dead for whom we are performing rituals might already living as a animal or human then what is the necessity we should perform yearly rituals and 'magic' for the soul of the dead to be at peace when it is already living as some other life form? So, they say that the Vedas are full of lies and rubbish with no rhyme or reason. A few points put by them in common are taken up as the chief concepts to analyze the problem.

Now we proceed on to enlist the main points given by them.

1. When they claim Vedic Mathematics to be a 'magic', it has more ulterior motives behind it than mathematics.

2. Vedic Mathematics uses 'tricks' to solve the problems – "tricks" cannot be used to solve all mathematical problems. Any person with some integrity never uses tricks. They may use tricks in "circus" or "street plays" to attract public and get money. Children cannot be misled by these tricks in their formative age, especially about sciences like mathematics that involves only truth.

3. Vedic Mathematics speaks of sutras not formulae but some Sanskrit words or phrases. This has the hidden motive of imposing caste and discrimination; especially birth-based discrimination of caste in the minds of youth. In fact Swami Vivekananda said that most of the caste discriminations and riots are due to Sanskrit which is from the north. If the Sanskrit books and the literature were lost it would certainly produce peace in the nation he says. He feels Sanskrit is the root cause of all social inequalities and problems in the south.

4. The very fact the Christian and Muslim educational institutions do not use Vedic Mathematics shows its standard and obvious religious motivation!

5. It has a pure and simple Hindutva agenda (the first page of the books I and II of Vedic Mathematics in Tamil is evidence for this). [85-6]

6. It is a means to globalize Hindutva.

7. It is a means to establish Aryan supremacy.

8. Vedic Mathematics is used only to disturb young non-Brahmin minds and make them accept their inferiority over the Brahmins.

9. It is more a political agenda to rule the nation by indoctrination and if Sanskrit literature were lost it would certainly produce peace in the nation.

These concepts are denoted by $P_1$ to $P_9$.

From several factors they gave us, we took these nine concepts after discussion with few experts. Further we had a problem on who should be an expert. If a person from other group were made an expert it would not be so proper, so we chose only members of this group to be the experts and chose the simple Fuzzy Cognitive Maps (FCMs) to be the model because they can give the existence or the nonexistence of a relation together with its influence.

So we would be using only simple FCMs and NCMs to analyze the problem.

Since the data used also is only an unsupervised one we are justified in using FCMs. Now using the 9 nodes we obtain the directed graph using the expert 1 who is a frontline leader of a renowned Dravidian movement.

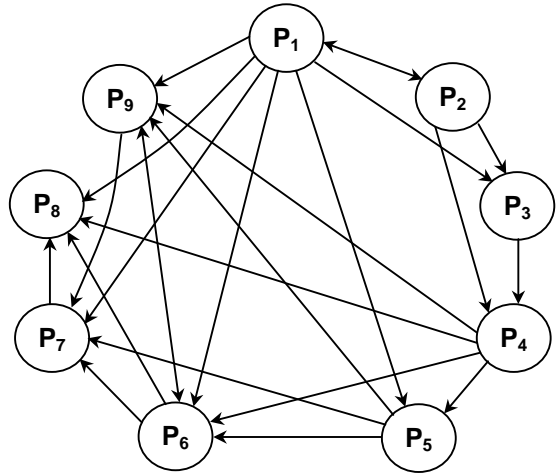

FIGURE 4.5.1

Using the directed graph given by the first expert we have the following relational matrix. Let $M_1$ denote the $9 \times 9$ fuzzy relational matrix.

$$M_1 = \begin{array}{c} \\ P_1 \\ P_2 \\ P_3 \\ P_4 \\ P_5 \\ P_6 \\ P_7 \\ P_8 \\ P_9 \end{array} \begin{array}{c} \begin{array}{ccccccccc} P_1 & P_2 & P_3 & P_4 & P_5 & P_6 & P_7 & P_8 & P_9 \end{array} \\ \left[ \begin{array}{ccccccccc} 0 & 1 & 1 & 0 & 1 & 1 & 1 & 1 & 1 \\ 1 & 0 & 1 & 1 & 0 & 0 & 0 & 0 & 0 \\ 0 & 0 & 0 & 1 & 0 & 0 & 0 & 0 & 0 \\ 0 & 0 & 0 & 0 & 1 & 1 & 0 & 1 & 1 \\ 0 & 0 & 0 & 0 & 0 & 1 & 1 & 0 & 1 \\ 0 & 0 & 0 & 0 & 0 & 0 & 1 & 1 & 1 \\ 0 & 0 & 0 & 0 & 0 & 0 & 0 & 1 & 0 \\ 0 & 0 & 0 & 0 & 0 & 0 & 0 & 0 & 0 \\ 0 & 0 & 0 & 0 & 0 & 1 & 1 & 0 & 0 \end{array} \right] \end{array}$$

Suppose the expert wants to study the state vector X when only the node 6 i.e. the globalization of Hindutva is the agenda of Vedic Mathematics is in the ON state and all other nodes are in the OFF state i.e.

$$X = (0\ 0\ 0\ 0\ 0\ 1\ 0\ 0\ 0);$$

Now the effect of X on the dynamical system $M_1$, is given by

$$XM_1 \rightarrow (0\ 0\ 0\ 0\ 0\ 1\ 1\ 1\ 1)$$
$$= X_1 \text{ (say)}$$

Now

$$X_1M_1 \rightarrow (0\ 0\ 0\ 0\ 0\ 1\ 1\ 1\ 1)$$
$$= X_2 = X_1.$$

Thus the hidden pattern of the state vector X gives a fixed point, which expresses, when the node globalization of Hindutva is the agenda of Vedic Mathematics alone is in the ON state we see the resultant is a fixed point and it makes nodes 7, 8, and 9 to ON state i.e. Vedic Mathematics establishes Aryan supremacy, Vedic Mathematics disturbs the young non-Brahmin minds and make them accept their inferiority over the Brahmins and Vedic Mathematics is more a political agenda to rule the nation.

Now the expert wants to study the effect of the node (1) i.e. Vedic Mathematics claims to be a 'magic' and this has ulterior motives than of mathematics; and all other nodes are in the OFF state. To study the effect of $Y = (1\ 0\ 0\ 0\ 0\ 0\ 0\ 0\ 0)$ on the dynamical system $M_1$.

$$YM_1 = (0\ 1\ 1\ 0\ 1\ 1\ 1\ 1\ 1)$$

after updating and thresholding we get

$$Y_1 = (1\ 1\ 1\ 0\ 1\ 1\ 1\ 1\ 1)$$
$$Y_1M_1 \rightarrow (1\ 1\ 1\ 0\ 1\ 1\ 1\ 1\ 1)$$

(where $\rightarrow$ denotes the resultant vector has been updated and thresholded).

Thus only the very notion that their claim of Vedic Mathematics being a magic is sufficient to make all the nodes to the ON state.

Further the hidden pattern is not a limit cycle but only a fixed point. Thus the experts claims, they made 'magic' rituals for people after death and now the non-Brahmins are leading a very miserable life in their own nation. Now, what this Vedic Mathematics magic will do to the school children is to be watched very carefully because if the innocent younger generation is ruined at that adolescent stage it is sure we cannot have any hopes to rejuvenate them says the expert. Further he adds that nowadays the students' population is so streamlined that they do not participate in any social justice protests; they only mind their own business of studying, which is really a harm to the nation because we do not have well-principled, young, educated politicians to make policies for our nation.

Thus we do not know that our nation is at a loss. However the Brahmins thrive for even today they are in all the post in which they are the policy makers for the 97% of us. How can they even do any justice to us in making policies for us? They say reservation for Dalits (SC/ST) and Other Backward Classes (OBCs) should not be given in institutes of national importance because these people lack quality. This is the kind of policy they make for the non-Brahmins at large.

Now we proceed on to work with the node (4) in the ON state and all other nodes in the OFF state.

Let
$Z = (0\ 0\ 0\ 1\ 0\ 0\ 0\ 0\ 0)$

be the state vector given by the expert. Effect of Z on the system $M_1$ is given by

$ZM_1 \rightarrow (0\ 0\ 0\ 1\ 1\ 1\ 0\ 1\ 1)$
$= Z_1$ (say)

$Z_1 M_1 \rightarrow (0\ 0\ 0\ 1\ 1\ 1\ 1\ 1\ 1)$
$= Z_2;$

a fixed point. Thus the hidden pattern in this case also is a fixed point. It makes ON all the state vectors except (1) (2) and (3).

Now we proceed on to take the second expert's opinion. He is a president of a small Christian organization. The directed graph given by the 2${}^{nd}$ expert is as follows.

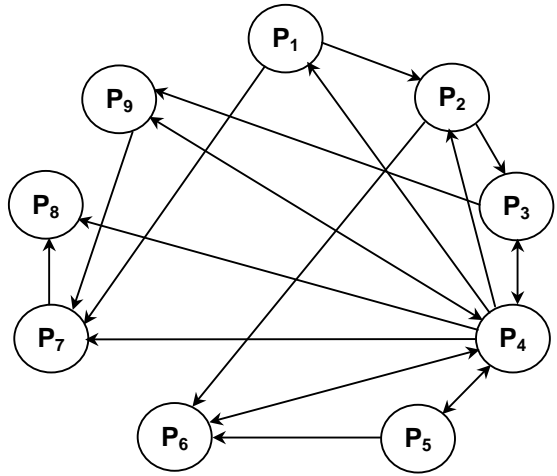

FIGURE 4.5.2

The related matrix of the directed graph given by the second expert is as follows:
We denote it by $M_2$

$$M_2 = \begin{array}{c} \\ P_1 \\ P_2 \\ P_3 \\ P_4 \\ P_5 \\ P_6 \\ P_7 \\ P_8 \\ P_9 \end{array} \begin{array}{c} P_1 \ P_2 \ P_3 \ P_4 \ P_5 \ P_6 \ P_7 \ P_8 \ P_9 \\ \begin{bmatrix} 0 & 1 & 0 & 0 & 0 & 0 & 1 & 0 & 0 \\ 0 & 0 & 1 & 0 & 0 & 1 & 0 & 0 & 0 \\ 0 & 0 & 0 & 1 & 0 & 0 & 0 & 0 & 1 \\ 1 & 1 & 1 & 0 & 1 & 1 & 1 & 1 & 1 \\ 0 & 0 & 0 & 1 & 0 & 1 & 0 & 0 & 0 \\ 0 & 0 & 0 & 1 & 0 & 0 & 0 & 0 & 0 \\ 0 & 0 & 0 & 0 & 0 & 0 & 0 & 1 & 0 \\ 0 & 0 & 0 & 0 & 0 & 0 & 0 & 0 & 0 \\ 0 & 0 & 0 & 1 & 0 & 0 & 1 & 0 & 0 \end{bmatrix} \end{array}$$

Using the dynamical system $M_2$ given by the second expert we study the same state vectors as given by the first expert, mainly for comparison purposes.
Let
X = (0 0 0 0 0 1 0 0 0)

be the state vector whose resultant we wish to study on the dynamical system $M_2$.

$$XM_2 \quad = \quad (0\ 0\ 0\ 1\ 0\ 0\ 0\ 0\ 0)$$

after updating the resultant state vector we get

$$X_1 \quad = \quad (0\ 0\ 0\ 1\ 0\ 1\ 0\ 0\ 0)$$

Now the effect of $X_1$ on the dynamical system $M_2$ is given by

$$X_1 M_2 \quad \rightarrow \quad (1\ 1\ 1\ 1\ 1\ 1\ 1\ 1\ 1)$$
$$\quad = \quad X_2.$$

Now the effect of $X_2$ on $M_2$ is

$$X_2 M_2 \quad \rightarrow \quad (1\ 1\ 1\ 1\ 1\ 1\ 1\ 1\ 1)$$
$$\quad = \quad X_3\ (=X_2).$$

Thus the resultant vector is a fixed point and all nodes come to ON state. The resultant vector given by the two experts of the dynamical systems $M_1$ and $M_2$ are distinctly different because in one case we get (0 0 0 0 0 1 1 1 1) and in case of the system $M_2$ for the same vector we get (1 1 1 1 1 1 1 1 1).

Now we study the same vector

$$Y \quad = \quad (1\ 0\ 0\ 0\ 0\ 0\ 0\ 0\ 0)$$

after updating and thresholding we get

$$YM_2 \quad = \quad Y_1$$
$$\quad = \quad (1\ 1\ 0\ 0\ 0\ 0\ 1\ 0\ 0)$$
$$Y_1 M_2 \quad \rightarrow \quad (1\ 1\ 1\ 0\ 0\ 1\ 1\ 1\ 0)$$
$$\quad = \quad Y_2\ (\text{say})$$
$$Y_2\ M_2 \quad \rightarrow \quad (1\ 1\ 1\ 1\ 0\ 1\ 1\ 1\ 1)$$
$$\quad = \quad Y_3\ (\text{say}).$$

Now

$$Y_3 M_2 \quad \rightarrow \quad (1\ 1\ 1\ 1\ 1\ 1\ 1\ 1\ 1)$$
$$\quad = \quad Y_4\ (\text{say}).$$
$$Y_4 M_2 \quad \rightarrow \quad Y_5 = (Y_4).$$

Thus we see all the nodes come to ON state. The resultant is the same as that of the first expert. Here also the hidden pattern is a fixed point that has made all the nodes to come to the ON state.

Now we take the 3$^{rd}$ state vector given by he first expert in which only the node ($P_4$) is in the ON state and all other nodes in the OFF state i.e., Z = (0 0 0 1 0 0 0 0 0).
Now we study the effect of Z on the dynamical system $M_2$,

$ZM_2$ = (1 1 1 0 1 1 1 1 1)

after updating and thresholding we get

$Z_1$ = (1 1 1 1 1 1 1 1 1);

which is a fixed point which has made all other nodes to come to the ON state. The reader can see the difference between the two resultant vectors and compare them.

Now we take the 3$^{rd}$ expert who is a Muslim activist working in minority political party; we have asked him to give his views and converted it to form the following directed graph:

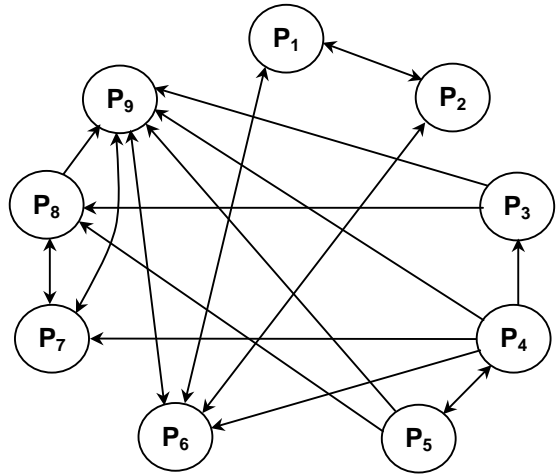

FIGURE 4.5.3

The related matrix of the directed graph given by the third expert is $M_3$

$$M_3 = \begin{array}{c} \\ P_1 \\ P_2 \\ P_3 \\ P_4 \\ P_5 \\ P_6 \\ P_7 \\ P_8 \\ P_9 \end{array} \begin{array}{c} \begin{array}{ccccccccc} P_1 & P_2 & P_3 & P_4 & P_5 & P_6 & P_7 & P_8 & P_9 \end{array} \\ \begin{bmatrix} 0 & 1 & 0 & 0 & 0 & 1 & 0 & 0 & 0 \\ 1 & 0 & 0 & 0 & 0 & 1 & 0 & 0 & 0 \\ 0 & 0 & 0 & 0 & 0 & 0 & 0 & 1 & 1 \\ 0 & 0 & 1 & 0 & 1 & 1 & 1 & 0 & 1 \\ 0 & 0 & 0 & 1 & 0 & 0 & 0 & 1 & 1 \\ 1 & 1 & 0 & 0 & 0 & 0 & 0 & 0 & 1 \\ 0 & 0 & 0 & 0 & 0 & 0 & 0 & 1 & 1 \\ 0 & 0 & 0 & 0 & 0 & 0 & 1 & 0 & 1 \\ 0 & 0 & 0 & 0 & 1 & 1 & 0 & 0 \end{bmatrix} \end{array}$$

Now we study the effect of same three state vectors given by the first expert. This is mainly done for comparison purposes.

Let $X = (0\ 0\ 0\ 0\ 0\ 1\ 0\ 0\ 0)$ be the state vector in which only the node (6) i.e., $P_6$ is in the ON state and all other nodes are in the OFF state. To study the effect of this vector on the dynamical system $M_3$.

$$XM_3 \quad = \quad (1\ 1\ 0\ 0\ 0\ 0\ 0\ 0\ 1)$$

after updating the resultant vector we get

$$X_1 \quad = \quad (1\ 1\ 0\ 0\ 0\ 1\ 0\ 0\ 1).$$

The effect of $X_1$ on the dynamical system $M_3$ is given by

$$\begin{aligned} X_1 M_3 &\rightarrow (1\ 1\ 0\ 0\ 0\ 1\ 1\ 1\ 1) \\ &= X_2 \text{ (say)} \\ X_2 M_3 &\rightarrow (1\ 1\ 0\ 0\ 0\ 1\ 1\ 1\ 1) \\ &= X_3\ (= X_2). \end{aligned}$$

Thus the hidden pattern of the resultant of the state vector X is a fixed point in which all the nodes have come to ON state. Thus resultant vector is the same as that of the second experts views and different from the first expert.

Now consider the state vector Y = (1 0 0 0 0 0 0 0 0) where all nodes are in the OFF state except the first node we wish to find the hidden pattern of Y using the dynamical system $M_3$

$YM_3$ = (0 1 0 0 0 1 0 0 0).

After updating we get

$Y_1$ = (1 1 0 0 0 1 0 0 0).

Now the effect of $Y_1$ on the dynamical system $M_3$ is given by

$Y_1 M_3$ → (1 1 0 0 0 1 0 0 1)
= $Y_2$ (say).

Effect of $Y_2$ on the dynamical system $M_3$ is given by

$Y_2 M_3$ → (1 1 0 0 0 1 1 0 1)
= $Y_3$ (say).

The resultant given by $Y_3$ is

$Y_3 M_3$ → (1 1 0 0 0 1 1 1 1)
= $Y_4$ (say).

Now the hidden pattern given by $Y_4$ using the dynamical system $M_3$ is

$Y_4 M_3$ → (1 1 0 0 0 1 1 1 1)
= $Y_5 (= Y_4)$.

Thus the hidden pattern is a fixed point. The resultant vector given by the third dynamical system $M_3$ is different from $M_1$ and $M_2$.

Now we study the effect of the state vector

Z = (0 0 0 1 0 0 0 0 0)

on the system $M_3$

$ZM_3$ = (0 0 1 0 1 1 1 0 1).

After updating we get

$Z_1$ = (0 0 1 1 1 1 1 0 1).

The effect of $Z_1$ on $M_3$ is given by

$Z_1 M_3$ → (1 1 1 1 1 1 1 1 1)
= $Z_2$ (say).

$Z_2 M_3$ → (1 1 1 1 1 1 1 1 1)
= $Z_3$ (= $Z_2$).

Thus we get a fixed point as the hidden pattern in which all the nodes come to ON state.

Now we take the views of the fourth expert, an old man who has involved himself in several political struggles and also has some views of Vedic Mathematics that some of his grandchildren studied. He heavily condemns the Hindutva policy of polluting the syllabus. We have taken his views as a public person.

Now using the directed graph

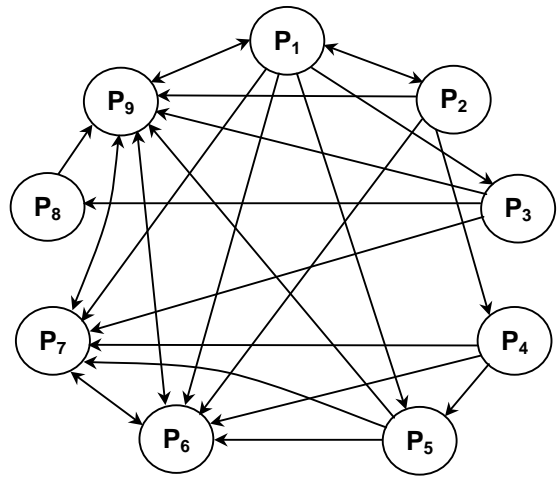

FIGURE 4.5.4

given by this expert we obtain the associated fuzzy matrix $M_4$ of the FCM.

$$M_4 = \begin{array}{c} \\ P_1 \\ P_2 \\ P_3 \\ P_4 \\ P_5 \\ P_6 \\ P_7 \\ P_8 \\ P_9 \end{array} \begin{array}{c} P_1 \; P_2 \; P_3 \; P_4 \; P_5 \; P_6 \; P_7 \; P_8 \; P_9 \\ \begin{bmatrix} 0 & 1 & 1 & 0 & 1 & 1 & 1 & 0 & 1 \\ 1 & 0 & 0 & 1 & 0 & 1 & 0 & 0 & 1 \\ 0 & 0 & 0 & 0 & 0 & 0 & 1 & 1 & 1 \\ 0 & 0 & 0 & 0 & 1 & 1 & 1 & 0 & 0 \\ 0 & 0 & 0 & 0 & 0 & 1 & 1 & 0 & 1 \\ 0 & 0 & 0 & 0 & 0 & 0 & 1 & 0 & 1 \\ 0 & 0 & 0 & 0 & 0 & 1 & 0 & 0 & 1 \\ 0 & 0 & 0 & 0 & 0 & 0 & 0 & 0 & 1 \\ 1 & 0 & 0 & 0 & 0 & 1 & 1 & 0 & 0 \end{bmatrix} \end{array}$$

Now using the matrix $M_4$ we obtain the resultant of the three state vectors viz.

1) $X = (0\;0\;0\;0\;0\;1\;0\;0\;0)$
2) $Y = (1\;0\;0\;0\;0\;0\;0\;0\;0)$
3) $Z = (0\;0\;0\;1\;0\;0\;0\;0\;0)$.

Consider the state vector $X = (0\;0\;0\;0\;0\;1\;0\;0\;0)$ given by the first expert in which only the node (6) is in the ON state and all other nodes are in the off state. The effect of $X$ on the dynamical system $M_4$ is given by

$XM_4 \quad = \quad (0\;0\;0\;0\;0\;0\;1\;0\;1)$.

after updating we get

$X_1 \quad = \quad (0\;0\;0\;0\;0\;1\;1\;0\;1)$.

The effect of $X_1$ on $M_4$ is given by
$X_1 M_4 \quad \rightarrow \quad (1\;0\;0\;0\;0\;1\;1\;0\;1)$
$\quad = \quad X_2$ (say).

Now $X_2$ acts on the dynamical system $M_4$ and gives
$X_2 M_4 \quad \rightarrow \quad (1\;1\;1\;1\;1\;1\;1\;0\;1)$
$\quad = \quad X_3$ (say).

Now the effect of $X_3$ is given by

$$X_3 M_4 \rightarrow (1\ 1\ 1\ 1\ 1\ 1\ 1\ 1\ 1)$$
$$= X_4 \text{ (say)}.$$

Now when $X_4$ is passed through $M_4$ we get

$$X_4 M_4 \rightarrow (1\ 1\ 1\ 1\ 1\ 1\ 1\ 1\ 1)$$
$$= X_5 \,(= X_4).$$

Thus the hidden pattern of the state vector X is given by (1 1 1 1 1 1 1 1 1), which is a fixed point. All nodes come to ON state. This resultant is different from the other experts' opinions.

Now we proceed on to study the effect of Y on the dynamical system $M_4$, where
$$Y = (1\ 0\ 0\ 0\ 0\ 0\ 0\ 0\ 0)$$
all nodes except node (1) is in the ON state.

$$Y M_4 = (0\ 1\ 1\ 0\ 1\ 1\ 1\ 0\ 1)$$

after updating we get
$$Y_1 = (1\ 1\ 1\ 0\ 1\ 1\ 1\ 0\ 1).$$

Now we study the effect of $Y_1$ on $M_4$
$$Y_1 M_4 \rightarrow (1\ 1\ 1\ 1\ 1\ 1\ 1\ 1\ 1)$$
$$= Y_2 \text{ (say)}.$$
$$Y_2 M_4 \rightarrow (1\ 1\ 1\ 1\ 1\ 1\ 1\ 1\ 1)$$
$$= Y_3 \,(= Y_2).$$

Thus the hidden pattern of Y is a fixed point. This resultant is also different from that of the others.

Now we proceed on to study the effect of the state vector
$$Z = (0\ 0\ 0\ 1\ 0\ 0\ 0\ 0\ 0);$$

where all nodes are in the off state except the node (4).
Now
$$Z M_4 = (0\ 0\ 0\ 0\ 1\ 1\ 1\ 0\ 0)$$

After updating we get

$Z_1 = (0\ 0\ 0\ 1\ 1\ 1\ 1\ 0\ 0)$

$Z_1 M_4 \rightarrow (0\ 0\ 0\ 1\ 1\ 1\ 1\ 0\ 1)$
$= Z_2$ (say)

$Z_2 M_4 \rightarrow (1\ 0\ 0\ 1\ 1\ 1\ 1\ 0\ 1)$
$= Z_3$ (say)

$Z_3 M_4 \rightarrow (1\ 1\ 1\ 1\ 1\ 1\ 1\ 0\ 1)$
$= Z_4$ (say)

$Z_4 M_4 \rightarrow (1\ 1\ 1\ 1\ 1\ 1\ 1\ 1\ 1)$
$= Z_5$ (say)

$Z_5 M_4 \rightarrow (1\ 1\ 1\ 1\ 1\ 1\ 1\ 1\ 1)$
$= Z_6 (= Z_5)$.

Thus the hidden pattern of Z using the dynamical system $M_4$ is the fixed point given by $(1\ 1\ 1\ 1\ 1\ 1\ 1\ 1\ 1)$. The reader can study the differences and similarities from the other four experts.

Now we have taken the 5$^{th}$ expert who is a feminist and currently serves as the secretary of a women association and who showed interest and enthusiasm in this matter. The directed graph given by this expert is as follows:

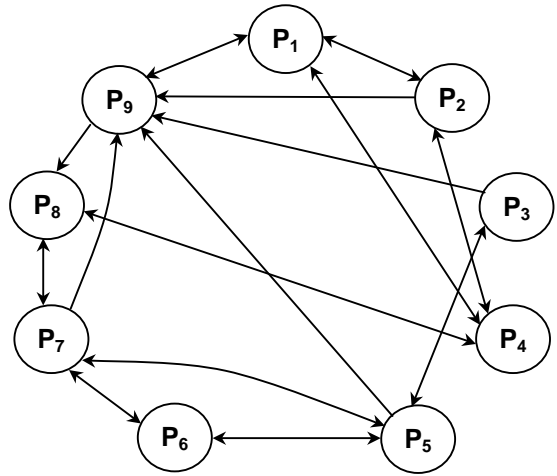

FIGURE 4.5.5

The connection matrix related to the directed is given by the matrix $M_5$

$$M_5 = \begin{array}{c} \\ P_1 \\ P_2 \\ P_3 \\ P_4 \\ P_5 \\ P_6 \\ P_7 \\ P_8 \\ P_9 \end{array} \begin{array}{c} P_1\ P_2\ P_3\ P_4\ P_5\ P_6\ P_7\ P_8\ P_9 \\ \begin{bmatrix} 0 & 1 & 0 & 1 & 0 & 0 & 0 & 0 & 1 \\ 1 & 0 & 0 & 1 & 0 & 0 & 0 & 0 & 1 \\ 0 & 0 & 0 & 0 & 1 & 0 & 0 & 0 & 1 \\ 1 & 1 & 0 & 0 & 0 & 0 & 0 & 1 & 0 \\ 0 & 0 & 1 & 0 & 0 & 1 & 1 & 0 & 1 \\ 0 & 0 & 0 & 0 & 1 & 0 & 1 & 0 & 0 \\ 0 & 0 & 0 & 0 & 1 & 1 & 0 & 1 & 1 \\ 0 & 0 & 0 & 1 & 0 & 0 & 1 & 0 & 0 \\ 1 & 0 & 0 & 0 & 0 & 0 & 0 & 1 & 0 \end{bmatrix} \end{array}$$

Now consider the state vector $X = (0\ 0\ 0\ 0\ 0\ 1\ 0\ 0\ 0)$ as given by the first expert, where only the node (6) is in the ON state and all other nodes are in the OFF state. The effect of $X$ on the dynamical system $M_5$ is given by

$XM_5 \quad = \quad (0\ 0\ 0\ 0\ 1\ 0\ 1\ 0\ 0)$

after updating the resultant state vector we get

$X_1 \quad = \quad (0\ 0\ 0\ 0\ 1\ 1\ 1\ 0\ 0)$.

The effect of $X_1$ on $M_5$ is given by

$X_1 M_5 \quad \to \quad (0\ 0\ 1\ 0\ 1\ 1\ 1\ 1\ 1)$
$\quad\quad\quad = \quad X_2$ (say)
$X_2 M_5 \quad \to \quad (1\ 0\ 1\ 1\ 1\ 1\ 1\ 1\ 1)$
$\quad\quad\quad = \quad X_3$ (say)
$X_3 M_5 \quad \to \quad (1\ 1\ 1\ 1\ 1\ 1\ 1\ 1\ 1)$
$\quad\quad\quad = \quad X_4$ (say)
$X_4 M_5 \quad \to \quad (1\ 1\ 1\ 1\ 1\ 1\ 1\ 1\ 1)$
$\quad\quad\quad = \quad X_5\ (= X_4)$.

The hidden pattern happens to be a fixed point in which all the nodes have come to ON state. Next we study the effect of the state vector

$$Y = (1\ 0\ 0\ 0\ 0\ 0\ 0\ 0\ 0)$$

on the dynamical system $M_5$.

$$YM_5 = (0\ 1\ 0\ 1\ 0\ 0\ 0\ 0\ 1)$$

After updating we get the resultant as

$$Y_1 = (1\ 1\ 0\ 1\ 0\ 0\ 0\ 0\ 1).$$

The effect of $Y_1$ on $M_5$ is given by

$$Y_1 M_5 \rightarrow (1\ 1\ 0\ 1\ 0\ 0\ 0\ 1\ 1)$$
$$= Y_2 \text{ (say)}$$
$$Y_2 M_5 \rightarrow (1\ 1\ 0\ 1\ 1\ 0\ 1\ 1\ 1)$$
$$= Y_3 \text{ (say)}$$
$$Y_3 M_5 \rightarrow (1\ 1\ 1\ 1\ 1\ 1\ 1\ 1\ 1)$$
$$= Y_4 \text{ (say)}$$
$$Y_4 M_5 \rightarrow (1\ 1\ 1\ 1\ 1\ 1\ 1\ 1\ 1)$$
$$= Y_5 (= Y_4).$$

Thus the hidden pattern is a fixed point. We see that when the node (1) alone is in the ON state all other nodes come to ON state there by showing when Vedic Mathematics is based on magic it has several ulterior motives and no one with any common sense will accept it as mathematics according this expert.

Now we study the effect of the node $Z = (0\ 0\ 0\ 1\ 0\ 0\ 0\ 0\ 0)$ where only the node (4) is in the ON state and all other nodes are in the OFF state. The effect of Z on the dynamical system $M_5$ is given by

$$ZM_5 = (1\ 1\ 0\ 0\ 0\ 0\ 0\ 1\ 0)$$

after updating we obtain the following resultant vector;

$X_1$ = (1 1 0 1 0 0 0 1 0).

The effect of $X_1$ on $M_5$ is given by

$X_1 M_5$ → (1 1 0 1 1 0 0 1 1 1)
= $X_2$ (say)
$X_2 M_5$ → (1 1 1 1 1 1 1 1 1)
= $X_3$ (say)
$X_3 M_5$ → (1 1 1 1 1 1 1 1 1)
= $X_4$ (=$X_3$).

Thus the hidden pattern of the vector Z is a fixed point. When the nodes Christians and Muslims do not accept Vedic Mathematics shows all the nodes came to ON state it is a Hindutva agenda it is not mathematics to really improve the students, it has all ulterior motives to saffronize the nation and there by establish the supremacy of the Aryans.

Now we seek the views of the sixth expert who is a political worker.

The directed graph given by the 6$^{th}$ expert is as follows:

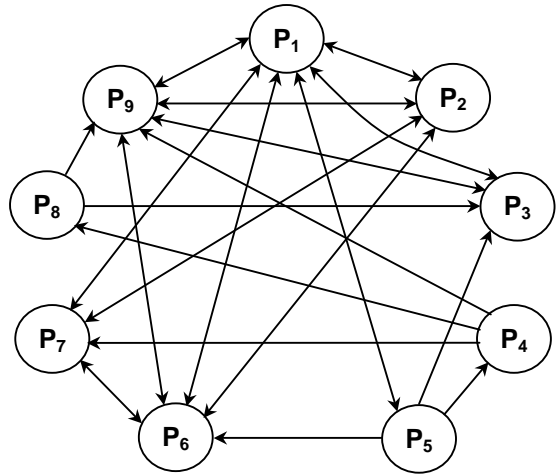

FIGURE 4.5.6

Using the directed graph given by the expert we obtain the following fuzzy matrix $M_6$.

$$M_6 = \begin{array}{c} \\ P_1 \\ P_2 \\ P_3 \\ P_4 \\ P_5 \\ P_6 \\ P_7 \\ P_8 \\ P_9 \end{array} \begin{array}{c} P_1 \; P_2 \; P_3 \; P_4 \; P_5 \; P_6 \; P_7 \; P_8 \; P_9 \\ \begin{bmatrix} 0 & 1 & 1 & 0 & 1 & 1 & 1 & 0 & 1 \\ 1 & 0 & 0 & 0 & 0 & 1 & 1 & 0 & 1 \\ 1 & 0 & 0 & 0 & 0 & 0 & 0 & 0 & 1 \\ 0 & 0 & 0 & 0 & 0 & 0 & 1 & 1 & 1 \\ 1 & 0 & 1 & 1 & 0 & 1 & 0 & 0 & 0 \\ 1 & 1 & 0 & 0 & 0 & 0 & 1 & 0 & 1 \\ 1 & 1 & 0 & 0 & 0 & 1 & 0 & 0 & 0 \\ 0 & 0 & 1 & 0 & 0 & 0 & 0 & 0 & 1 \\ 1 & 1 & 1 & 0 & 0 & 1 & 0 & 0 & 0 \end{bmatrix} \end{array}$$

Using this dynamical system we obtain the resultant of the three vectors.

1) $X$ = (0 0 0 0 0 1 0 0 0)
2) $Y$ = (1 0 0 0 0 0 0 0 0)
and 3) $Z$ = (0 0 0 1 0 0 0 0 0).

The effect of

$X$ = (0 0 0 0 0 1 0 0 0)

on the dynamical system $M_6$ is given by

$XM_6$ = (1 1 0 0 0 0 1 0 1)

after updating we get

$X_1$ = (1 1 0 0 0 1 1 0 1).

Now the effect of $X_1$ on $M_6$ is given by

$X_1M_6$ → (1 1 1 0 1 1 1 0 1)
     = $X_2$ (say)

$$X_2 M_6 \rightarrow (1\ 1\ 1\ 1\ 1\ 1\ 1\ 0\ 1)$$
$$= X_3 \text{ (say)}$$
$$X_3 M_6 \rightarrow (1\ 1\ 1\ 1\ 1\ 1\ 1\ 1\ 1)$$
$$= X_4 = (X_3).$$

Thus when only the node (6) is in the ON state we get the hidden pattern of the resultant vector to be a fixed point which makes all the other nodes come to the ON state.

Now we study the effect of $Y = (1\ 0\ 0\ 0\ 0\ 0\ 0\ 0\ 0)$ i.e only the node (1) is in the ON state and all other nodes are in the OFF state; effect of Y on the dynamical system $M_6$ is given by

$$Y M_6 = (0\ 1\ 1\ 0\ 1\ 1\ 1\ 0\ 1).$$

After updating we get the resultant

$$Y_1 = (1\ 1\ 1\ 0\ 1\ 1\ 1\ 0\ 1).$$

Now the resultant of $Y_1$ on the dynamical system $M_6$ is given by

$$Y_1 M_6 \rightarrow (1\ 1\ 1\ 1\ 1\ 1\ 1\ 0\ 1)$$
$$= Y_2 \text{ (say)}$$
$$Y_2 M_6 \rightarrow (1\ 1\ 1\ 1\ 1\ 1\ 1\ 1\ 1)$$
$$= Y_3 (= Y_2).$$

Thus the hidden pattern is a fixed point we see that when the concept 'Vedic Mathematics is a magic according to their claims' is alone in the ON state, all the other nodes come to the ON state by which it is evident that Vedic Mathematics has more ulterior motives and it is not Mathematics because mathematics cannot be magic. Mathematics is a science of down to earth reality.

Now we study the effect of the vector

$$Z = (0\ 0\ 0\ 1\ 0\ 0\ 0\ 0\ 0)$$

where only the node (4) is in the ON state and all other nodes are in the OFF state.

$ZM_6 = (0\ 0\ 0\ 0\ 0\ 0\ 1\ 1\ 1)$

After updating we got the resultant vector to be

$Z_1 = (0\ 0\ 0\ 1\ 0\ 0\ 1\ 1\ 1)$

$Z_1 M_6 \to (1\ 1\ 1\ 1\ 0\ 1\ 1\ 1\ 1)$
$= Z_2$ (say)

$Z_2 M_6 \to (1\ 1\ 1\ 1\ 1\ 1\ 1\ 1\ 1)$
$= Z_3$ (say)

$Z_3 M_6 \to (1\ 1\ 1\ 1\ 1\ 1\ 1\ 1\ 1)$
$= Z_4 (= Z_3)$.

Thus the hidden pattern of this vector Z is a fixed point that makes all the nodes into ON state, i.e., when the Christians and Muslims of India do not accept Vedic Mathematics it means that it has ulterior motives and above all shows that it is a Hindutva agenda.

Thus, now we have seen the same set of vectors by all three experts. It is left for the reader to make comparisons. Now we give the opinion of the 7$^{th}$ expert who is a human rights activists working in an NGO in the form of the directed graph.

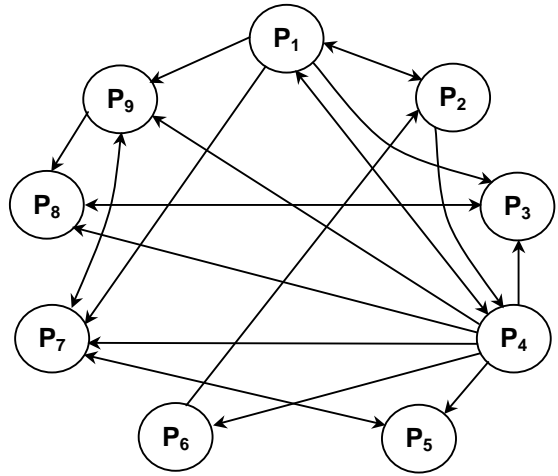

FIGURE 4.5.7

Now we obtain the connection matrix $M_7$ using the directed graph.

$$M_7 = \begin{array}{c} \\ P_1 \\ P_2 \\ P_3 \\ P_4 \\ P_5 \\ P_6 \\ P_7 \\ P_8 \\ P_9 \end{array} \begin{array}{c} P_1\ P_2\ P_3\ P_4\ P_5\ P_6\ P_7\ P_8\ P_9 \\ \left[ \begin{array}{ccccccccc} 0 & 1 & 1 & 1 & 0 & 0 & 1 & 0 & 1 \\ 1 & 0 & 0 & 1 & 0 & 0 & 0 & 0 & 0 \\ 0 & 0 & 0 & 0 & 0 & 0 & 0 & 1 & 0 \\ 1 & 1 & 1 & 0 & 1 & 1 & 1 & 1 & 1 \\ 0 & 0 & 0 & 0 & 0 & 0 & 1 & 0 & 0 \\ 0 & 1 & 0 & 0 & 0 & 0 & 0 & 0 & 0 \\ 0 & 0 & 0 & 0 & 1 & 0 & 0 & 0 & 1 \\ 0 & 0 & 1 & 0 & 0 & 0 & 0 & 0 & 0 \\ 0 & 0 & 0 & 0 & 0 & 0 & 1 & 1 & 0 \end{array} \right] \end{array}$$

This expert wanted to work with some other set of three vectors so we start to work with state vectors as suggested by him. He wants the node (9) alone to be in the ON state and all other nodes to be in the OFF state. Let

X    =    (0 0 0 0 0 0 0 0 1).

Now we study the effect of X on the dynamical system $M_7$,

XM$_7$    =    (0 0 0 0 0 0 1 1 0)

after updating we get,

X$_1$    =    (0 0 0 0 0 0 1 1 1).

The effect of X$_1$ on M is given by

X$_1$ M$_7$    →    (0 0 1 0 1 0 1 1 1)
         =    X$_2$ (say)

X$_2$ M$_7$    →    (0 0 1 0 1 0 1 1 1)
         =    X$_3$ (= X$_2$).

Thus the hidden pattern of the dynamical system is a fixed point. Now we proceed on to work with the state vector (0 0 0 0 0 0 1 0 0) where only the node (7) is in the ON state and all other nodes are in the OFF state.

The effect of Y on the dynamical system $M_7$ is given by

$\quad$ $YM_7$ $\quad = \quad$ (0 0 0 0 1 0 0 0 1)

after updating we get

$\quad$ $Y_1$ $\quad = \quad$ (0 0 0 0 1 0 1 0 1)

Now the effect of $Y_1$ on $M_7$ is given by

$\quad$ $Y_1 M_7$ $\quad \to \quad$ (0 0 0 0 1 0 1 1 1)
$\quad\quad\quad\quad = \quad$ $Y_2$ (say)
$\quad$ $Y_2 M_7$ $\quad \to \quad$ (0 0 0 0 1 0 1 1 1)
$\quad\quad\quad\quad = \quad$ $Y_3$ (= $Y_2$).

Thus the hidden pattern of the dynamical system is a fixed point.

Now we study the state vector

$\quad$ Z $\quad = \quad$ (0 1 0 0 0 0 0 0 0)

here the node (2) i.e., Vedic Mathematics is 'trick' alone is in the ON state and all other nodes are in the OFF state.

The effect of Z on the dynamical system $M_7$ is given by

$\quad$ $ZM_7$ $\quad = \quad$ (1 0 0 1 0 0 0 0 0).

After updating we get

$\quad$ $Z_1$ $\quad = \quad$ (1 1 0 1 0 0 0 0 0).

Now the effect of $Z_1$ on dynamical system $M_7$ is given by

$Z_1 M_7$ → (1 1 1 1 1 1 1 1 1)
= $Z_2$ (Say)

$Z_2 M_7$ → (1 1 1 1 1 1 1 1 1)
= $Z_3$ (= $Z_2$).

Thus the hidden pattern is a fixed point and all the nodes come to ON state. Thus according to this expert Vedic Mathematics uses 'trick' to solve arithmetical problems is enough to condemn Vedic Mathematics as a tool which has ulterior motives to make the nation come under the influence of revivalist and fundamentalist Hindutva.

Next we take the opinion of an expert who is a union leader, who has studied up to the 10[th] standard and belongs to a socially and economically backward community.

The opinion of the 8[th] expert is given by the following directed graph:

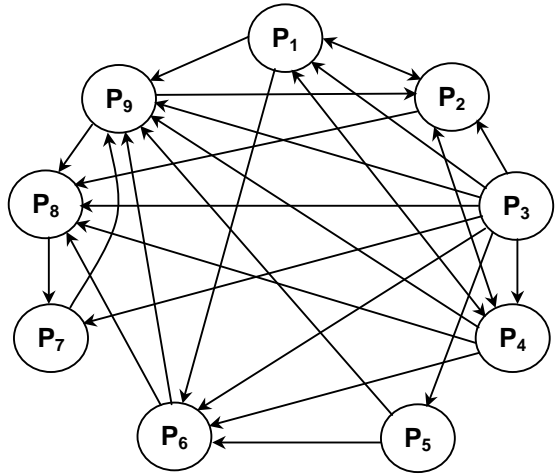

FIGURE 4.5.8

The related relational matrix $M_8$ is given in the following:

$$M_8 = \begin{array}{c} \\ P_1 \\ P_2 \\ P_3 \\ P_4 \\ P_5 \\ P_6 \\ P_7 \\ P_8 \\ P_9 \end{array} \begin{array}{c} \begin{array}{ccccccccc} P_1 & P_2 & P_3 & P_4 & P_5 & P_6 & P_7 & P_8 & P_9 \end{array} \\ \left[ \begin{array}{ccccccccc} 0 & 1 & 0 & 1 & 0 & 1 & 0 & 0 & 1 \\ 1 & 0 & 0 & 1 & 0 & 0 & 0 & 1 & 0 \\ 1 & 1 & 0 & 1 & 1 & 1 & 1 & 1 & 1 \\ 1 & 1 & 0 & 0 & 0 & 1 & 0 & 1 & 1 \\ 0 & 0 & 0 & 0 & 0 & 1 & 0 & 0 & 1 \\ 0 & 0 & 0 & 0 & 0 & 0 & 0 & 1 & 1 \\ 0 & 0 & 0 & 0 & 0 & 0 & 0 & 0 & 1 \\ 0 & 0 & 0 & 0 & 0 & 0 & 1 & 0 & 0 \\ 0 & 1 & 0 & 0 & 0 & 0 & 0 & 1 & 0 \end{array} \right] \end{array}$$

Now we study the effect of the same state vectors as given by the 8$^{th}$ expert.
Given

$X = (0\ 0\ 0\ 0\ 0\ 0\ 0\ 0\ 1)$.

Now

$XM_8 = (0\ 1\ 0\ 0\ 0\ 0\ 0\ 1\ 0)$.

After updating we get

$X_1 = (0\ 1\ 0\ 0\ 0\ 0\ 0\ 1\ 1)$.

The effect of $X_2$ on the dynamical system $M_8$ is given by

$X_2 M_8 \rightarrow (1\ 1\ 0\ 1\ 0\ 0\ 1\ 1\ 1)$
$\quad\quad\ = X_3$ (say).
$X_3 M_8 \rightarrow (1\ 1\ 0\ 1\ 1\ 1\ 1\ 1\ 1)$
$\quad\quad\ = X_4$ (say).
$X_4 M_8 = (1\ 1\ 0\ 0\ 1\ 1\ 1\ 1\ 1)$
$\quad\quad\ = X_5\ (= X_4)$.

Thus the hidden pattern is a fixed point. Except for the nodes (3) and (4) all other nodes come to the ON state. Now we study the effect of the state vector $Y = (0\ 0\ 0\ 0\ 0\ 0\ 1\ 0\ 0)$ given by the $7^{th}$ expert. On the dynamical system $M_8$ given by the 8th expert

$$YM_8 = (0\ 0\ 0\ 0\ 0\ 0\ 0\ 0\ 1).$$

After updating we get the resultant to be

$$Y_1 = (0\ 0\ 0\ 0\ 0\ 0\ 1\ 0\ 1).$$

Now

$$\begin{aligned}
Y_1 M_8 &\to (0\ 1\ 0\ 1\ 0\ 0\ 1\ 1\ 1) \\
&= Y_2 \text{ (say)} \\
Y_2 M_8 &\to (1\ 1\ 0\ 1\ 0\ 1\ 1\ 1\ 1) \\
&= Y_3 \text{ (say)} \\
Y_3 M_8 &\to (1\ 1\ 0\ 1\ 1\ 1\ 1\ 1\ 1) \\
&= Y_4 \text{ (say)} \\
Y_4 M_8 &\to (1\ 1\ 0\ 1\ 1\ 1\ 1\ 1\ 1) \\
&= Y_5\ (= Y_4).
\end{aligned}$$

Thus the hidden pattern of the state vector Y given by the dynamical system $M_8$ is a fixed point.

Now we study the effect of the state vector $Z = (0\ 1\ 0\ 0\ 0\ 0\ 0\ 0\ 0)$ on $M_8$

$$ZM_8 = (1\ 0\ 0\ 1\ 0\ 0\ 0\ 1\ 0).$$

After updating we get

$$Z_1 = (1\ 1\ 0\ 1\ 0\ 0\ 0\ 1\ 0)$$

The effect of $Z_1$ on $M_8$ is given by

$$\begin{aligned}
Z_1 M_8 &\to (1\ 1\ 0\ 1\ 0\ 1\ 1\ 1\ 1) \\
&= Z_2. \\
Z_2 M_8 &\to (1\ 1\ 0\ 1\ 0\ 1\ 1\ 1\ 1) \\
&= Z_3\ (= Z_2).
\end{aligned}$$

Thus the resultant vector is a fixed point. According to this expert the notion of Christians and Muslims not following Hindutva and Vedic Mathematics is only due to Sanskrit phrases and words.

Now we proceed on to work with the 9th expert who is a freelance writer in Tamil. He has failed in his 10th standard examination, that too in mathematics. He is now in his late fifties. He writes about social issues, poems and short stories in Tamil. Having failed in mathematics, he has spent his whole life being scared of mathematics. He says he was asked by a weekly magazine to write about Vedic Mathematics and they gave him two Tamil books in Vedic Mathematics so that he could make use of them for writing his article. He studied both the books and says that most of the arithmetical problems are very simple and elementary, like the primary school level. He says that he wrote an essay in which he strongly criticized the Swamiji for writing such stuff and calling it Vedic Mathematics. He said there was nothing Vedic in that book and even with his standard he could find any mathematics in it. So he very strongly opposed it and viewed it in the angle of an attempt to saffronize the nation. When the editor of the journal took the article he was upset about the way it was written and said they could not publish it and suggested many changes. However this writer refused to do a positive review.

Now we catch his opinion as a directed graph.

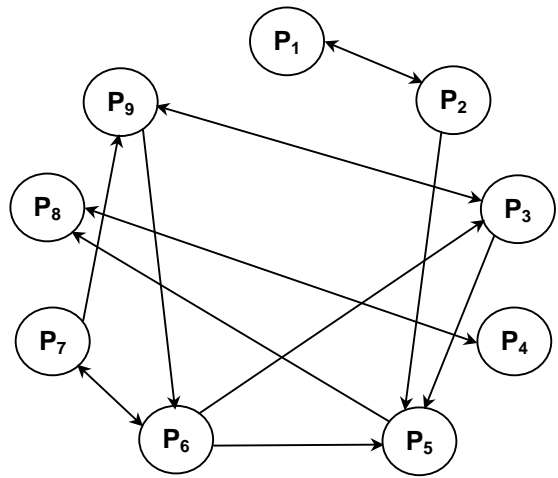

FIGURE 4.5.9

He is taken as the 9$^{th}$ expert to give views about Vedic Mathematics.

Using the directed graph we have the following connection matrix $M_9$ given in the following:

$$M_9 = \begin{array}{c} \\ P_1 \\ P_2 \\ P_3 \\ P_4 \\ P_5 \\ P_6 \\ P_7 \\ P_8 \\ P_9 \end{array} \begin{array}{c} P_1\ P_2\ P_3\ P_4\ P_5\ P_6\ P_7\ P_8\ P_9 \\ \begin{bmatrix} 0 & 1 & 0 & 0 & 0 & 0 & 0 & 0 & 0 \\ 1 & 0 & 0 & 0 & 1 & 0 & 0 & 0 & 0 \\ 0 & 0 & 0 & 0 & 1 & 0 & 0 & 0 & 1 \\ 0 & 0 & 0 & 0 & 0 & 0 & 0 & 1 & 0 \\ 0 & 0 & 0 & 0 & 0 & 0 & 0 & 1 & 0 \\ 0 & 0 & 1 & 0 & 1 & 0 & 1 & 0 & 0 \\ 0 & 0 & 0 & 0 & 0 & 1 & 0 & 0 & 1 \\ 0 & 0 & 0 & 1 & 0 & 0 & 0 & 0 & 0 \\ 0 & 0 & 1 & 0 & 0 & 1 & 0 & 0 & 0 \end{bmatrix} \end{array}$$

Now we study the resultant of the three state vectors given by the 7$^{th}$ expert

$$X = (0\ 0\ 0\ 0\ 0\ 0\ 0\ 0\ 1)$$
$$Y = (0\ 0\ 0\ 0\ 0\ 0\ 1\ 0\ 0)$$

and

$$Z = (0\ 1\ 0\ 0\ 0\ 0\ 0\ 0\ 0).$$

The effect of X on the dynamical system $M_9$ is given by

$XM_9 \quad = \quad (0\ 0\ 1\ 0\ 0\ 1\ 0\ 0\ 0).$

After updating we get

$X_1 \quad = \quad (0\ 0\ 1\ 0\ 0\ 1\ 0\ 0\ 1).$

Now the effect of $X_1$ on the system $M_9$ is given by

$X_1 M_9 \quad \rightarrow \quad (0\ 0\ 1\ 0\ 1\ 1\ 1\ 0\ 1)$
$\quad\quad\quad\quad = \quad X_2$ (say).

The effect of $X_2$ on the dynamical system $M_9$ is given by

$$
\begin{aligned}
X_2 M_9 &\to (0\ 0\ 1\ 0\ 1\ 1\ 1\ 0\ 1) \\
&= X_3 \text{ (say)} \\
X_3 M_9 &\to (0\ 0\ 1\ 0\ 1\ 1\ 1\ 0\ 1) \\
&= X_4 (= X_3).
\end{aligned}
$$

Thus the hidden pattern of the dynamical system is a fixed point.

Now we study the effect of

$$Y = (0\ 0\ 0\ 0\ 0\ 0\ 1\ 0\ 0)$$

on the system $M_9$ where only the node (7) is in the ON state i.e., Vedic Mathematics imposes Aryan supremacy on the non-Brahmins and all other nodes are in the OFF state.

The effect of Y on the system $M_9$ is given by

$$YM_9 = (0\ 0\ 0\ 0\ 0\ 1\ 0\ 0\ 1)$$

after updating we get.

$$Y_1 = (0\ 0\ 0\ 0\ 0\ 1\ 1\ 0\ 1).$$

Now the resultant vector when $Y_1$ is passed into the dynamical system $M_9$ is given by

$$
\begin{aligned}
Y_1 M_9 &\to (0\ 0\ 1\ 0\ 1\ 1\ 1\ 0\ 1) \\
&= Y_2 \text{ (say)}. \\
Y_2 M_9 &\to (0\ 0\ 1\ 0\ 1\ 1\ 1\ 1\ 1) \\
&= Y_3 \\
Y_3 M_9 &\to (0\ 0\ 1\ 1\ 1\ 1\ 1\ 1\ 1) \\
&= Y_4.
\end{aligned}
$$

Thus the resultant is a fixed point.

Now we proceed on to find the effect of the state vector.

$$Z = (0\ 1\ 0\ 0\ 0\ 0\ 0\ 0)\text{ on }M_9$$

$$ZM_9 = (1\ 0\ 0\ 0\ 1\ 0\ 0\ 0).$$

After updating we get

$$Z_1 = (1\ 1\ 0\ 0\ 1\ 0\ 0\ 0)$$

$$Z_1M_9 \to (1\ 1\ 0\ 0\ 1\ 0\ 0\ 1\ 0)$$
$$= Z_2\text{ (say)}$$
$$Z_2 M_9 \to (1\ 1\ 0\ 1\ 1\ 0\ 0\ 1\ 0)$$
$$= Z_3\text{ (say)}$$
$$Z_3 M_9 \to (1\ 1\ 0\ 1\ 1\ 0\ 0\ 1\ 0)$$
$$= Z_4\ (= Z_3).$$

Thus the hidden pattern is a fixed point.

Now we proceed on to take the 10$^{th}$ expert who is a social worker. She failed in her 12$^{th}$ standard but does social work without any anticipation for public recognition or honour. She is in her late forties. As she was also taking adult education classes besides helping children in their studies we have taken her views. She was aware of Vedic Mathematics and said that she used it to find shortcut methods but it was not of much use to her. The reason for its non-usefulness according to her is because for every individual type of problem we have to remember a method or some of its properties that did not apply uniformly. So she did not like it. She also came down heavily upon the cover pages of the Vedic Mathematics books (1) and (2) in Tamil [85-6]. She says that though she is a religious Hindu yet as a social worker she does not want to discriminate anyone based on religion.

Also she said that she has faced several problems with the Brahmin priest of the temple and his family members. Though they were only one family yet they were always opposed to her because they did not like the villagers in their village to be reformed or educated and live with a motive and goal. They had started giving her several problems when she began to educate people of good things. Now she is educating the people not to visit temples and put money for him. Now we give the directed graph given by this woman who is our $10^{th}$ expert.

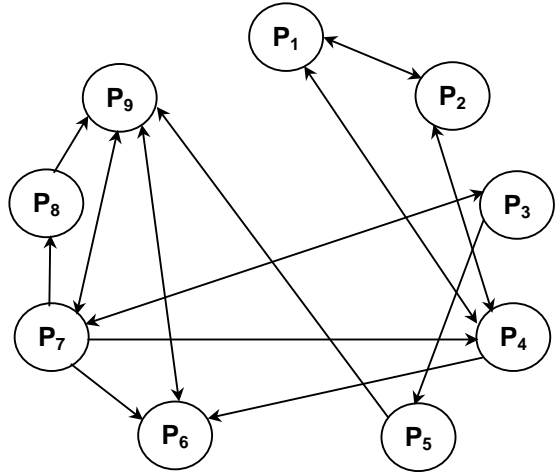

FIGURE 4.5.10

Now using this directed graph we have the following connection matrix $M_{10}$:

$$M_{10} = \begin{array}{c} \\ P_1 \\ P_2 \\ P_3 \\ P_4 \\ P_5 \\ P_6 \\ P_7 \\ P_8 \\ P_9 \end{array} \begin{array}{c} P_1 \ P_2 \ P_3 \ P_4 \ P_5 \ P_6 \ P_7 \ P_8 \ P_9 \\ \begin{bmatrix} 0 & 1 & 0 & 1 & 0 & 0 & 0 & 0 & 0 \\ 1 & 0 & 0 & 1 & 0 & 0 & 0 & 0 & 0 \\ 0 & 0 & 0 & 0 & 1 & 0 & 1 & 0 & 0 \\ 1 & 1 & 0 & 0 & 0 & 1 & 0 & 0 & 0 \\ 0 & 0 & 0 & 0 & 0 & 0 & 0 & 0 & 1 \\ 0 & 0 & 0 & 0 & 0 & 0 & 0 & 0 & 1 \\ 0 & 0 & 1 & 1 & 0 & 1 & 0 & 1 & 1 \\ 0 & 0 & 0 & 0 & 0 & 0 & 0 & 0 & 1 \\ 0 & 0 & 0 & 0 & 0 & 1 & 1 & 0 & 0 \end{bmatrix} \end{array}$$

Now using this dynamical system $M_{10}$ we study the effect of the vectors X, Y and Z given by the $7^{th}$ expert.

Let

$X \quad = \quad (0\ 0\ 0\ 0\ 0\ 0\ 0\ 0\ 1)$

be the state vector which has only node 9 in the ON state and all other nodes are in the OFF state.

The effect of X on the system $M_{10}$ is given by;

$X \quad = \quad (0\ 0\ 0\ 0\ 0\ 0\ 0\ 0\ 1)$

$XM_{10} \quad = \quad (0\ 0\ 0\ 0\ 0\ 1\ 1\ 0\ 0).$

After updating we get

$X_1 \quad = \quad (0\ 0\ 0\ 0\ 0\ 1\ 1\ 0\ 1)$
$X_1 M_{10} \quad \rightarrow \quad (0\ 0\ 1\ 1\ 0\ 1\ 1\ 1\ 1)$
$\quad \quad = \quad X_2.$
$X_2 M_{10} \quad \rightarrow \quad (1\ 1\ 1\ 1\ 1\ 1\ 1\ 1\ 1)$
$\quad \quad = \quad X_3.$
$X_3 M_{10} \quad \rightarrow \quad (1\ 1\ 1\ 1\ 1\ 1\ 1\ 1\ 1)$
$\quad \quad \quad \quad X_4\ (= X_3).$

Thus the hidden pattern is a fixed point and the node (9) alone that Vedic Mathematics has the political agenda to rule the nation is sufficient to make all the other nodes to come to the ON state.

Now we consider the state vector $Y = (0\ 0\ 0\ 0\ 0\ 0\ 1\ 0\ 0)$; i.e only the node (7) alone is in the ON state and all other nodes are in the OFF state. The effect of Y on the dynamical system $M_{10}$ is given by

$YM_{10} \quad = \quad (0\ 0\ 1\ 1\ 0\ 1\ 0\ 1\ 1).$

After updating we get

$Y_1 \quad = \quad (0\ 0\ 1\ 1\ 0\ 1\ 1\ 1\ 1).$

Now the effect of $Y_1$ on the dynamical system $M_{10}$ is given by

$Y_1 M_{10}$ → (1 1 1 1 1 1 1 1 1)
= $Y_2$ (say).

$Y_2 M_{10}$ → (1 1 1 1 1 1 1 1 1)
= $Y_3 (= Y_2)$.

Thus the resultant is a fixed point and all nodes come to ON state, when the agenda of Vedic Mathematics is to establish the superiority of Aryans.

Now we proceed on to find the effect of the state vector

$Z$ = (0 1 0 0 0 0 0 0)

where only the node (2) is in the ON state and all other nodes are in the OFF state.

The effect of Z on $M_{10}$ is given by

$Z M_{10}$ = (1 0 0 1 0 0 0 0 0).

After updating the resultant vector we get

$Z_1$ = (1 1 0 1 0 0 0 0 0).

The effect of $Z_1$ on the system $M_{10}$ is given by

$Z_1 M_{10}$ → (1 1 0 1 0 1 0 0 0)
= $Z_2$ (say).
$Z_2 M_{10}$ → (1 1 0 1 0 1 0 0 0)
= $Z_3 (= Z_2)$.

Thus the resultant is a fixed point.
    Having obtained the views of experts now we proceed on to find the consolidated view of them and find the effect of state vectors on this combined dynamical system.

Let $M = M_1 + M_2 + M_3 + \ldots + M_{10}$ i.e., we add the 10 matrices where the first column corresponds to the node 1 and first row of all the 10 matrices correspond to node 1. Now we divide each and every term of the matrix M by 10 we obtain a matrix, which is a not a simple FCM, the entries invariably in the matrix $\dfrac{M}{10}$ are values from the interval [0, 1].

$$M = \begin{array}{c} \\ P_1 \\ P_2 \\ P_3 \\ P_4 \\ P_5 \\ P_6 \\ P_7 \\ P_8 \\ P_9 \end{array} \begin{array}{c} P_1 \; P_2 \; P_3 \; P_4 \; P_5 \; P_6 \; P_7 \; P_8 \; P_9 \\ \begin{bmatrix} 0 & 10 & 4 & 4 & 3 & 5 & 5 & 1 & 6 \\ 9 & 0 & 2 & 6 & 1 & 4 & 1 & 1 & 3 \\ 2 & 1 & 0 & 3 & 4 & 1 & 3 & 4 & 7 \\ 5 & 5 & 3 & 0 & 5 & 7 & 5 & 7 & 6 \\ 1 & 0 & 2 & 3 & 0 & 6 & 4 & 2 & 6 \\ 2 & 3 & 1 & 1 & 2 & 0 & 5 & 2 & 6 \\ 1 & 1 & 1 & 1 & 2 & 5 & 0 & 5 & 7 \\ 0 & 0 & 2 & 2 & 0 & 0 & 3 & 0 & 4 \\ 3 & 2 & 2 & 1 & 0 & 6 & 6 & 3 & 0 \end{bmatrix} \end{array}$$

Let $M/10 = N$, N is FCM; which is not simple for the entries belong to [0, 1]

$$N = \begin{array}{c} \\ P_1 \\ P_2 \\ P_3 \\ P_4 \\ P_5 \\ P_6 \\ P_7 \\ P_8 \\ P_9 \end{array} \begin{array}{c} P_1 \; P_2 \; P_3 \; P_4 \; P_5 \; P_6 \; P_7 \; P_8 \; P_9 \\ \begin{bmatrix} 0 & 1 & 0.4 & 0.4 & 0.3 & 0.5 & 0.5 & 0.1 & 0.6 \\ 0.9 & 0 & 0.2 & 0.6 & 0.1 & 0.4 & 0.1 & 0.1 & 0.3 \\ 0.2 & 0.1 & 0 & 0.3 & 0.4 & 0.1 & 0.3 & 0.4 & 0.7 \\ 0.5 & 0.5 & 0.3 & 0 & 0.5 & 0.7 & 0.5 & 0.7 & 0.6 \\ 0.1 & 0 & 0.2 & 0.3 & 0 & 0.6 & 0.4 & 0.2 & 0.6 \\ 0.2 & 0.3 & 0.1 & 0.1 & 0.2 & 0 & 0.5 & 0.2 & 0.6 \\ 0.1 & 0.1 & 0.1 & 0.1 & 0.2 & 0.5 & 0 & 0.5 & 0.7 \\ 0 & 0 & 0.2 & 0.2 & 0 & 0 & 0.3 & 0 & 0.4 \\ 0.3 & 0.2 & 0.2 & 0.1 & 0 & 0.6 & 0.6 & 0.3 & 0 \end{bmatrix} \end{array}$$

Now consider the state vector

$$X = (0\ 0\ 0\ 0\ 0\ 1\ 0\ 0\ 0).$$

Only the node (6) is in the ON state and all other nodes are in the OFF state.

The effect of X on N using the max, min composition rule.

$$\begin{aligned}
XN &= (0.2, 0.3, 0.1, 0.1, 0.2, 0, 0.5, 0.2, 0.6) \\
&= X_1 \text{ (say)} \\
X_1 N &= (0.3, 0.2, 0.2, 0.3, 0.2, 0.6, 0.6, 0.5, 0.5) \\
&= X_2 \text{ (say)}. \\
X_2 N &= (0.3, 0.3, 0.3, 0.3, 0.3, 0.5, 0.5, 0.5, 0.6) \\
&= X_3 \\
X_3 N &= (0.3, 0.3, 0.3, 0.3, 0.3, 0.5, 0.5, 0.5, 0.6) \\
&= X_4 = (X_2).
\end{aligned}$$

Thus we get the fixed point and all the nodes come to ON state.

Now we study the effect of

$$Y = (1\ 0\ 0\ 0\ 0\ 0\ 0\ 0\ 0)$$

on the system N.

$$\begin{aligned}
YN &= (0, 1, 0.4, 0.4, 0.3, 0.5, 0.5, 0.1, 0.6) \\
&= Y_1 \text{ (say)} \\
Y_1 N &= (0.9, 0.3, 0.3, 0.4, 0.4, 0.6, 0.6, 0.5, 0.5) \\
&= Y_2 \\
Y_2 N &= (0.4, 0.9, 0.4, 0.4, 0.4, 0.5, 0.5, 0.5, 0.6) \\
&= Y_3 \\
Y_3 N &= (0.9, 0.4, 0.4, 0.6, 0.4, 0.6, 0.6, 0.5, 0.5) \\
&= Y_4 \\
Y_4 N &= (0.4, 0.9, 0.4, 0.4, 0.5, 0.6, 0.5, 0.6, 0.6) \\
&= Y_5 \\
Y_5 N &= (0.9, 0.4, 0.4, 0.4, 0.4, 0.6, 0.6, 0.5, 0.5) \\
&= Y_6.
\end{aligned}$$

Thus it fluctuates in which case only upper bounds would be taken to arrive at the result.

Now we proceed on to study the effect of

$$Z = (0\ 0\ 0\ 1\ 0\ 0\ 0\ 0\ 0)$$

on N.

$$\begin{aligned}
ZN &= (0.5, 0.5, 0.3, .0, 0.5, 0.7, 0.5, 0.7, 0.6) \\
&= Z_1 \text{ (say)} \\
Z_1 N &= (0.5, 0.5, 0.4, 0.5, 0.3, 0.6, 0.6, 0.5, 0.6) \\
&= Z_2 \text{ (say)} \\
Z_2 N &= (0.5, 0.5, 0.4, 0.5, 0.5, 0.6, 0.6, 0.5, 0.6) \\
&= Z_3 \text{ (say)} \\
Z_{30} N &= (0.5, 0.5, 0.4, 0.5, 0.5, 0.6, 0.6, 0.5, 0.6) \\
&= Z_4 (= Z_3).
\end{aligned}$$

Thus we arrive at a fixed point and all nodes come significantly to a value in [0 1].
Let
$$T = (0\ 1\ 0\ 0\ 0\ 0\ 0\ 0\ 0).$$

The effect of T on N is given by

$$\begin{aligned}
TN &= (0.9, 0, 0.2, 0.6, 0.3, 0.4, 0.1, 0.1, 0.3) \\
&= T_1 \text{ (say)} \\
T_1 N &= (0.6, 0.9, 0.4, 0.6, 0.5, 0.6, 0.5, 0.6, 0.6) \\
&= T_2 \text{ (say)} \\
T_2 N &= (0.9, 0.6, 0.4, 0.6, 0.5, 0.6, 0.6, 0.6, 0.6) \\
&= T_3 \text{ (say)} \\
T_3 N &= (0.6, 0.9, 0.4, 0.6, 0.5, 0.6, 0.5, 0.6, 0.6) \\
&= T_4 \text{ (say)} \\
T_4 N &= (0.9, 0.4, 0.4, 0.6, 0.5, 0.6, 0.6, 0.6, 0.5) \\
&= T_5 = (T_3).
\end{aligned}$$

We see the resultant is a limit cycle fluctuating between $T_3$ and $T_5$. Now consider the state vector

V        =        (0 0 0 0 0 0 0 0 1).

The effect of V on N is given by

| | | |
|---|---|---|
| VN | = | (0.3, 0.2, 0.2, 0.1, 0, 0.6, 0.6, 0.3, 0) |
| | = | $V_1$ (say) |
| $V_1N$ | = | (0.2, 0.3, 0.3, 0.3, 0.3, 0.5, 0.5, 0.5, 0.6) |
| | = | $V_2$ (say) |
| $V_2N$ | = | (0.3, 0.3, 0.3, 0.3, 0.3, 0.6, 0.6, 0.5, 0.5) |
| | = | $V_3$ (say) |
| $V_3N$ | = | (0.3, 0.3, 0.3, 0.3, 0.3, 0.5, 0.5, 0.5, 0.6) |
| | = | $V_4$ (say) |
| $V_4N$ | = | (0.3, 0.3, 0.3, 0.3, 0.3, 0.6, 0.6, 0.5, 0.5) |
| | = | $V_5$ (say = $V_3$). |

Thus the resultant is a fixed point. Now we work with the state vector

W        =        (0 0 0 0 0 0 1 0 0).

Now we study the effect of W on the system N.

| | | |
|---|---|---|
| WN | = | (0.1, 0.1, 0.1, 0.1, 0.2, 0.5, 0, 0.5, 0.7) |
| | = | $W_1$ (say) |
| $W_1N$ | = | (0.3, 0.3, 0.3, 0.3, 0.2, 0.6, 0.6, 0.3, 0.5) |
| | = | $W_2$ (say) |
| $W_2N$ | = | (0.3, 0.3, 0.3, 0.3, 0.3, 0.5, 0.5, 0.5, 0.6) |
| | = | $W_3$ (say) |
| $W_3N$ | = | (0.3, 0.3, 0.3, 0.3, 0.3, 0.6, 0.6, 0.5, 0.5) |
| | = | $W_4$ (say) |
| $W_4N$ | = | (0.3, 0.3, 0.3, 0.3, 0.3, 0.5, 0.5, 0.5, 0.6) |
| | = | $W_5$ (= $W_3$) |

Thus the resultant is again a fixed point

Now we see the

Min of row 1        =        0.1
Min of column 1        =        0.0

| | | |
|---|---|---|
| Min of row 2 is | = | 0.1 |
| Min of column 2 | = | 0.0 |
| Min of row 3 | = | 0.1 |
| Min of column 3 is | = | 0.1 |
| Min of row 4 is | = | 0.3 |
| Min of column 4 is | = | 0.1 |
| Min of row is 5 | = | 0.0 |
| Min of column 5 is | = | 0.0 |
| Min of row six is | = | 0.1 |
| Min of column 6 | = | 0.1 |
| Min of row 7 is | = | 0.1 |
| Min of column 7 is | = | 0.1 |
| Min of row 8 is | = | 0.0 |
| Min of column 8 is | = | 0.1 |
| Min of row 9 | = | 0.0 |
| Min of column 9 | = | 0.3. |

Now one can compare and see the resultant.

For in case of the resultant vector W when the node 7 is the ON state i.e. W = (0 0 0 0 0 0 1 0 0); we see the resultant is (0.3, 0.3, 0.3, 0. 3, 0.5, 0.5, 0.5, 0.5) the nodes 6, 7, 8 and 9 take value 0.5, and these three nodes 6, 8 and 9 are equally affected and also the nodes 1, 2, 3, 4 and 5 are affected and all of them to the same degree taking the value 0.3.

Likewise one can make observations about the state vectors X, Y, Z V and T and arrive at conclusions. However we have worked out these conclusions and have put them under the title 'observations' in the last chapter of this book. We requested the experts that if they had any form of dissatisfaction while giving membership to the nodes and if they felt in some cases the relation (i.e., membership grade) was an indeterminate they can use NCMs and described it to them (Section 3.4). A few agreed to work with it. Majority of them did not wish to work with it. However we have used the NCM models given by them and worked with the state vectors given by them and included the analysis in chapter 5. Now the working is identical with that of FCMs. Here we give a typical associated neutrosophic matrix of the NCM given by an expert.

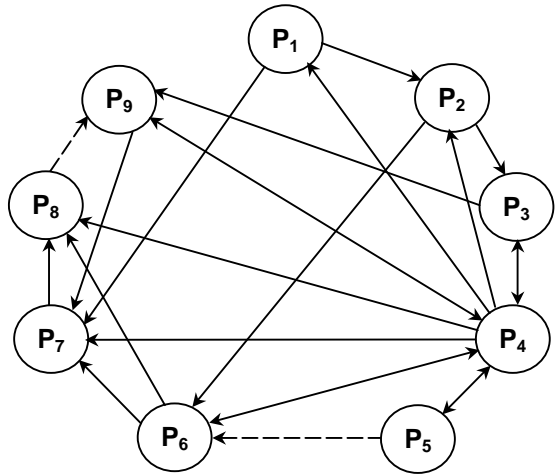

FIGURE 4.5.11

$$M_n = \begin{array}{c} \\ P_1 \\ P_2 \\ P_3 \\ P_4 \\ P_5 \\ P_6 \\ P_7 \\ P_8 \\ P_9 \end{array} \begin{array}{c} P_1 \; P_2 \; P_3 \; P_4 \; P_5 \; P_6 \; P_7 \; P_8 \; P_9 \\ \begin{bmatrix} 0 & 1 & 0 & 0 & 0 & 0 & 1 & 0 & 0 \\ 0 & 0 & 1 & 0 & 0 & 1 & 0 & 0 & 0 \\ 0 & 0 & 0 & 1 & 0 & 0 & 0 & 0 & 1 \\ 1 & 1 & 1 & 0 & 1 & 1 & 1 & 1 & 1 \\ 0 & 0 & 0 & 1 & 0 & I & 0 & 0 & 0 \\ 0 & 0 & 0 & 1 & 0 & 0 & 1 & 1 & 0 \\ 0 & 0 & 0 & 0 & 0 & 0 & 0 & 1 & 0 \\ 0 & 0 & 0 & 0 & 0 & 0 & 0 & 0 & I \\ 0 & 0 & 0 & 1 & 0 & 0 & 1 & 0 & 0 \end{bmatrix} \end{array}$$

Now we study the effect of $X = (1\ 0\ 0\ 0\ 0\ 0\ 0\ 0\ 0)$ on $M_n$ i.e., only the node 'when they claim Vedic Mathematics is magic has more ulterior motives' is in the ON state and all other nodes are in the OFF state. Effect of X on the dynamical system $M_n$ is given by

$$XM_N \rightarrow (1\ 1\ 0\ 0\ 0\ 0\ 1\ 0\ 0)$$
$$= X_1 \text{ (say)}$$
$$X_1 M_N \rightarrow (1\ 1\ 1\ 0\ 0\ 1\ 1\ 1\ 0)$$
$$= X_2.$$
$$X_2 M_N \rightarrow (1\ 1\ 1\ 1\ 1\ 1\ 1\ 1\ I)$$
$$= X_3.$$
$$X_3 M_N \rightarrow (1\ 1\ 1\ 1\ 1\ I\ 1\ 1\ I)$$
$$= X_4 \text{ (say)}.$$
$$X_4 M_N \rightarrow (1\ 1\ 1\ 1\ 1\ I\ 1\ 1\ I).$$

Thus the hidden pattern of the dynamical system is a fixed point which is interpreted as: "if the node Vedic Mathematics is 'magic' then Vedic Mathematics has more ulterior motives" alone is in the ON states all nodes come to the ON state except the nodes 6 and 9 "It is a means to globalize Hindutva" and "It is a more a political agenda to rule the nation and if Sanskrit literature is lost it could produce peace in the nation" alone are in the indeterminate state. Likewise the dynamical system $M_n$ can be worked with any node or nodes in the ON states and the resultant effect can be derived!

**Chapter Five**

# OBSERVATIONS

This chapter gives the observations that were obtained from our mathematical research. It is listed under 5 heads. In the first section we give the views of students and the observations made by the teachers is given in section two. Section three gives the views of the parents, and observations of the educated elite are given in section four. Public opinion is recorded in section five.

## 5.1 Students' Views

1. Almost all students felt that Vedic Mathematics has no mathematical content except at the level of primary school arithmetic.

2. All of them strongly objected to the fact that Vedic Mathematics classes wasted their time.

3. None of the students ever felt that Vedic Mathematics would help them in their school curriculum.

4. Many students said that in this modernized world, Vedic Mathematics was an utter waste because calculators could do all the arithmetical tricks given in that textbook in a fraction of a second.

5. Students criticized heavily that they were forced to learn by rote topics like the Vedic Mathematics with its 16 sutras in these days of globalization and modernization. Without any mathematical significance, just reading these sutras made them feel as if they were the laughing stock of the world.

6. Non-Hindu students felt it difficult to accept the subject, because they were made to feel that they have to be Hindus to read Vedic Mathematics. For instance, the cover of the two Vedic Mathematics books (Books 1 and 2) in Tamil had the picture of Hindu Goddess of Learning, Saraswathi [85-6]. Some of the parents objected because they did not want their children to be forcefully made to take up some other religion using mathematics.

7. Some students frankly said, "our younger brothers and sisters will be made to attend classes on Vedic chemistry, Vedic physics, Vedic zoology, Vedic history, Vedic geography and so on. As our main aim was to obtain their unrestricted views we did not curtail them and in fact recorded the height of their creative imagination!

8. A group of boys said, "Give us just one day's time, we will also write one problem like Swamiji and give a mental solution in a line or two." Students of one particular school said that Mohan, their class topper in Mathematics placed one such simple elementary arithmetic problem and a single line solution within a span of five minutes, and he had told that this is his own Vedic Mathematics for fun. Their teacher got furious and slapped him. The students said, "We all thought the Vedic Mathematics teacher will praise him but his action made us hate Vedic Mathematics all the more. We also hated the meaningless 'sutras', which has nothing in it." Their contention was that everyone could invent or write such sutras, which are very simple and have no Vedic notions about it. They felt that everything was so simple and unscientific, and just $5^{th}$ standard mathematics was sufficient to invent these problems and sutras. They even said that they could invent any form of word in 'Sanskrit' and say it means such-and-such-a-thing; and they came up with some Sanskrit sounding names that could not be easily pronounced!

9. In conclusion, over 90% of the students visibly showed their rationalistic views on the subject and condemned Vedic Mathematics as useless. They felt it would do only more harm to them than any good because they feel that their scientific temperament is caged by being made to repeat sutras that they really do not understand. They said that at least when they repeated rhymes in UKG or LKG they knew at least 90% of the meaning, but this one or two-word Sanskrit sutras never conveyed anything to them, mathematical or scientific. They said, "to please our teacher we had to do the monkey tricks. When the language of communication in the classroom is English what was the relevance of the 16 sutras in Sanskrit, which is an alien language to us and does not convey any meaning?" Even French or German (that are foreign languages) was more appealing to these students than these sutras that they treated with utmost contempt. The younger generation was really very open-minded and frank in its views and choices. They were not clouded by caste or religion. They exhibited a scientific approach which was unbiased and frank!

10. We also met a group of 9$^{th}$ class students who were undergoing Vedic Mathematics training. We asked them to give their true feelings. Most of them said that it was boring compared to their usual mathematics classes. Several of the students strongly disposed of the idea because when have mini-calculators to help them with calculations why did they need Vedic Mathematics for simple multiplication? But any way we have to waste money both buying the book as well as waste time by attending the classes. It would be better if they teach us or coach us in any of the entrance test than in making us study this bore; was the contention of the majority. Some said our parents have no work they in their enthusiasm have even bought the teachers manual for us but we see manual is more interesting with pictures; for when we see book it is just like a primary school mathematics text. This with calculators in hand we don't need all this for our career they said in a single voice.

11. None of the rural school students have heard about Vedic Mathematics. When we illustrated certain illustrations from this book, a few of them said that their mathematics teachers knew much more simpler methods than the ones shown by us. Most of the rural mathematics teachers were unaware about the Vedic Mathematics book. A few of them did know more simple and easy calculations than the ones given in that book. The teachers said that if 'multiplication tables' were taught in the primary class and more arithmetic problems given, then students themselves would invent more such formulae. Awareness about Vedic Mathematics was almost totally absent. In rural areas, the question about parents' opinion does not come up because they are either uneducated or totally ignorant of the book on Vedic Mathematics. They are involved in the struggle to make both ends meet to support the education of their children.

12. A 9-year-old boy from a very remote village claims that he has never heard about Vedic Mathematics, but however wanted to know what it was. He asked us whether it was taught in Sanskrit/ Hindi? When we explained one or two illustrations, within 10 minutes time he came to us and said that he has discovered more such Vedic Mathematics and said he would give answer to all multiplication done by 9, 99 and 999 mentally. We were very much surprised at his intelligence. From this the reader is requested to analyze how fast he has perceived Vedic Mathematics. Further each person has a mathematical flair and his own way of approach in doing arithmetical problems, especially addition, multiplication and division. In fact if such a boy had been given a week's time he would have given us more than 10 such sutras to solve arithmetic problems very fast. He said he did not know Sanskrit or Hindi or English to name the sutras in Sanskrit. On the whole, students of government-run corporation schools were bright and quick on the uptake but fortunately or unfortunately they have not heard or seen any book on Vedic Mathematics. Might be most of the students who study in such schools are the lesser children of God, so Vedic Mathematics has not yet reached them or the school authorities.

## 5.2 Views of Teachers

1. "As a Sankaracharya, who is a Hindu religious leader, wrote the book, neither the mathematical community nor the teachers had the courage to refute it. But we had to accept it as a great work," says one teacher. He continues, "If a teacher like me had written a book of this form, I would have been dismissed from my job and received a mountain of criticism which I would not be in a position to defend." Thus when a religious man professes foolish things, Indians follow it just like goats and are not in a position to refute it. It is unfortunate that Indians do not use reasoning mainly when it comes from the mouth of a religious leader. Thus he says this book is an insane method of approaching mathematics because even to multiply 9 by 7 he uses several steps than what is normally required. Thus, this retired teacher, who is in his late sixties, ridicules this book.

2. Most of the mathematics teachers in the 50+ age group are of the opinion that while doing arithmetical calculations the teachers' community uses most of the methods used by the book of Vedic Mathematics. They claim each of them had a shortcut method, which was their own invention or something which they had observed over years of practice. So they just disposed of the book Vedic Mathematics as only a compilation of such methods and said that it has nothing to do with Vedas. Because the Jagadguru Sankaracharya was a religious man, he had tried to give it a Vedic colour. This has faced criticism and ridicule from mathematicians, students and teachers.

3. Similarly, a vast majority of the teachers felt that a group of people has made a lot of money by using this book. They further feel that such methods of simplification are of no use in the modernized world where calculators can do the job in a fraction of a second. They felt that instead of teaching haphazard techniques, it would be better to teach better mathematics to children who fear mathematics. In their opinion, most rural children do well in mathematics. But these methods of Vedic Mathematics will certainly not wipe out fear from their mind but only further repel them from mathematics. They are of the opinion that the Swamiji who has studied up to a M.Sc. or M.A. in Mathematics did not show any talent but just the level of a middle school mathematics teacher. They still felt sad because several parents who do not have any knowledge of mathematics force their children to read and solve problems using the methods given in the Vedic Mathematics book.

4. We discussed about Jagadguru Swami Sankaracharya of Puri with a Sanskrit pundit (now deceased) hailing from Tanjore who had served as headmaster, and was well versed in Sanskrit and Hindi and had even worked in the Kanchi mutt. At the first place, he came down heavily on the Swami Sankaracharya of Puri because he had crossed the seas and gone abroad which was equivalent to losing ones caste. He cited the example of how the Sankaracharya of Kanchi was not permitted by other religious leaders to visit China or even Tibet. Under these conditions, his visit abroad, that too, to the Western countries under any pretext was wrong and against all religious dharma. He said that the Vedic Mathematics book written by that Sankaracharya was humbug. He said that as a retired headmaster he also knew too well about the mathematics put forth in that book. He said that being a Sanskrit scholar he too could give some sutras and many more shortcuts for both multiplication and division. He asked us, "Can my sutras be appended to the Vedas?" He was very sharp and incomprehensible so we could only nod for his questions. Finally he said that the Sankaracharya of Puri failed to give any valuable message of Vedanta to the people and had wasted nearly 5 decades. It is pertinent to mention here that this Pundit's father had taught all the Vedas to the Sankaracharyas of the Kanchi Math.

5. A mathematics teacher with over 30 years of experience and still in service made the following comments: He said that he has seen the three books Book 1, Book 2 and Book 3 of Vedic Mathematics for schools [148-150, 42-4]. He adds he has also seen the Vedic Mathematics teachers' manual level I, II and III. The Indian edition of the teachers' manual appeared only in the year 2005. He has read all these books. He asked us why Vedic Mathematics books were written first and only recently the Teachers Manual was written. Why was the procedure topsy-turvy? Does Vedic Mathematics teach topsy-turvy procedure? Secondly he says he is utterly displeased to see that the foreword for all these six books was given by Dr. L.M. Singhvi, High Commissioner for India in the UK. Does he hold a doctorate in mathematics? What made him give preface or foreword to all these books? What has made him appreciate Vedic Mathematics: Is it Vedas? Is it the Jagadguru? Or does the publisher try to get some popularity and fame in the west by choosing the Indian High Commissioner in the UK to give the foreword/ preface? If a mathematical expert had reviewed the book in the foreword/ preface it would have been 100 times more authentic. It actually seems to hold ulterior motives. The teacher points out that one of the obvious factors is that Dr. Singhvi writes in his foreword in the Vedic Mathematics Teachers' Manual [148], "British teachers have prepared textbooks of Vedic Mathematics for British schools. Vedic Mathematics is thus a bridge across countries, civilizations, linguistic barriers and national frontiers." This teacher construes that being a High Commissioner Dr. Singhvi would have had a major role in propagating Vedic Mathematics to British schools. The teacher said, "when Vedic scholars (i.e. the so-called Brahmins) do not even accept the rights of Sudras (non-Brahmins) and ill-treat them in all spheres of life and deny them all economical, social, religious, educational and political equality, it is a mockery that Dr. L.M Singhvi says that Vedic Mathematics is a bridge across civilizations and linguistic barriers. They have always spoken not only about their superiority but also about the superiority of their Sanskrit language. When they cannot treat with equality other Indians with whom they have lived for so many centuries after their entry into India by the Khyber Pass how can this book on Vedic Mathematics now profess equality with British, whom they chased out of India at one point of time?" He claims that all this can be verified from the books by Danasekar, Lokamanya Thilak Popular Prakashan, Bombay, p.442 and Venkatachalapathi, VOC & Bharathi p.124, People Publication, Chennai 1994. He feels that Vedic Mathematics is a modern mathematical instrument used by a section of the so-called upper castes i.e. Brahmins to make India a Hindu land and Vedic Mathematics would help in such a Hindu renaissance. The minute somebody accepts Vedic Mathematics, it makes him or her unconditionally accept Hinduism and the Hindu way of life. Certainly modern youth will not only be cheated but they will have to lead a life of slavery, untouchability and Sudrahood. So this teacher strongly feels Vedic Mathematics is a secret means to establish India as a Hindutva land.

6. Next, we wanted to know the stand of good English medium schools in the city that were run by Christian missionaries. So we approached one such renowned school. We met the Principal; she said she would fix an appointment for us with her school mathematics teachers. Accordingly we met them and had many open discussions. Some of the nuns also participated in these discussions. Their first and basic objection was that Vedic Mathematics was an attempt to spread Hindutva or to be more precise Brahminism. So they warned their teachers and students against the use of it. They criticized the cover-page of the Vedic Mathematics books in Tamil. The cover page is adorned with a picture of Saraswathi, the Hindu goddess of education. She has four hands and holds a veena. Underneath the photo a Sanskrit sloka is written in Tamil that prays for her blessings. The first question they put to us was, "Is Vedic Mathematics Hindutva mathematics? It would be more appropriate if they could call it "Hindutva Mathematics" because it would not be misleading in that case. Is it for unity or for diversity? Can a Muslim or a Christian be made to accept the cover of the book? Have you ever seen a mathematics book with the cover page of Jesus or Mohammad or Mary? How can Vedic Mathematics books have such a cover if they are really interested in spreading mathematics for children? Their main mission is this: They have come to know that because of the lack of devoted teachers in the recent days, mathematics has become a very difficult subject especially in private non-government city schools. The present trend of parents and students is to get good marks and get a seat in a good professional institution. So, to capture both the students and parents in the Hindutva net, they have written such books with no mathematical value." Then they said that there are many good mathematics teachers who do more tricks than the tricks mentioned in the Vedic Mathematics books. They also started criticizing the 'trick' aspects of mathematics. They asked, "Can a perfect and precise science like mathematics be studied as lessons of trick? How can anyone like a subject that teaches performing tricks? If somebody dislikes performing tricks or does not know to perform such tricks can he or she be categorized as a dull student? If one accepts Vedic Mathematics, he accepts his Hindu lineage thereby he becomes either a Sudra or an Untouchable? Can they apply the universalism that they use for Vedic Mathematics and declare that the four Varnas do not exists, all are equal and that no caste is superior? Are we Christians from Europe? We were the true sons of the Indian soil and were forced into embracing Christianity because we were very sensitive and did not want to accept ourselves as Sudras or Untouchables. We wanted to say to them that we were equal and in fact superior to the Brahmins. Our self-respect prompted us to become Christians. So Vedic

Mathematics is only an instrument to spread Hindutva and not mathematics. Also the mathematics given in Vedic Mathematics is of no use because our school children are brighter and can invent better shortcut methods to arithmetic than what is given in that book." Finally they asked us whether any relation existed between Motilal Banarsidass and the author of the book because the company seems to have made a lot of money selling these books?

We then visited a reputed boys school run by a Christian missionary. We had a four-hour long discussion with mathematics teachers of that school. The principal and the vice principal were also present. They had a collection of 8 books displayed on the table: Vedic Metaphysics, Vedic Mathematics, Book 1, Book 2, Book 3 of Vedic Mathematics for Schools, Vedic Mathematics Teacher's Manual, for the elementary level, intermediate level and advanced level.

From the intermediate level teachers' manual, they showed us p.145 of [51].

$$13.\ \text{Solve}\quad x + y = 6$$
$$x - y = 2$$

"The formulae by addition and by subtraction and by alternate elimination and retention can be used to solve simultaneous equations." Everyone said that such trivial equations could be solved mentally and need not find its place in the Teachers' manual for the intermediate level! [150] "If a teacher solves or gives hints to solve this problem the way it is described in page 145 of that manual, he will be sent home by my students the same day," said the principal of the school.

Next, they showed us an example from p.30 of the same book [51] *Nikhilam Navatascharaman Dasatah* (All from 9 and the last from 10) is (14) $88 \times 98$. They said that such mental calculations are done at the primary school level and need not find place in the teachers' manual. They also added, "We have hundreds of such citations from the three books of the Vedic Mathematics teachers' manual—all of them are substandard examples."

The principal said vehemently, "We have kept these books as if they are specimen items in a museum and are not for educational use. In the first place, the Vedic Mathematics book has no mathematical value and secondly it imparts not mathematics but only destructive force like casteism. For instance, it is said in the book, "*Vedic Mathematics is not a choice for slow learners. It demands a little briskness.*" So, the Brahmins will go on to say that all of the other castes are slow learners, and they might declare that we cannot read mathematics.

I remember what a student here mentioned about a Brahmin teacher in his previous school who had said: "Even if the Durba Grass is burnt and kept into the tongue of the Sudras, then also they cannot get mathematics." Why do they write Vedic Mathematics books for school children? Is it not the height of arrogance and cunning to declare first that Vedas cannot be imparted to non-Brahmins, so also Vedic Mathematics? One should analyze Vedic Mathematics, not as a mathematics book but for its underlying caste prejudice of Vedas ingrained in it. As a mathematics book even a $10^{th}$ grader would say it is elementary!"

He continued, "Can a Christian pontiff write a Christian mathematics book for school children stating a few Hebrew phrases and say that they mean "one less than the existing one" "one added to the previous one" and so on. Will Hindus all over the world welcome it? Suppose we put the cover picture of Jesus or Mary in that Christian mathematics book what will be their first reaction? They will say, "Christian fanatics are trying to spread Christianity; in due course of time India would become a Christian nation, so ban the book." Likewise, if a Maulana writes a book on Islam mathematics saying some words in the Kuran are mathematical sutras; what will be the Brahmins' reactions? They will say, "The nation is at stake. Terrorism is being brought in through mathematics. Ban the book, close down all minority institutions. Only Hindu institutions should be recognized by the Government". "If the Hindutva Government was in power, Government Orders would have been passed to this effect immediately." So, in his opinion Vedic Mathematics has no mathematical content. Secondly, it is of no use to slow learners (this is their own claim) so in due course of time it would be doing more harm to people than any good. Thirdly, it is a sophisticated tool used to reestablish their lost superiority and identity.

7. Next, we discussed the Vedic Mathematics Teachers' manuals with a group of school teachers. They put forth the following points:

   1. The manuals cost Rs.770/- totally. They are so highly priced only to make money and not for really spreading Vedic Mathematics.
   2. The intermediate manual itself looks only like primary school mathematics.

   One example given from the manual [p.3, Intermediate, 150]: Finding digit sum i.e. digit sum of 42 is 6 is first practice given in the manual. There follows very simple first standard addition and multiplication up to p.43 [149-150]. Then there is simple primary school division. There ends the teachers' manual for the intermediate level. When we come to Vedic Mathematics Teachers' Manual Advanced level we have the following: First few sections are once again primary school level addition, multiplication, division and subtraction. Solution to equation page 79 is nothing more than what the usual working does. So is the following exercise [148]. Page 126, osculation [148]. Find out if 91 is divisible by 7. The method by Ekadhika is longer and cumbersome than the usual long division of 91 by 7. Now we come to analyze Vedic Mathematics Teachers' Manual of elementary level [149]. Page 98 Vedic Mathematics [51] The first by the first and the last by the last. He says $27 \times 87$ = 23/49. The condition are satisfied here as $2 + 8 = 10$ and both numbers end in 7. So we multiply the first figure of each number together and add the last figure. $2 \times 8 = 16$, 16

+ 7 = 23 which is the first part of the answer. Multiplying the last figures together 7 × 7 = 49, which is the last part of the answer. The teacher feels the same method cannot be applied for finding the value of

47 × 97 for 47 × 97 ≠ 43/49
    43 / 49 is got by applying the formula
4 × 9 = 36, 36 + 7 = 43
    7 × 7 = 49 so 43/49. The true value of 47 × 97 = 4559

So the formula cannot be applied. Everyone can find product 27 × 87 and 47 ×97 if they remember that

1. One condition is the first figures should add to 10
2. The 2$^{nd}$ digit must be the same.

    How will a student remember this while carrying out multiplication that too only by two digits in an exam hall? The product is not defined for three digit × two digit or four digit × two digit …
    How could one claim that Vedic Mathematics is fast and wipes out fear in students? As teachers we feel it is not only a waste of time but will also scare children from mathematics because it requires more memory than intelligence whereas the reverse is required for mathematics. Thus true intelligence will be lost in children. Also the sharpness of the mind is at stake by teaching them Vedic Mathematics.

8. Next, we met the teachers working in a school run by a Muslim minority educational trust. There were 6 mathematics teachers: one Muslim woman, the rest were Hindus. At the first instance, all of them said it would be appropriate to term it Hindutva mathematics because the term 'Vedic Mathematics' was a misnomer because in Vedic times no one would have had the facility or time or above all the need to find the values of1/17 or 1/19 and so on.

In the second place, when Vedas are thought to be so religious that they should be read only by the Brahmins; how is it that such trivial arithmetic is included in it? Above all, why should Shudras read this trivial mathematics today? Will not this pollute the Vedas and the Vedic principle? The Muslim lady teacher said that if a Maulana came up with these simple arithmetic formulae after some eight minutes of meditation, they would say he was mad and send him to Erwadi. She wondered how he could occupy the highest place and be the Jagadguru Puri Sankaracharya. In her opinion, their religious leaders hold a high place and by no means would they poke their nose into trivialities like easy arithmetic for school children. They only strive to spread their religion and become more and more proficient in religious studies.

All the teachers had a doubt whether the Swamiji wanted to spread Hindutva through Vedic Mathematics? They asked will we soon have Islamic mathematics, Christian mathematics, and Buddhist mathematics in India. Another teacher pointed out, "Will any secular/ common Mathematics book be adorned on the cover page with Goddess Saraswathi? Is this not proof enough to know whose mission is Vedic Mathematics?" They all concluded our brief interview by saying, "We don't follow any trash given in that book because it has no mathematical content. We have many more shortcuts and easy approaches to solving problems."

9. We met a group of teachers who are believers in the ideology of Tamil rationalist leader Periyar and his self-respect movement. Some of them were retired teachers, while others where still in service. These teachers were very angry about Vedic Mathematics. They were uniformly of the opinion that it was a means to spread Hindutva. They claim that in due course of time, these people may even forbid non-Brahmins from reading Mathematics just as they forbade them from learning the Vedas. Also in due course of time they may claim mathematics itself as a Vedanta and then forbid non-Brahmins from learning it. The book has no mathematical content and only a religious mission viz. spread of Hindutva. That is why simple things like addition, subtraction and multiplication are given the name of 'Vedic Mathematics'. They felt that anyone who accepts Vedic Mathematics accepts Hindutva.

They proceeded to give us examples

(1) 'Antyayoreva' – only the last digits
(2) 'Vilokanam' – by mere inspection
(3) 'Paravartya Yojayet' – transpose and 'adjust'
(4) "Nikhilam Navatascaraman Dasatah' – All from nine and the last from ten –

They began their arguments in an unexpected angle: "Suppose we write such sutras in Tamil, what will be our position? Who will accept it? Why are the non-Brahmins who are the majority so quiet? Prof. Dani was great to warn us of the stupidity of Vedic Mathematics and appealed to the saner elements to join hands and educate people on the truth of this so-called Vedic Mathematics and prevent the use of public money and energy on its propagation. He said it would result in wrong attitudes to both history and mathematics especially where the new generation was concerned."

"Periyar has warned us of the cunningness of Brahmins, so we must be careful! It is high time we evaluate the Vedic Mathematics and ban its use beyond a limit because 'magic' cannot be mathematics. The tall claims about Vedic Mathematics made by some sections like applying it to advanced problems such as Kepler's problem etc. are nothing more than superficial tinkering. It offers nothing of interests to professionals in the area."

Then they said, "Why did it take nearly a decade for a Swamiji to invent such simple sutras in arithmetic? Sharma says "intuitional visualization of fundamental mathematical truths born after eight years is the highly concentrated endeavour of Jagadguru"—but does anyone have to spend such a long time.

10. One woman teacher spoke up: "As teachers we feel if we spell out the sutras like 'By mere inspection', 'only the last digits', our students will pelt stones at us in the classroom and outside the classroom. What is a sutra? It must denote some formula. Just saying the words, 'By mere inspection' cannot be called as a sutra! What are you going to inspect? So each and every sutra given by the Swamiji does not look like sutra at all. We keep quiet over this, because even challenging Vedic Mathematics will give undue publicity to Hindutva."

11. The hidden pattern given by dynamical system FRM used by the teachers revealed that the resultant was always a fixed binary pair. In most cases only the nodes Vedic Mathematics is primary level mathematics, Vedic Mathematics is secondary level mathematics, Vedic Mathematics is high school level mathematics and it has neither Vedic value nor mathematical value remained as 0, that is unaffected by the ON state of other nodes because teachers at the first stage itself did not feel that Vedic Mathematics had any mission of teaching mathematics. None of them admitted to finding new short-cuts through the book. Teachers were also very cautious to answer questions about the "Vedic value of Vedic Mathematics" and the "religious value of Vedic Mathematics" for reason best known to them. The study reveals that teachers totally agree with the fact that Vedic Mathematics has a major Hindutva/ Hindu rightwing, revivalist and religious agenda.

## 5.3 Views of Parents

We interviewed a cross-section of parents (of school-going children) for their opinion on Vedic Mathematics.

1. Several parents whose wards were studying in schools run by Hindu organizations spoke of the ill-treatment faced by their children in Vedic Mathematics classrooms. The students were forced to learn Sanskrit sutras by rote and repeat it. Some of them faced difficulties in the pronunciation for which they had been ridiculed by their teachers. Some of the parents even alleged that their children had been discriminated on caste basis by the teachers. One parent reported that after negotiations with some powerful members of the school, she got her child an exemption from attending those classes. She expressed how her son used to feel depressed, when he was ill-treated. She added that because of her son's dark complexion, the teacher would always pounce on him with questions and put him down before his classmates.

2. They uniformly shared the opinion that Vedic Mathematics was more about teaching of Sanskrit sutras than of mathematics, because their children did the problems given in that textbook in no time. Most of the children had told their parents that it was more like primary school mathematics. They said it was just like their primary school mathematics. Yet, the Vedic Mathematics classes were like language classes where they were asked to learn by rote Sanskrit sutras and their meaning.

3. A section of the parents felt that it was more a religious class than a mathematics class. The teachers would speak of the Jagadguru Puri Sankaracharya and of the high heritage of the nation that was contained in Vedas. Actual working of mathematics was very little, so the young minds did not appreciate Vedic Mathematics. Parents expressed concern over the fact that they were compulsorily made to buy the books which cost from Rs.95 to Rs.150. In some schools, in classes 5 to 8 students were given exams and given grades for studying Vedic Mathematics. A few parents said that the classroom atmosphere spoilt their child's mental make-up. Some of them had made their children to switch schools. Thus most non-Brahmins felt that Vedic Mathematics made their children feel discriminated and indirectly helped in developing an inferiority complex.

    A small boy just in his sixth standard had asked his parents what was meant by the word Sudra. Then he had

wanted to know the difference in meaning between the words Sutra (formulae) and Sudra (low caste Hindu). His teacher often said in the classroom that Sudras cannot learn mathematics quickly and to learn Vedic Mathematics one cannot be a slow learner. Thus they felt that caste creeps in indirectly in these Vedic Mathematics classes.

## 5.4 Views of the Educated

We interviewed over 300 educated persons from all walks of life: doctors, judges, senior counsels, lawyers, engineers, teachers, professors, technicians, secretarial workers and psychiatrists. The minimum educational qualification stipulated by us was that they should at least be graduates. In fact several of them were post-graduates and doctorates; some of them were vice-chancellors, directors, educationalists, or employed in the government cadre of Indian Administrative Service (IAS), Indian Revenue Service (IRS) and Indian Police Service (IPS) also.

They showed a lot of enthusiasm about this study, but for their encouragement and cooperation it would not have been possible for us to write this book. Further, they made themselves available for discussions that lasted several hours in some instances. They made many scientific and psychological observations about the effect created by Vedic Mathematics in young minds. Some people said that Vedic Mathematics was an agenda of the right-wing RSS (Rashtriya Swayamsevak Sangh) which planned to 'catch them young' to make them ardent followers of Hindutva. They suggested several points as nodal concepts in our models, we took the common points stressed by several of them. Now we enlist the observations both from the discussions and mathematical analysis done in chapter 4.

1. All of them felt that Vedic Mathematics had some strong ulterior motives and it was not just aimed to teach simple arithmetic or make mathematics easy to school students.

2. Most of them argued that it would create caste distinction among children.

3. All of them dismissed Vedic Mathematics as simple arithmetic calculations!

4. Many of them came down heavily on the Puri Sankaracharya for writing this book by lying that it has its origin in Vedas. All the 16 sutras given in the Vedic Mathematics book had no mathematical content of that sort [31,32].

5. A few of the scholars came down heavily on the title. They felt that when the Vedas cannot be read or even heard by the non-Brahmins, how did Jagadguru Sankaracharya have the heart to write Vedic Mathematics for students when the non-Brahmin population is over 90% in India. They said, "If Vedic Mathematics was really derived from the Vedas, will Brahmins ever share it with others?" Further, they said that Jagadguru Sankaracharya himself was fully aware of the fact that the 16 sutras given by him in pages 17-18 of the book [51] were coined only by him. Those phrases have no deep or real formula value. They were of the opinion that because someone wanted to show that "mathematics: the queen of sciences" was present in the Vedas this book was written. This had been done so that later on they could make a complete claim that all present-day inventions were already a part of the Vedas. But the poor approach of the Jagadguru had made them fail miserably.

6. They were totally against the imposition of Vedic Mathematics in schools run by pro-Hindutva schools. They condemned that teaching Vedic Mathematics also involved discrimination on caste basis. Some backward class and Dalit students were put down under the pretext that they were not concentrating on the subject. Their parents disclosed this during the discussions. Questions like "how many of you do '*Sandhyavadana*'?" were put to the students. Such tendencies will breed caste discrimination.

7. Majority of them did not comment about the Vedic content in Vedic Mathematics but were of the opinion that the 16 sutras were rudimentary and had no relation with Vedas. Further, they agreed that the mathematics described in that book was elementary school arithmetic.

8. All of them agreed upon the fact that Vedic Mathematics had an ulterior motive to establish that Brahmins were superior to non-Brahmins and that Sanskrit was superior to Tamil. This was slowly injected in the minds of the children in the formative age.

9. A section of the interviewed people said that Vedas brought the nasty caste system to India, and they wondered what harm Vedic Mathematics was going to bring to this society. They also questioned the reasons why Vedic Mathematics was being thrown open to everybody, whereas the Vedas had been restricted to the Brahmins alone and the 'lower' castes had been forbidden from even hearing to the recitations.

10. From mathematical analysis we found out that all the educated people felt that Vedic Mathematics was a tool used by the Brahmins to establish their supremacy over the non-Brahmins.

11. Vedic Mathematics was the Hindutva agenda to saffronize the nation.

12. Nobody spoke about Vedic contents in Vedic Mathematics. This node always took only the zero value in our mathematical analysis.

13. In this category, a strong view emerged that Vedic Mathematics would certainly spoil the student-teacher relationship.

14. Some of them said that they are selling this book to make a quick buck and at the same time spread the agenda of Hindutva. Students in urban areas, generally tend to be scared of mathematics. They have exploited this weakness and have aimed to spread Vedic Mathematics.

15. A section of the people interviewed in this category said that Vedic Mathematics was being taught in schools for nearly a decade but has it reduced the fear of mathematics prevailing among students?" The answer is a big NO. Even this year students complained that the mathematics paper in entrance tests was difficult. If Vedic Mathematics was a powerful tool it should have had some impact on the students ability after so many years of teaching.

16. Many of the respondents in this category said that it was very surprising to see Vedic Mathematics book talk of Kamsa and Krishna. Examples cited from the book were: p.354 of the book [51] says, "During the reign of King Kamsa" read a Sutra, "rebellions, arson, famines and unsanitary conditions prevailed". Decoded, this little piece of libelous history gave the decimal answer to the fraction 1/17; sixteen processes of simple mathematics reduced to one." Most of them felt that this is unwarranted in a mathematics text unless it was written with some other ulterior motive. A Sudra king Kamsa is degraded. Can anyone find a connection between modern mathematics and a religious Brahmin pontiff like Sankaracharya of Kanchi? Why should Brahmins find mathematical sutras in sentences degrading Sudras? At least if some poetic allegory was discovered, one can accept it, but it was not possible to understand why the decimal answer to the fraction of 1/17 was associated with a Sudra king Kamsa. Moreover, decimal representation was invented only in the 17th century, so how can an ancient sloka be associated with it? If some old Islamic/Christian phrase as given mathematical background, will it be accepted in India?

17. Some of them wanted to debate the stand of the media with regard to Vedic Mathematics. While a major section of the media hyped it, there was a section that sought to challenge the tall claims made by the supporters of Vedic Mathematics. This tiny section, which opposed Vedic Mathematics, consisted notably of leftist magazines that carried articles by eminent mathematicians like [31-2].

18. People of this category shared a widespread opinion that like the tools of yoga, spirituality, this Vedic Mathematics also was introduced with the motivation of impressing the West with the so-called Hindu traditions. They feared that these revivalists would say that all discoveries are part of the Vedas, or they might go ahead and say that the Western world stole these discoveries from them. They rubbished the claims that the Vedas contained all the technology or mathematics of the world. Already, the Brahmins / Aryans in those ages had appropriated all the indigenous tradition and culture and with a little modification established their superiority. Perhaps Vedic Mathematics is a step in that direction because p.XXXV of the book states, "*(1) The sutras (aphorisms) apply to and cover each and every part of each and every chapter of each and every branch of mathematics (including arithmetic, algebra, geometry – plane and solid, trigonometry plane and spherical, conics – geometrical and analytical, astronomy, calculus – differential and integral etc). In fact there is no part of mathematics pure or applied which is beyond their jurisdiction*". Thus they felt that by such a broad, sweeping statement, the Swamiji had tried his level best to impress everybody about the so-called powers of Vedic Mathematics. A few of them said that pages XXXIII to XXXIX of the book on Vedic Mathematics should be read by everybody to understand its true objective and mission which would show their fanatic nature. They merely called it an effort for the globalization of Hindutva. [51]

19. A Sanskrit Pundit whom we interviewed claims that Swamiji (with his extraordinary proficiency in Sanskrit)

could not invent anything mathematically, so he indulged in extending the Vedas. A similar instance can be the story of how the Mahabharata grew from a couple of hundred verses into tens of thousands of verses added by later composers. He said that such a false propagation of Vedic Mathematics would spoil both Vedas and Mathematics done by the Indians.

20. A principal of a renowned college said that this book showed the boastful nature of the Aryan mind because they have proclaimed, "*I am the giver and source of knowledge and wisdom*." He added, "Ideology (philosophy) and Reality (accurate science) couldn't be compared or combined. Vedic Mathematics is only a very misleading concept, it is neither Vedic nor mathematics for such a combination cannot sustain. Further ideology (philosophy) varies from individual to individual depending on his or her faith, religion and living circumstances. But a reality like mathematics is the same for everybody irrespective of religion, caste, language, social status or circumstance. Magic or tricks are contradictory to reality. Vedic Mathematics is just a complete bundle of empty noise made by Hindutva to claim their superiority over others." It has no mathematics or educational value.

21. A sociologist said throughout the book they do not even say that zero and the number system belongs to Indians, but they say that it belongs to Hindus—this clearly shows their mental make-up where they do not even identify India as their land. This shows that they want to profess that Vedic Mathematics belongs to Aryans and not to the people of India.

22. We interviewed a small group of 6 scholars who were doing their doctorate in Hindu Philosophy and religion. They were given a copy of the work of Jagadguru Puri Sankaracharya for their comments and discussions about Vedic Mathematics and its authenticity as a religious product. We met them two weeks later.

We had a nice discussion over this topic for nearly three hours. The scholars showed enthusiasm over the discussions. We put only one question "Does a Jagadguru Sankaracharya of Puri need 10 years to invent or interpret the 16 sutras that too in mathematics? Is it relevant to religion? Can a religious head extend the sacred Vedas?" This was debated and their views were jotted down with care. All of them said that it was not up to the standard for the Jagadguru Sankaracharya of Puri to reflect for 10 year about Vedic Mathematics and the 16 sutras when the nation was in need of more social and ethical values. His primary duty was to spread the philosophy of Vedanta. Instead, the discovery of the sutras, his own interpretations about the use in Calculus or Algebra or Analytical geometry and so on which are topics of recent discovery puts Vedas in a degrading level. Swamiji should have reflected only upon Vedanta and not on Vedic Mathematics that is practically of no use to humanity or world peace. They added that Vedic Mathematics caters only to simple school level mathematics though tall claims have been made about its applications to other subjects. They felt the biggest weakness of the Swamiji was that he was not in a position to completely come away from Academics and become a pontiff. He was unable to come out of the fascination of working with arithmetic because he found more solace and peace with mathematics rather than Vedanta. That is why he wasted 10 years. He was able to renounce everything but was not able to renounce simple arithmetic, only this led him to write that book. They felt that Vedic Mathematics would take the student community towards materialism than towards philosophy. The only contention of these students was, "Swamiji has heavily failed to do his duty. His work on Vedic Mathematics is of no value but it is only a symbol of disgrace." They asked us to record these statements. They concluded that he was more an ordinary than an extraordinary saint or mathematician.

23. We met a retired Educational Officer of schools in Tamil Nadu. He was a post-graduate in mathematics. As soon as the Vedic Mathematics book was published, several Brahmin officials had wanted to include it in the school syllabus just as it had been included in the school curriculum in states ruled by BJP and RSS like Uttar Pradesh, Madhya Pradesh, Rajasthan and Himachal Pradesh. But in Tamil Nadu with its history of rationalism there was no possibility of introducing Vedic Mathematics into the school syllabus. He lamented that they have succeeded in unofficially teaching Vedic Mathematics in all schools run by Hindutva forces.

    He also added that young, non-Brahmin children face a lot of ill treatment and harassment in the Vedic Mathematics classroom on account of their caste. He felt that the state should intervene and ensure that is not made compulsory for children to learn Vedic Mathematics in any school. Persons who accept Vedic Mathematics will be led to believe in caste superiority, so it is just a powerful attempt to impose Hindutva. He wondered how so-called experts like Dr. Singhvi, Dr.V.S.Narasimhan, Mr.Mayilvanan, Dr.P.K.Srinivasan, S.Haridas Kadayil, S.C.Sharma (NCERT Ex-Chairman) had the heart to recommend this book with no serious mathematical content. In my opinion this misleading of Vedic Mathematics cannot penetrate in south India for we are more rationalistic than the north. They can only spread this rubbish in the north that too only as long as Hindutva forces rule these states!

24. We interviewed well-placed persons working in banks, industries and so on. Most of them said that when Vedic Mathematics was introduced they came to know about it through their children or friends. A section of them said that they were able to teach the contents of the book to their children without any difficulty because the standard was only primary school level.

They said it was recreational and fun, but there was no relevance in calling it as Vedic Mathematics. We are not able to understand why it should be called Vedic Mathematics and we see no Vedas ingrained in it. The sutras are just phrases, they seem to have no mathematical flavour. This book could have been titled "*Shortcut to Simple Arithmetical Calculations*" and nothing more. Some of them said an amateur must have written the book! Few people felt that the Swamiji would have created these phrases and called them sutras; then he would have sought some help from others and made them ghost-write for him. Whatever the reality what stands in black and white is that the material in the book is of no mathematical value or Vedic value!

25. The new fuzzy dynamical system gives results with membership degrees 0.9 or 0.8, which in fact is very high. The least degree 0 corresponds to the node "Vedic Mathematics has no Vedic content." No other node ever gets its membership degree to be too low. In almost all the cases the resultant vector gets a membership grade greater than or equal to 0.6. Thus all the nodes given by the educated under the nine categories happens to give more than 0.5 membership degree. The largest number of persons belonged to this group and everyone's views were taken to form the new fuzzy dynamical system. We took their views on the 14 attributes.

    We divided the educated into 9 groups according to their profession and the type of education. The conclusions reflected uniformity, because all the 9 categories of people held the same opinion. At no time 0 or less than 0.8 was obtained from the representatives of the educated category, which clearly shows they all held a common view, this is evident from the detailed working given in chapter 4, section 4.4.

26. We interviewed a mathematical expert who was associated for a few years with the textbook committees and advisory

member in the NCERT and who came out of it because he felt that he could not accept several of the changes made by them. He felt sad that plane geometry was given least importance and so on. We asked his opinion about the book on Vedic Mathematics. He gave a critical analysis of the views given on the back cover and asked us if we had the guts to put his views in our book. He was very critical but also down to earth. Here are some of his views: "Dr. S.C. Sharma, Ex-Head of Department of Mathematics, NCERT does not know the difference between *subjects* in mathematics and *tools* in mathematics when he wrote the sentence, "All subjects in mathematics—Multiplication, Division, Factorization etc. are dealt in 40 chapters vividly working out problems in the easiest ever method discovered so far". These operations, especially multiplication, division and factorization (of numbers) can be only categorized under arithmetical or algebraic operations and not as subjects in mathematics. It is unfortunate that an Ex-Head of the Department of Mathematics in NCERT does not know the difference between simple arithmetic operations and subjects in mathematics. Next, he feels no one needs eight long years to find these fundamental mathematical tricks because most of the school mathematics teachers who are devoted to teaching and imparting good mathematics to school students frequently discover most of these shortcut in calculations by themselves. So these techniques and many more such techniques (that are not explained in this book) are used by the good teachers of mathematics. Further, this simple arithmetical calculation cannot be called as any "magic" (which S.C. Sharma calls). Also, our methods are no match to modern western methods. We are inferior to them in mathematics too." Then, he took our permission and quoted from the editorial on *Shanghai Rankings and Indian Universities* published in Current Science, Vol. 87, No. 4 dated 25 August 2004 [7]. "The editorial is a shocking revelation about the fate of higher education and the slide down of scientific research in India. None of the reputed '5 star' Indian Universities qualifies to find a slot among the top 500 at the global level. IISC Bangalore and IITs at Delhi and Kharagpur provide some redeeming feature and put India on the scoreboard with a rank between 250 and 500. Some of the interesting features of the Shanghai Rankings are noteworthy

(i) Among the top 99 in world, we have universities from USA (58), Europe (29), Canada (14), Japan (5), Australia (2) and Israel (1)[*] ([*] numbers in brackets show the number of universities in the respective countries. India has no such university),
(ii) On the Asia-Pacific list of top 90 universities we have maximum number of universities from Japan (35), followed by China (18) including Taiwan (5) and Hong Kong (5), Australia (13), South Korea (8), Israel (6), India (3), New Zealand (3), Singapore (2) and Turkey (2)
(iii) Indian universities lag behind even small Asian countries viz. South Korea, Israel, Taiwan and Hong Kong in ranking. Sadly the real universities in India are limping, with the faculty disinterested in research outnumbering those with an academic bent of mind. The malaise is deep-rooted and needs a complete overhaul analysis of the Indian educational system.

Balaram P Curr. Sci 2004 86 (1347 – 1348) .http-ed-sjtueduin/ranking.htm says H.S. Virk. 360 Sector 71 SAS Nagar 160071 India e-mail virkhs@yahoo.com.

What is the answer to Virk's question? What do we have to boast of greatness?" He strongly feels that if Vedic Mathematics as mentioned in the book [51] was taught to students it would only make them fail to think or reason. He concluded by saying, "As a teacher I can say that thousands of students, that too from the rural areas are very bright and excel in mathematics. If all this bunkum is taught, it will certainly do more harm than any good to them. It is high time the Indian government bans the use of this book in schools from northern states.

## 5.5 Observations from the Views of the Public

The public was a heterogeneous group consisting of political party leaders, trade union leaders, activists from women organizations, social workers, NGO representatives, religious leaders especially from the Indian minority communities of Christianity and Islam and so on. Most of them did not boast of great educational qualifications. But they were in public life for over 2 decades fighting for social justice. Some of them were human right activists, some of them worked in people's movements or political parties.

Because we had no other option we had to choose a very simple mathematical tool that could be explained to these experts for mathematical purposes. Further we always had the problem of mathematical involvement. Now we give the results of mathematical analysis and the views of them as observations.

1. All of them were very against the fact that Vedic Mathematics was "magic" because when it has been claimed (that too by the ex-Head of Mathematics NCERT), our experts felt that Vedic Mathematics cannot be considered mathematics at all. According to the best of their knowledge, mathematics was a real and an accurate science. In fact it was the queen of sciences. So they all uniformly said Vedic Mathematics was not mathematics, it had some ulterior motives and aims. They also said that in several places this 'trick' must be used. They criticized it because tricks cannot be mathematics; also they condemned the use of terms like "secret of solving" because there was no need of secrets in learning mathematics. They felt that such things would unnecessarily spoil the rationalism in children. Not only would they be inhibited but also forced to think in a particular direction that would neither be productive nor inventive.

2. The experts in this category wanted to ban the book because it contains and creates more caste feelings and discrimination that are alien to the study of mathematics. It makes us clearly aware of the political agenda ingrained it. They suggested that otherwise the book should be re-titled as "Mathematical Shortcuts".

   All the16 sutras should be removed from the book because the sutras and the calculations have no significant relation. Swamiji has invented the interpretations of the Sanskrit words or phrases and calls them sutras. These do not signify any precise mathematical term or formula. When we explained in detail about the other aspects, they said that the reference to the rule of king Kamsa was unwarranted. They were quick to point out that the Swamiji had a certain criminal genius.

3. They pointed out that the major drawback in Vedic Mathematics was that it forced some sort of memorizing, only then the students could apply the shortcut methods. They argued that anything that caters to memory in elementary mathematics is only a waste and would certainly spoil the mental ability of children. In some cases children who lack such rote memory may be extremely bright as mathematics students. Several of them expressed discontent with the way the book taught the children to think.

4. A few people said that it was utter foolishness to say that these elementary calculations are found in Vedas. They argued that it not only degraded the Vedas on one side but also harmed the young children in their very formative age by making them irrational. They also extended another argument: If everything is found in the Vedas, why should children be taught anything other than the Vedas? It can very well be made the sole school curriculum. While the whole world is progressing, why should India go back to the Stone Age, they asked us rhetorically. They also questioned why the mathematical contents of the Vedas came to light only in the year 1965?

5. They were unanimous in their opinion that the popularization of Vedic Mathematics was done to implement a strong Hindutva agenda of establishing Aryan supremacy over the world. They condemned the cover of the Vedic Mathematics books in Tamil that were adorned with a picture of Saraswati, the Hindu Goddess of Learning.

    They felt that the Aryans were at an identity crisis, because they were only migrants and not natives of India. Consequently, they chose to show themselves as superior in order to subjugate the native people. That is why they kept boasting about their intellectual superiority and invented fabricated things like 'Vedic Mathematics.' They felt that Hindutva agenda was very visible because a book published in 1965, and which remained in cold storage for two decades, was dusted up and introduced in the school syllabus in the 1990s when the right-wing, Hindutva party came to power in the northern states.

6. They positively quoted the statements of Swami Vivekananda, a Sudra who emphatically said that if riots and caste clashes should not take place in South India, all Sanskrit books must be lost! He said that all the names used in them are in the northern language which is alien to the Dravidians; and that the natural differences of their culture and habits led to all these clashes. So they felt that Vedic Mathematics should not be allowed to further cause discrimination between people.

7. When we viewed the opinion of public using Neutrosophic Cognitive Maps we saw that the ON state of the only node, "Vedic Mathematics is magic has ulterior motives" made all the nodes ON except the two nodes, "It means globalization of Hindutva" and "It is a political agenda to rule the nation and if Sanskrit literature is lost it would certainly produce peace in the nation" which was in an indeterminate state.

Further this shows Vedic Mathematics is not accepted as mathematics and the text of the book that calls mathematical methods as tricks and as "magic" by the reviewers had caused suspicions among the public, especially amidst the educated and socially aware people.

8.  Further none of the dynamical system gave the value of node 2 as an indeterminate for the book always mentioned that they used tricks to solve problems. The study led to the conclusion that the popularity of Vedic Mathematics was primarily due to the capture of power in the northern states by the right wing, revivalist Hindutva forces such as BJP/RSS/VHP. The experts feel that currently the popularity of Vedic Mathematics is in the downward trend.

9.  From the combined effect of the 10 FCM matrices given by the 10 experts one sees that all the nodes come to ON state and most of them take higher degrees of membership and each and every node contributes to some degree or so.

    Most of the resultants show fixed point, except the case when node 2 is in the ON state we see the hidden pattern is a limit cycle. Further the experts were not at home to work with NCMs. Only two of them agreed to work with NCMs the rest preferred to work only with FCMs.

10. The very fact that most of the resultants gave the Hidden pattern of the FCM to be a fixed point shows the concreteness of the views do not vary!

Department of Mathematics, Indian Institute of Technology, Chennai, April 2001.

Kharagpur on Dec. 20-22, 2001, published by Narosa Publishing House, (2001) 161-168.

141. **Vasantha Kandasamy, W.B., Pathinathan, and Narayan Murthy.** Child Labour Problem using Bi-directional Associative Memories (BAM) Model, *Proc. of the 9$^{th}$ National Conf. of the Vijnana Parishad of India on Applied and Industrial Mathematics* held at Netaji Subhas Inst. of Tech. on Feb. 22-24, 2002.

142. **Vasantha Kandasamy, W.B., S. Ramathilagam and N.R. Neelakantan.** Fuzzy Optimisation Techniques in Kiln Process, *Proc. of the National Conf. on Challenges of the 21$^{st}$ century in Mathematics and its allied topics*, Feb. 3-4 (2001), Univ. of Mysore, (2001) 277-287.

143. **Vasantha Kandasamy, W.B., and Smarandache, F.**, *Fuzzy Cognitive Maps and Neutrosophic Cognitive Maps*, Xiquan, Phoenix, 2003.

144. **Vazquez, A.,** *A Balanced Differential Learning Algorithm in Fuzzy Cognitive Map*
http://monet.aber.ac.uk:8080/monet/docs/pdf_files/qr_02/qr2002alberto-vazquez.pdf

145. **Venkatbabu, Indra.** *Mathematical Approach to the Passenger Transportation Problem using Fuzzy Theory*, Ph.D. Dissertation, Guide: Dr. W. B. Vasantha Kandasamy, Department of Mathematics, Indian Institute of Technology, Chennai, June 1998.

146. **Virk, H.S.** Shanghai Rankings and Indian Universities, *Current Science*, **87**, (Aug 2004), 416.

147. **Vysoký, P.** *Fuzzy Cognitive Maps and their Applications in Medical Diagnostics*.
http://www.cbmi.cvut.cz/lab/publikace/30/Vys98_11.doc

148. **Williams, Kenneth R.** *Vedic Mathematics*, *Teachers Manual (Advanced level),* Motilal Banarsidass Publishers, Pvt. Ltd., Delhi, 2005.
213

# INDEX

## E

Edge weights of FRM, 73
Edge weights, 67
Ekādhikena Pūrvena, 11, 13-4
Ekanyūnena Pūrvena, 11, 24
Equilibrium of FRM, 74-5
Equilibrium of NRM, 90
Equilibrium state of FCM, 68
Equilibrium state of NCM, 80

## F

FCM with the feed back, 68
Fixed point of FCM, 68, 71
Fixed point of FRM, 74-5
Fixed point of NCM, 80-1
Fixed point of NRM, 90
FRM with feed back, 74-5
Fuzzy Cognitive Maps (FCMs), 65-72
Fuzzy matrices, 77
Fuzzy neutrosophic dynamical system, 92-3
Fuzzy neutrosophic matrix, 65
Fuzzy neutrosophic matrix, 92-3
Fuzzy neutrosophic multi expert system, 92-3
Fuzzy neutrosophic number, 92-3
Fuzzy neutrosophical interval, 92-3
Fuzzy nodes, 67, 73
Fuzzy Relational Maps (FRMs), 65, 72-74, 87

## G

Gunakasamuccayah, 27
Gunitasamuccayah Samuccayagunitah, 11, 26
Gunitasamuccayah, 11, 26-7

## H

Hidden pattern of FRM, 74

# ABOUT THE AUTHORS

**Dr.W.B.Vasantha Kandasamy** is an Associate Professor in the Department of Mathematics, Indian Institute of Technology Madras, Chennai. In the past decade she has guided 11 Ph.D. scholars in the different fields of non-associative algebras, algebraic coding theory, transportation theory, fuzzy groups, and applications of fuzzy theory of the problems faced in chemical industries and cement industries. Currently, four Ph.D. scholars are working under her guidance.

She has to her credit 636 research papers. She has guided over 51 M.Sc. and M.Tech. projects. She has worked in collaboration projects with the Indian Space Research Organization and with the Tamil Nadu State AIDS Control Society. This is her 29$^{th}$ book.

On India's 60th Independence Day, Dr.Vasantha was conferred the Kalpana Chawla Award for Courage and Daring Enterprise by the State Government of Tamil Nadu in recognition of her sustained fight for social justice in the Indian Institute of Technology (IIT) Madras and for her contribution to mathematics. (The award, instituted in the memory of Indian-American astronaut Kalpana Chawla who died aboard Space Shuttle Columbia). The award carried a cash prize of five lakh rupees (the highest prize-money for any Indian award) and a gold medal.
She can be contacted at vasanthakandasamy@gmail.com
You can visit her on the web at: http://mat.iitm.ac.in/~wbv or: http://www.vasantha.net

**Dr. Florentin Smarandache** is an Associate Professor of Mathematics at the University of New Mexico in USA. He published over 75 books and 100 articles and notes in mathematics, physics, philosophy, psychology, literature, rebus. In mathematics his research is in number theory, non-Euclidean geometry, synthetic geometry, algebraic structures, statistics, neutrosophic logic and set (generalizations of fuzzy logic and set respectively), neutrosophic probability (generalization of classical and imprecise probability). Also, small contributions to nuclear and particle physics, information fusion, neutrosophy (a generalization of dialectics), law of sensations and stimuli, etc.
He can be contacted at smarand@unm.edu